\documentclass[10pt, a4]{article}

\usepackage{graphicx}
\usepackage{amsmath}
\usepackage{amssymb}
\usepackage{mathptmx}
\usepackage{multicol}
\usepackage{ifpdf}
\usepackage{hyperref}
\usepackage{url}

\usepackage{amsmath}
\usepackage{amssymb}
\usepackage{graphicx}
\usepackage{color}
\usepackage{ifpdf}
\usepackage{url}
\usepackage{hyperref}

\usepackage{booktabs}
\usepackage{empheq}
\usepackage[ruled,vlined]{algorithm2e}
\usepackage{algorithmicx}
\usepackage{algpseudocode}
\usepackage{multirow}
\usepackage{enumitem}

\newtheorem{definition}{Definition}[section]
\newtheorem{lemma}{Lemma}[section]
\newtheorem{theorem}{Theorem}[section]
\newtheorem{remark}{Remark}[section]
\newtheorem{corollary}{Corollary}[section]

\newtheorem{assumption}{Assumption A.}

\newcommand{\R}{\mathbb{R}}

\newcommand{\abs}[1]{\left\vert#1\right\vert}
\newcommand{\set}[1]{\left\{#1\right\}}
\newcommand{\norm}[1]{\Vert #1\Vert}

\newcommand{\Eproof}{\hfill$\square$}

\newcommand{\xb}{\mathbf{x}}
\newcommand{\yb}{\mathbf{y}}
\newcommand{\zb}{\mathbf{z}}
\newcommand{\ub}{\mathbf{u}}
\newcommand{\vb}{\mathbf{v}}
\renewcommand{\sb}{\mathbf{s}}

\newcommand{\X}{{\bf X}}
\newcommand{\Xc}{\mathcal{X}}
\newcommand{\Yc}{\mathcal{Y}}

\newcommand{\Wc}{\mathcal{W}}
\newcommand{\Xb}{\mathbf{X}}
\newcommand{\Yb}{\mathbf{Y}}

\newcommand{\Db}{\mathbf{D}}
\newcommand{\Sb}{\mathbf{S}}

\newcommand{\dom}[1]{\mathrm{dom}\left(#1\right)}

\renewcommand{\vec}[1]{\mathrm{vec}\left(#1\right)}

\newcommand{\Lc}{\mathcal{L}}

\newcommand{\rb}{\mathbf{r}}
\newcommand{\wb}{\mathbf{w}}
\newcommand{\bb}{\mathbf{b}}
\newcommand{\Ab}{\mathbf{A}}

\title{A Primal-Dual Algorithmic Framework for Constrained Convex Minimization}

\author{\centerline{Quoc Tran-Dinh ~~$\bullet$ ~~ Volkan Cevher} \\\\
\centerline{Laboratory for Information and Inference Systems (LIONS),}\\
\centerline{\'{E}cole Polytechnique F\'{e}d\'{e}rale de Lausanne (EPFL), CH1015 - Lausanne, Switzerland.}\\
\texttt{\{quoc.trandinh,  volkan.cevher\}@epfl.ch}.
}

\begin{document}
\maketitle

\begin{abstract}
We present a primal-dual algorithmic framework to obtain approximate solutions to a prototypical constrained convex optimization problem, and  rigorously characterize how common structural assumptions affect the numerical efficiency. Our main analysis technique provides a fresh perspective on Nesterov's excessive gap technique in a structured fashion and unifies it with smoothing and primal-dual methods.  For instance, through the choices of a dual smoothing strategy and a center point, our framework subsumes  decomposition algorithms, augmented Lagrangian as well as the alternating direction method-of-multipliers methods  as its special cases, and provides optimal convergence rates on the primal objective residual as well as the primal feasibility gap of the iterates for all. 


\end{abstract}

\section{Introduction}\label{sec:intro}
This article is concerned about the following constrained convex minimization problem, which captures a surprisingly broad set of problems in various disciplines \cite{Bertsekas1989b,Chandrasekaranm2012,McCoy2014,Wainwright2014}:
\begin{equation}\label{eq:separable_convex}
  f^{\star} := \min_\xb\left\{ f(\xb): \Ab\xb = \bb, \xb\in\Xc \right \},
\end{equation}
where $f : \R^n \to \R\cup\set{+\infty}$ is a proper, closed and convex function; $\Xc \subseteq\R^n$ is a nonempty, closed and convex set; and $\Ab\in\R^{m\times n}$ and $\bb\in\R^m$ are known. 
In the sequel, we develop efficient numerical methods to approximate an optimal solution $\xb^\star$ to \eqref{eq:separable_convex} and rigorously characterize how common structural assumptions on \eqref{eq:separable_convex} affect the efficiency of the methods. 

\subsection{Scalable numerical methods for \eqref{eq:separable_convex} and their limitations}
In principle, we can obtain high accuracy solutions to \eqref{eq:separable_convex} through an equivalent unconstrained problem \cite{Boyd2004,Nocedal2006}. For instance, when $\Xc$ is absent and $f$ is smooth, we can eliminate the linear constraint $\Ab\xb = \bb$ by using a projection onto the null-space of $\Ab$ and then applying well-understood smooth minimization techniques. Whenever available, we can also exploit barrier representations of the constraints $\Xc$ and avoid non-smooth $f$ via reformulations, such as \emph{lifting}, as in the interior point method  using {disciplined convex programming} \cite{Boyd2004,Grant2004,Nemirovski2008,Nesterov2004}. While the resulting smooth and unconstrained problems are simpler than \eqref{eq:separable_convex} in theory, the numerical efficiency of the overall strategy severely  suffers  from the curse-of-dimensionality as well as the loss of the numerical structures in the original formulation.

Alternatively, we can obtain low- or medium-accuracy solutions when we augment the objective $f(\xb)$ with simple penalty functions on the constraints. For instance, we can solve
\begin{equation}\label{eq: aug}
  \min_\xb\left\{ f(\xb) + (\rho/2)\|\Ab\xb - \bb\|_2^2 ~:~ \xb\in\Xc \right \},
\end{equation}
where $\rho>0$ is a penalty parameter. 
Despite the fundamental difficulties in choosing the penalty parameter, this approach enhances our computational capabilities as well as numerical robustness since we can apply modern proximal gradient, alternating direction, and primal-dual methods. Intriguingly, the scalability of virtually all these solution algorithms rely on three key structures that stand out among many others: 

% Structure 1.
\paragraph{Structure 1 (\textbf{Decomposability}): } We say that the constrained problem \eqref{eq:separable_convex} is $p$-\textit{decomposable} if the objective function $f$ and the feasible set $\Xc$ can be represented as follows
\begin{equation}\label{eq:separable_structure}
f(\xb) := \sum_{i=1}^pf_i(\xb_i), ~~\textrm{and}~~ \Xc := \prod_{i=1}^p\Xc_i,
\end{equation}
where $\xb_i\in\R^{n_i}$, $\Xc_i\in\R^{n_i}$, $f_i : \R^{n_i}\to\R\cup\set{+\infty}$ is proper, closed and convex for $i=1,\dots, p$, and $\sum_{i=1}^pn_i = n$.
Decomposability immediately supports parallel and distributed implementations in synchronous hardware architectures. This structure arises naturally in linear programming, network optimization, multi-stages models and distributed systems \cite{Bertsekas1989b}. With decomposability, the problem \eqref{eq:separable_convex} is also referred to as a monotropic convex program \cite{Rockafellar1985}.

% Structure 2.
\paragraph{Structure 2 (\textbf{Proximal tractability}):} Unconstrained problems can still pose significant difficulties in numerical optimization when they include non-smooth terms. However, many non-smooth problems (e.g., of the form \eqref{eq: aug}) can be solved nearly as efficiently as smooth problems, provided that the computation of the proximal operator is \textbf{tractable}\footnote{It can be solved in a closed form, low computational cost or polynomial time.} \cite{Beck2009,Parikh2013,Rockafellar1976b}:  
\begin{equation}\label{eq: prox}
  \mathrm{prox}_{\lambda f}(\xb) := \mathrm{arg}\min_{\zb\in\Xc}\set{ f(\zb) + (1/(2\lambda))\Vert\zb - \xb\Vert_2^2},
\end{equation}
where $\lambda>0$ is a constant. While the proximal operators simply use $\Xc = \mathbb{R}^n$ in the canonical setting, we employ \eqref{eq: prox} to do away with the $\Xc$-feasibility of the algorithmic iterates. Many smooth and non-smooth functions  support efficient proximal operators \cite{Chandrasekaranm2012,Combettes2005,McCoy2014,Wainwright2014}. 
Clearly, decomposability proves useful in the computation of \eqref{eq: prox}.

% Structure 3.
\paragraph{Structure 3 (\textbf{Special function classes}):} Often times, the function $f$ in \eqref{eq:separable_convex} or the individual terms $f_i$ in \eqref{eq:separable_structure} possess additional properties that can enhance numerical efficiency. Table \ref{eq: structured properties} highlights common properties that are typically (but not necessarily) associated with function smoothness. These structures provide iterative algorithms with analytic upper and lower bounds on the objective (or its gradient), and aid the theoretical design of their iterations as well as their practical step-size and momentum parameter selection \cite{Beck2009,Boyd2004,Nesterov2004,Nocedal2006,Tran-Dinh2013a}.
\begin{table}[!ht]
\begin{center}
\caption{Special convex function classes. In the optimization literature, we refer to $L,\sigma$, and $\nu$ as the Lipschitz, strong convexity, and barrier parameters, respectively.}\label{eq: structured properties}
\vskip0.25cm
\begin{tabular}{c|l|c} \toprule
\textbf{Class}   &  \textbf{Name}              &  \textbf{Property} \\
                             &      &  $\xb, \yb \in \dom{f}, ~\vb \in\mathbb{R}^n, 0\leq \sigma \leq L <+\infty$ \\ \cmidrule{1-3}
$\mathcal{F}_L$    & Lipschitz gradient     &  $\norm{\nabla f(\xb) - \nabla f(\yb)}_{*} \leq L \norm{\xb-\yb} $\\ \cmidrule{1-3}
$\mathcal{F}_\sigma$  & Strong convexity   & $\frac{\sigma}{2}\norm{\xb-\yb}^2 + f(\xb) + \nabla f(\xb)^T(\yb-\xb)  \leq f(\yb)$\\ \cmidrule{1-3}
$\mathcal{F}^2$    & Standard self-concordant     &  $\abs{\varphi^{\prime\prime\prime}(t)} \leq 2 \varphi^{\prime\prime}(t)^{3/2}$: $\varphi(t) = f(\xb+t\mathbf{v}),~ t\in\mathbb{R}$\\ \cmidrule{1-3}
$\mathcal{F}^{2,\nu}$ &  Self-concordant barrier &$\mathcal{F}^2$ and  $\sup_{\vb \in\mathbb{R}^n}\set{2\nabla{f}(\xb)^T\vb - \norm{\vb}_{\xb}^2} \leq \nu$ \\ \bottomrule
\end{tabular}
\end{center}
\end{table}

On the basis of these structures, we can design algorithms featuring a full spectrum of (nearly) dimension-independent, global convergence rates for composite convex minimization problems with well-understood analytical complexities \cite{Beck2009,Nesterov2004,Nesterov2007,Nesterov2011c,Tran-Dinh2013a}. Unfortunately, the scalable, penalty-based approaches above invariably feature one or both of the following two drawbacks which blocks their full impact. %significantly reduce numerical efficiency: 

% Limitation 1.
\paragraph{Limitation 1 (\textbf{Non-ideal convergence characterizations}):} Ideally, the convergence characterization of an algorithm for solving \eqref{eq:separable_convex} must establish rates  both on absolute value of the  primal objective residual $\abs{f(\xb^k)-f^\star}$ as well as the primal feasibility of its linear constraints $\Vert\Ab\xb^k-\bb\Vert$, simultaneously on its iterates $\xb^k\in\mathcal{X}$.  The constraint feasibility  is critical so that the primal convergence rate has any significance. Rates on weighted primal objective residual and feasibility gap is not necessarily meaningful since \eqref{eq:separable_convex} is a constrained problem and $f(\xb^k)-f^\star$ can easily be negative at all times as compared to the unconstrained setting where we trivially have $f(\xb^k)-f^\star\ge 0$. 

Table \ref{tbl: convergence} demonstrates that the convergence results for some existing methods are far from ideal. 
Most algorithms have guarantees in the ergodic sense (i.e., on the averaged history of iterates without any weight) \cite{Chambolle2011,He2012a,He2012,Ouyang2014,Shefi2014,Wang2013a} with non-optimal rates, which diminishes the practical performance; they rely on special function properties to improve convergence rates on the function and feasibility \cite{Ouyang2013,Ouyang2014}, which reduces the scope of their applicability; they provide rates on dual functions \cite{Goldstein2012}, or a weighted primal residual and feasibility score \cite{Shefi2014}, which does not necessarily imply convergence on the absolute value of the primal residual or the feasibility; or they obtain convergence rate on the gap function value sequence composed both the primal and dual variables via variational inequality and gap function characterizations \cite{Chambolle2011,He2012a,He2012}, where the rate is scaled by a diameter parameter which is not necessary bounded.\footnote{We refer to the standard ADMM (see, e.g., \cite{Boyd2011}) and not the parallel ADMM variant or multi-block ADMM, which can have convergence guarantees given additional assumptions. }  
\begin{table}[ht!]
\newcommand{\cell}[1]{{\!\!\!}#1{\!\!\!}}
\begin{center}
\caption{Illustrative convergence guarantees for solving \eqref{eq:separable_convex} under the proximal tractability assumption. Note that most convergence rate results in the table are in the ergodic or \emph{averaged} sense, where $\widehat{\xb}^k = k^{-1}\sum_{i=1}^k \xb^i$.}\label{tbl: convergence}
\vskip0.25cm
%\begin{footnotesize}
\begin{tabular}{p{3.4cm}|p{3.4cm}|p{5.5cm}|c}\hline\hline
\cell{Method name} & \cell{Assumptions} &  \cell{Convergence }  & \cell{\!\! References\!\!} \\ \hline\hline
\cell{ADMM} & \cell{$\le 2$-decomposable} & \cell{ $\mathcal{O}(1/k)$ on the joint $(\xb^k,\yb^k)$ using a gap function} & \cell{\cite{Chambolle2011,He2012a,He2012}} \\\hline
\cell{Decomposition method} &  \cell{$p$-decomposable}                           & \cell{$f(\widehat{\xb}^k)-f^{\star} + r\Vert\Ab\widehat{\xb}^k-\bb\Vert_2\le \mathcal{O}(1/k)$ ($r > 0$)} & \cell{\cite{Shefi2014}} \\ \hline
\cell{[Fast] ADMM} & \cell{ $\le 2$-decomposable  and $f_1$ or $f_2\in\mathcal{F}_{\mu}$} & \cell{ $\mathcal{[O}(1/k^2)] ~\mathcal{O}(1/k)$ on the dual-objective } & \cell{\cite{Goldstein2012}} \\ \hline
\cell{Bregman ADMM} & \cell{$\le 2$-decomposable}  & \cell{$f(\widehat{\xb}^k)-f^{\star}\le \mathcal{O}(1/k) $ and $\Vert\Ab\widehat{\xb}^k-\bb\Vert_2 \le \mathcal{O}(1/\sqrt{k})$} & \cell{\cite{Wang2013a}} \\\hline
\cell{Fast Linearized ADMM} & \cell{ $\le 2$-decomposable and $f_1$ or $f_2\in\mathcal{F}_{L}$  } &  \cell{$f(\widehat{\xb}^k)-f^{\star}\le \mathcal{O}(1/k)$ and $\Vert\Ab\widehat{\xb}^k-\bb\Vert_2\le \mathcal{O}(1/k)$} & \cell{\cite{Ouyang2014}} \\  \hline
\cell{Primal-Dual Hybrid Gradient (PDHG)} & \cell{Saddle point problem} & \cell{$\mathcal{O}(1/k)$ based on gap function values composed both primal-dual variables} & \cell{\cite{Goldstein2013}} \\ \hline\hline
\cell{[Inexact] augmented Lagrangian method} & \cell{$\le 2$-decomposable} & \cell{$\abs{f(\xb^k)-f^{\star}}\le \mathcal{O}(1/k^2)$ and $\Vert\Ab\xb^k-\bb\Vert_2\le \mathcal{O}(1/k^2)$ (\emph{non-ergodic})}&  \emph{This work}\\\hline
\cell{{Decomposition methods} [Inexact] 1P2D and 2P1D} &  \cell{$p$-decomposable}     & \cell{$\abs{f(\xb^k)-f^{\star}}\le \mathcal{O}(1/k)$ and $\Vert\Ab\xb^k-\bb\Vert_2\le \mathcal{O}(1/k)$ (\emph{non-ergodic})} & \emph{This work} \\ \cmidrule{2-3}
 &  \cell{ $p$-decomposable and $f_i\in\mathcal{F}_{\sigma}$   }  & \cell{$\abs{f(\xb^k)-f^{\star}}\le \mathcal{O}(1/k^2)$,  $\Vert\Ab\xb^k-\bb\Vert_2\le \mathcal{O}(1/k^2)$, and $\Vert \xb^k-\xb^\star \Vert_2\le \mathcal{O}(1/k)$ (\emph{non-ergodic})} &  \\ \hline
\cell{{New ADMM} and its preconditioned variants} &  \cell{$\leq 2$-decomposable}     & \cell{$\abs{f(\xb^k)-f^{\star}}\le \mathcal{O}(1/k)$ and $\Vert\Ab\xb^k-\bb\Vert_2\le \mathcal{O}(1/k)$ (\emph{non-ergodic})} & \emph{This work} \\
\hline\hline
\end{tabular}
%\end{footnotesize}
\end{center}
\end{table}

% Limitation 2.
\paragraph{Limitation 2 (\textbf{Computational inflexibility}):}  Recent theoretical developments customize algorithms to exploit special function classes for scalability. We have indeed moved away from the black-box model of optimization, which forms the foundation of the interior point method's flexibility, where, for instance, we restrict ourselves to compute solely the values and the (sub)gradients of the objective and the constraints at a point. 

Unfortunately, specialized algorithms requires knowledge of function class parameters, do not address the full scope of \eqref{eq:separable_convex} (e.g., with self-concordant functions or fully non-smooth decompositions), and often have complicated algorithmic implementations with backtracking steps, which create computational bottlenecks. Moreover, these issues are further compounded by their penalty parameter selection, such as $\rho$ in \eqref{eq: aug} (cf., \cite{Boyd2011} for an extended discussion), which can significantly decrease numerical efficiency, as well as the inability to handle $p$-decomposability in an optimal fashion, which rules out parallel architectures for their computation. 

% Our contributions.
\subsection{Our contributions}
To this end, we address the following two questions in this paper: ``Is it possible to efficiently solve \eqref{eq:separable_convex} using only the proximal tractability assumption with global convergence guarantees?'' and 
``Can we actually characterize the convergence rate of the primal objective residual and primal feasibility gap separately?''
 The answer is indeed positive provided that there exists a solution in a bounded primal feasible set $\Xc$. 
 
 Surprisingly, we can still exploit favorable function classes, such as $\mathcal{F}_L$ and  $\mathcal{F}_\sigma$ when available, optimally exploit $p$-decomposability and its special $2$-decomposable sub-case, and have a penalty parameter-free black-box optimization method.  
The second question is also important since in primal-dual framework, trade-off between the primal objective residual and the primal feasibility gap is crucial, which makes algorithm numerically stable, see, e.g., \cite{Goldstein2013} for numerical examples.

To achieve the desiderata, we unify primal-dual methods \cite{Bertsekas1996d,Rockafellar1976a}, smoothing \cite{Nesterov2005c,Rockafellar1976a}, and the excessive gap function technique introduced in \cite{Nesterov2005d} in convex optimization. 

% Primal dual methods.
\paragraph{\textbf{Primal-dual methods}:} 
Primal-dual methods rely on strong duality in convex optimization \cite{Rockafellar1970} and are also related to many other methods for solving  saddle points, monotone inclusions and variational inequalities \cite{Facchinei2003}.
In our approach, we reformulate the optimality condition of \eqref{eq:separable_convex} as a mixed-variational inequality and use the  gap function as our main tool to develop the algorithms.

% Smoothing.
\paragraph{\textbf{Smoothing}:}
Smoothing techniques are widely used in optimization to replace non-smooth functions with differentiable approximations. In this work, we describe two smoothing strategies for the dual function of  \eqref{eq:separable_convex} in the Lagrange formulation based on Bregman distances and the augmented Lagrangian technique. We show that the augmented Lagrangian smoother preserves convergence properties for the algorithm to solve \eqref{eq:separable_convex} and feature a convergence rate independent of the spectral norm of $\Ab$. In addition, the Bregman smoother allows us to handle $p$-decomposability by only relying on the proximal tractability assumption.

% Excessive gap function.
\paragraph{\textbf{Excessive gap function}:}
Excessive gap technique was introduced by Nesterov in \cite{Nesterov2005d} and has been used to develop primal-dual solution methods for solving nonsmooth unconstrained problems.
In this paper, we exploit the same excessive gap idea but in a structured form for a variational inequality characterizing the optimality condition of  \eqref{eq:separable_convex}. 
We then combine these three existing techniques in order to develop a unified primal-dual framework for solving \eqref{eq:separable_convex} and analyze the convergence of its algorithmic instances under mild assumptions.

Our specific theoretical and practical contributions are as follows: 

$\mathrm{i)}$~We present a unified primal-dual framework for solving constrained convex optimization problems of the form \eqref{eq:separable_convex}. 
This framework covers augmented Lagrangian method \cite{Lan2013b,Nedelcu2014}, (preconditioned) ADMM \cite{Chambolle2011}, proximal-based decomposition \cite{Chen1994} and decomposition method \cite{TranDinh2012a} as special cases, which we make explicit in Section \ref{sec: explicit}. 

$\mathrm{ii)}$~We prove the convergence and establish rates for three variants (cf., Theorem \ref{th:convergence_A1}) of our algorithmic framework without any need to select a penalty parameter. 
An important result is the convergence rate in a non-ergodic sense of both primal objective residual $\abs{f(\bar{\xb}^k) - f^{\star}} \leq \mathcal{O}(1/k^{\alpha})$ and the primal feasibility gap $\Vert\Ab\bar{\xb}^k - \bb\Vert \leq \mathcal{O}(1/k^{\alpha})$, where $\alpha = 1$ or $2$. Our rates are considered optimal given our particular assumptions (cf., Table \ref{tbl: convergence}).

$\mathrm{iii)}$~We consider an inexact variant of our algorithmic framework for the special case of $2$-decomposability, which allows one to solve the subproblems up to given predetermined accuracy so that it still maintains the same worst-case analytical complexity as in the exact case provided that the accuracy of solving the subproblems is controlled appropriately. This variant allows us to handle $2$-decomposability with only proximal tractability assumption. 

$\mathrm{iv)}$~We show how special function classes can be exploited and describe their convergence implications. 

Our characterization is radically different from existing results such as in \cite{Beck2014,Chambolle2011,Deng2012,He2012a,He2012,Ouyang2014,Shefi2014}.
We clarify the importance of this result in Section \ref{sec:algorithm} as well as Section \ref{sec: explicit} in the context of existing convergence results for ADMM and its variants. For the $p$-decomposability, the variants corresponding to our Bregman smoothing technique can be implemented in a fully parallel and distributed manner, where the feasibility guarantee acts as a consensus rate. In special case, where $p = 2$, we propose a strategy to enhance the practical convergence rate by trading off the objective residual with the feasibility gap. 

On the computational front, we test our algorithms on several well-studied numerical problems using both synthetic and real-world data, compare them to other existing state-of-the-art methods, and provide open-source code for each application. We also discuss the update of the smoothness parameters in order to enhance the performance of the algorithms by trading-off between the optimality gap and the feasibility gap. Numerical results show the advantages of our methods on several numerical tests.

% 1.3. Related work.
\subsection{Related work} 
Due to the generality of \eqref{eq:separable_convex}, there has been an explosion of interest in the convex optimization in developing solution algorithms for it. Unfortunately, it is impossible to provide a comprehensive summary of the ever-expanding literature in any reasonable space. Hence, this subsection attempts to relate some important algorithmic frameworks for solving \eqref{eq:separable_convex} to our work with selected, representative citations in each. 

% Methods-of-multipliers/primal-dual methods
\paragraph{\textbf{Methods-of-multipliers/primal-dual methods}:}
One of the oldest primal-dual methods for solving \eqref{eq:separable_convex} is the method-of-multipliers (MoM), which is based on Lagrange dualization \cite{Bertsekas1996d}. 
Without further assumptions on $f$ and $\Xc$, the dual step of this method can be viewed as a subgradient iteration, which features a provably slow convergence rate, i.e., $\mathcal{O}(1/\sqrt{k})$, where $k$ is the iteration count. MoM is also known to be  sensitive to the step-size selection rules for damping the search direction.

In order to overcome the difficulty of nonsmoothness in the dual function, several attempts have been made. For instance, we can add either a proximal term or an augmented term to the Lagrange function of \eqref{eq:separable_convex}  to smooth the dual function \cite{Chen1994,Hamdi2005a,Hamdi2005,Necoara2008,Nedelcu2014,Rockafellar1976a}. Intriguingly, 
while the specific methods studied in \cite{Chen1994,Hamdi2005a,Hamdi2005,Rockafellar1976a} are quite borad, no global convergence rate has been established so far. 

The works in \cite{Necoara2008,Nedelcu2014}  provide convergence rates by applying Nesterov's accelerated scheme to the dual problem of \eqref{eq:separable_convex}. In recent paper \cite{Shefi2014}, the authors shows that the method proposed in \cite{Chen1994} has convergence rate $\mathcal{O}(1/k)$. However, this convergence rate is a joint between the objective residual and the primal feasibility gap, i.e., $f(\xb^k) - f^{\star} + r\Vert\Ab\xb^k - \bb\Vert_2 \leq \mathcal{O}(1/k)$ for $r  > 0$ given. We note that this convergence rate on the weighted measure does not imply the convergence rate of  $\abs{ f(\xb^k) - f^{\star}}$ and  $\Vert\Ab\xb^k - \bb\Vert_2$ separately in  constrained optimization.

In \cite{Esser2010a} the author studies several variants of the primal-dual algorithm and presented several applications in image processing. Convergence analysis of these variants are also presented in \cite{Esser2010a}, however the global convergence rate has not been provided. In \cite{Goldstein2013}, the authors describe a primal-dual hybrid gradient (PDHG) algorithm, which can be considered as a variant of the same primal-dual algorithm. 
%We note that updating parameters is a crucial point in primal-dual algorithms. 
In \cite{Goldstein2013}, the authors also studied several heuristic strategies to update the parameters, and show that the convergence rate of this algorithm is $\mathcal{O}(1/k)$ in an ergodic sense with respect to a VIP gap function values.

%% Methods from monotone inclusions and variational inequalities
\paragraph{\textbf{Methods from monotone inclusions and variational inequalities}:}
The optimality condition of \eqref{eq:separable_convex} can be viewed as a monotone inclusion or a mixed variational inequality (VIP) corresponding to both the primal and dual variables $[\xb, \yb] \in \Xc\times\R^m$. As a result, we can leverage algorithms from these two respective fields to solve \eqref{eq:separable_convex} \cite{Chambolle2011,Facchinei2003,He2012a,He2012}. 
For instance, the work in \cite{Chambolle2011} exploit the idea from variational inequality proposed in \cite{Nemirovskii2004,Nesterov2007a}. 
Splitting methods including Douglas-Rachford and predictor-corrector methods considered \cite{Combettes2005,Condat2012,Eckstein1992,He2012b,Odonoghue2012a} also belong to this direction.
However, since monotone inclusions or variational inequalities are much more general than \eqref{eq:separable_convex}, using methods tailored for optimization purposes may be more efficient in practice for solving the specific optimization problem \eqref{eq:separable_convex}. 

%% Augmented Lagrangian and alternating direction methods
\paragraph{\textbf{Augmented Lagrangian and alternating direction methods}:}
Augmented Lagrangian (AL) methods have come to offer an important computational perspective on a broad class of constrained convex problems of the form \eqref{eq:separable_convex}. In this setting, we first define the Lagrangian function associated with the linear constraint $\Ab\xb = \bb$ of \eqref{eq:separable_convex}  as $\mathcal{L}(\xb, \yb) := f(\xb) + \yb^T(\Ab\xb-\bb)$.
Then, we introduce the augmented Lagrangian function:  $\mathcal{L}_{\gamma}(\xb, \yb) := \mathcal{L}(\xb, \yb) + (\gamma/2)\norm{\Ab\xb - \bb}^2_2$ for a given penalty parameter $\gamma > 0$.
Classical augmented Lagrangian method \cite{Bertsekas1989b} solving \eqref{eq:separable_convex} produces a sequence $\set{(\xb^k, \yb^k)}_{k\geq 0}$ starting from $(\xb^0, \yb^0)\in\Xc\times\R^m$ as
\begin{equation}\label{eq:aug_Lag_alg}
\left\{\begin{array}{ll}
\xb^{k+1} &:= \mathrm{arg}\min_{\xb\in\Xc}\mathcal{L}_{\gamma}(\xb, \yb^k),\\
\yb^{k+1} &:= \yb^k + \gamma (\Ab\xb^{k+1} - \bb),
\end{array}\right.
\end{equation}
Under a suitable choice of $\gamma$, it is well-known that method \eqref{eq:aug_Lag_alg} converges to a global optimal $(\xb^{\star}, \yb^{\star})$ of \eqref{eq:separable_convex} at  $\mathcal{O}(1/k)$ rate under mild assumptions, i.e., $\mathcal{L}(\xb^k, \yb^k) - \mathcal{L}(\xb^{\star}, \yb^{\star}) \leq \mathcal{O}(1/k)$.
In fact, this method can be accelerated by applying Nesterov's accelerating scheme \cite{Nesterov2004} to obtain $\mathcal{O}(1/k^2)$ convergence rate.

Within the class of augmented Lagrangian methods, perhaps the most famous variant is the alternating direction method of multipliers (ADMM), which appears in many guises in the literature.
This method has been recognized as a special case of Douglas-Rachford splitting algorithm applying to its optimality condition \cite{Boyd2011,Eckstein1992,Goldstein2012}. 
In ADMM, given that $f$ and $\Xc$ are separable with $p = 2$. This case also covers the composite minimization problem of the form $\min_{\xb_1 \in\R^{n_1}}f_1(\xb_1) + f_2(\Ab\xb_1)$, where both $f_1$ and $f_2$ are convex. By using a slack variable, we can reformulate the composite problem into  \eqref{eq:separable_convex} as $\min_{\xb \in\R^n}f_1(\xb_1) + h(\xb_2)$ subject to $\Ab\xb_1 = \xb_2$. 
In the ADMM context,  the first problem in \eqref{eq:aug_Lag_alg} can be solved iteratively as
\begin{equation}\label{eq:alternate_step}
\left\{\begin{array}{ll}
\xb^{k+1}_1 &\!\!\!:= \displaystyle\mathrm{arg}\!\!\!\min_{\xb_1\in\Xc_1}\Big\{f_1(\xb_1) + (\yb^k)^T\Ab_1\xb_1 + (\gamma/2)\norm{\Ab_1\xb_1 - \xb_2^k  }_2^2\Big\},\\
\xb^{k+1}_2 &\!\!\!:= \displaystyle\mathrm{arg}\!\!\!\min_{\xb_2\in\Xc_2}\Big\{f_2(\xb_2) + (\yb^k)^T\Ab_2\xb_2 + (\gamma/2)\norm{\Ab_1\xb_1^{k+1} -\xb_2  }_2^2\Big\}.
\end{array}\right.
\end{equation}
The main computational difficulty of ADMM is the $\xb_1$-update problem (i.e., the first subproblem) in \eqref{eq:alternate_step}. Indeed, we have to numerically solve this step in general except when $\Ab^T\Ab$  is efficiently diagonalizable. Interestingly, the diagonalization step in many cases can be done via Fourier Transform. Many notable applications support this feature, such as matrix completion where $\Ab$ models sub-sampled matrix entries, image deblurring where $\Ab$ is a convolution operator, and total variation regularization where $\Ab$ is a differential operator with periodic boundary conditions. We can also circumvent this computational  difficulty by using a preconditioned ADMM variant \cite{Chambolle2011}. 

ADMM is one of the most popular method in practice. However, its efficiency depends significantly on the choice of the penalty parameter $\gamma$. Unfortunately, theoretical guarantee for choosing this parameter is still an open problem and is not yet well-understood. 
When $f_1$ is strongly convex, we can drop the quadratic term in the first line of \eqref{eq:alternate_step} in order to obtain an alternating minimization algorithm (AMA) \cite{Tseng1991a}. This method turns out to be a forward-backward splitting algorithm for its optimality inclusion \cite{Goldstein2012}.

\paragraph{\textbf{A note on \cite{Nesterov2005c}}:} 
We note that the approach presented in this paper builds upon the excessive gap idea in \cite{Nesterov2005c}. Technically, we use the same idea but in a much structured fashion, whereby we enforce a particular linear form in preserving the excessive gap as shown in Definition 3.2. This particular structure is key in obtaining our convergence rates.  

Moreover, since our problem setting \eqref{eq:separable_convex} is different from the general minmax formulation considered in \cite{Nesterov2005d}, there are still several differences between our algorithmic framework and the methods studied in \cite{Nesterov2005d} as a result of the excessive gap technique. 
First, we use augmented Lagrangian functions and Bregman distances for smoothing the dual problem of \eqref{eq:separable_convex}. 
Second, we consider the Lagrangian primal-dual formulation for \eqref{eq:separable_convex} where we do not have the boundedness of the feasible set of the dual variable. In this case the key estimate \cite[estimate (3.3)]{Nesterov2005c} does not apply to our setting.
Third, we update all algorithmic parameters simultaneously and do not need an odd-even switching strategy \cite[Method 1: b) and c)]{Nesterov2005d}.
Four, we do not assume that the objective function $f$ of \eqref{eq:separable_convex} has Lipschitz gradient which is required in \cite{Nesterov2005d}. 
Note that there are several important applications, where this assumption simply does not hold \cite{McCoy2014}.
Fifth, our method is applied to the constrained problem \eqref{eq:separable_convex}, which requires the feasibility gap characterization as opposed to unconstrained problems where we only need to worry about   the primal-dual optimality.

% 1.5. Paper organization.
\subsection{Paper organization}
The rest of this paper is organized as follows. In the next section, we recall basic concepts, and introduce a mixed-variational inequality formulation of \eqref{eq:separable_convex}. 
In Section \ref{sec:smoothing}, we propose two key smoothing techniques for \eqref{eq:separable_convex}, called the Bregman  and augmented Lagrangian smoothing techniques. 
We also provide a formal definition for  the excessive gap function from \cite{Nesterov2005c} and further investigate its properties. 
Section \ref{sec:algorithm} presents the main primal-dual algorithmic framework for solving \eqref{eq:separable_convex} and its convergence theory. 
Section \ref{sec:variants} specifies different instances of our algorithmic framework for \eqref{eq:separable_convex} under given assumptions. 
Section \ref{sec:existing_methods} makes further connections to existing methods in the literature. 
Section \ref{sec:impl_issues} is devoted to implementation issues and Section \ref{sec:num_results} presents numerical simulations.  
The appendix provides detail proofs of the theoretical results in the main text.

\section{Preliminaries}\label{sec:preliminaries}
First we recall the well-known definition of the Bregman distance, the primal-dual formulation for \eqref{eq:separable_convex}, and a variational inequality characterization for the optimality condition of \eqref{eq:separable_convex}, which will be used in the sequel.

%%%%%%%%%%%%%%%%%%%%%%%%%%%%%%%%%%%%
%% 2.1. Basic notation.
%%%%%%%%%%%%%%%%%%%%%%%%%%%%%%%%%%%%
\subsection{Basic notation}\label{subsec:basic_notation}
Given a proper, closed and convex function $f$, we denote $\dom{f} := \set{\xb\in\R^n~|~f(\xb) <+\infty}$ the domain of $f$, $\partial{f}(\xb) := \{\vb\in\R^n ~|~ f(\tilde{\xb}) - f(\xb) \geq \vb^T(\tilde{\xb} - \xb), ~\forall \tilde{\xb}\in\dom{f}\}$ the subdifferential of $f$ at $\xb$. If $f$ is differentiable, $\nabla{f}(\xb)$ denotes the gradient of $f$ at $\xb$.
For given vector $\xb\in\R^n$, we define $\norm{\xb}_2$ the Euclidean norm of $\xb$.
We use a superscripted notation $L^f > 0$ to denote the corresponding Lipschitz constant of a differentiable function $f$. Similarly,  we use a subscripted notation $\sigma_g > 0$ to denote the corresponding strong convexity constant of a convex function $g$. 

%%%%%%%%%%%%%%%%%%%%%%%%%%%%%%%%%%%%
%% 2.2. Proximity functions and Bregman distance.
%%%%%%%%%%%%%%%%%%%%%%%%%%%%%%%%%%%%
\subsection{Proximity functions and Bregman distances}\label{subsec:prox_bregman}
Given a nonempty, closed convex set $\Xc$, a nonnegative, continuous and $\sigma_b$-strongly convex function $b$ is called a \textit{proximity} function (or prox-function) of $\Xc$ if $\Xc\subseteq \dom{b}$. 
For example, the simplest prox-function is $b_{\Xc}(\xb) := (\sigma_b/2)\norm{\xb - \xb_c}^2_2$ for any $\sigma_b > 0$ and $\xb_c\in\Xc$. 
Whenever unspecified, we use this specific prox-function with $\sigma_b = 1$.

Given a smooth prox-function $b$ of $\Xc$ with the parameter $\sigma_b > 0$. We define 
\begin{equation}\label{eq:bregman_distance}
d_b(\xb, \yb) := b(\xb) - b(\yb) - \nabla{b}(\yb)^T(\xb - \yb), ~~\forall \xb, \yb\in\dom{b},
\end{equation}
the Bregman distance between $\xb$ and $\yb$ with respect to $b$.
Given a matrix $\Sb$, we also define the projected prox-diameter of a given set $\Xc$ with respect to $d_b$ as
\begin{equation}\label{eq:bregman_diameter}
D^{\Sb}_{\Xc} := \sup_{\xb, \xb_c\in\Xc} d_b(\Sb\xb, \Sb\xb_c).
\end{equation}
Here, we project the set $\Xc$ onto the range space of matrix $\Sb$. If $\Xc$ is bounded, then $0 \leq D^{\Sb}_{\Xc} < +\infty$.
For $b(\xb) := (\sigma_b/2)\norm{\xb - \xb_c}^2_2$, we have $d_b(\xb, \yb) = (\sigma_b/2)\norm{\xb-\yb}_2^2$, which is indeed the Euclidean distance.

%%%%%%%%%%%%%%%%%%%%%%%%%%%%%%%%%%%%
%% 2.3. Primal-dual formulation.
%%%%%%%%%%%%%%%%%%%%%%%%%%%%%%%%%%%%
\subsection{Primal-dual formulation}\label{subsec:primal_dual_formulation}
We write the min-max formulation of \eqref{eq:separable_convex} based on the Lagrange dualization as follows:
\begin{equation}\label{eq:min_max}
\max_{\yb\in\R^m}\min_{\xb\in\Xc}\Lc(\xb, \yb) \equiv \max_{\yb\in\R^m}\min_{\xb\in\Xc} \{ f(\xb) + \yb^T(\Ab\xb - \bb) \},
\end{equation}
where $\Lc$ is the Lagrange function and $\yb$ is the dual variable. We write the dual function $g(\yb)$ as
\begin{equation}\label{eq:dual_func}
g(\yb) := \min_{\xb\in\Xc}\{  f(\xb) + \yb^T(\Ab\xb - \bb) \},
\end{equation}
which leads to the following definition of the so-called dual problem
\begin{equation}\label{eq:dual_prob}
g^{\star} := \max_{\yb\in\R^m} g(\yb).
\end{equation}
%
%First, we can write the min-max problem of \eqref{eq:separable_convex} as follows:
%\begin{equation}\label{eq:min_max}
%\max_{\yb\in\R^m}\min_{\xb\in\Xc}\Lc(\xb, \yb) \equiv \max_{\yb\in\R^m}\min_{\xb\in\Xc} \{ f(\xb) + \yb^T(\Ab\xb - \bb) \}.
%\end{equation}
%Then, if we define the dual function $g(\yb)$ as
%\begin{equation}\label{eq:dual_func}
%g(\yb) := \min_{\xb\in\Xc}\{  f(\xb) + \yb^T(\Ab\xb - \bb) \},
%\end{equation}
%then the dual problem can be expressed as
%\begin{equation}\label{eq:dual_prob}
%g^{\star} := \max_{\yb\in\R^m} g(\yb).
%\end{equation}
Let $\xb^{\star}(\yb)$ be a solution of \eqref{eq:dual_func} at a given $\yb\in\R^m$.  Corresponding to $\xb^{\star}(\yb)$, we also define the domain of $g$ as 
\begin{equation}\label{eq:dom_g}
\dom{g} := \set{ \yb\in\R^m ~|~ \xb^{\star}(\yb) ~\textrm{exists}}.
\end{equation}
If $f$ is continuous on $\Xc$ and if $\Xc$ is compact, then $\xb^{\star}(\yb)$ exists for any $\yb\in\R^m$.
Unfortunately, the dual function $g$ is typically nonsmooth, and hence the numerical solutions of \eqref{eq:dual_prob} are usually difficult \cite{Nesterov2004}. 
In general, we have $g(\yb) \leq f(\xb)$, which is known as weak-duality in convex optimization. In order to guarantee strong duality, i.e., $f^{\star} = g^{\star}$ for \eqref{eq:separable_convex} and \eqref{eq:dual_prob}, we require the following assumption:

% Assumption A.1.
\begin{assumption}\label{as:A1}
The constraint set  $\mathcal{X}$ and the solution set $\mathcal{X}^{\star}$  of \eqref{eq:separable_convex} are nonempty.
The function $f$ is proper, closed and convex.
In addition, either $\mathcal{X}$ is a polytope or the following Slater condition holds:
\begin{equation}\label{eq:slater_cond}
\set{\xb \in\mathbb{R}^n ~|~ \Ab\xb = \bb} \cap \mathrm{relint}(\Xc) \neq \emptyset,
\end{equation}
where $\mathrm{relint}(\mathcal{X})$ is the relative interior of $\mathcal{X}$.  
\end{assumption}
%The solution set $\mathcal{X}^{\star}$ of \eqref{eq:separable_convex} is nonempty.
%The function $f$ is proper, closed and convex and $\mathcal{X}$ is also nonempty, closed and convex. 
%In addition, either $\mathcal{X}$ is a polytope or the following Slater condition holds:
%\begin{equation}\label{eq:slater_cond}
%\set{\xb \in\mathbb{R}^n ~|~ \Ab\xb = \bb} \cap \mathrm{relint}(\Xc) \neq \emptyset,
%\end{equation} 
%where $\mathrm{relint}(\mathcal{X})$ is the relative interior of $\mathcal{X}$.
%\end{assumption}

Under Assumption \ref{as:A1}, the solution set $\Yc^{\star}$ of the dual problem \eqref{eq:dual_prob} is also nonempty and bounded. Moreover, the strong duality holds, i.e., $f^{\star} = g^{\star}$. Any point $(\xb^{\star}, \yb^{\star})\in\Xc^{\star}\times\Yc^{\star}$ is a primal-dual solution to \eqref{eq:separable_convex} and \eqref{eq:dual_prob}, and is also a saddle point of the Lagrange function $\Lc$, i.e., $\Lc(\xb^{\star},\yb) \leq \Lc(\xb^{\star}, \yb^{\star}) \leq \Lc(\xb, \yb^{\star})$ for all $\xb\in \Xc$ and $\yb\in\R^m$. 
These inequalities lead to the following estimate
\begin{equation}\label{eq:pd_lower_bound}
f(\xb) - g(\yb) \geq f(\xb) -  f^{\star}  \geq -\norm{\yb^{\star}}_2\norm{\Ab\xb - \bb}_2, ~~\forall \xb\in\Xc, ~\yb\in\R^m.
\end{equation}
Our goal in this paper is to solve the primal constrained problem  \eqref{eq:separable_convex}, while numerical algorithms only give an approximate solution up to a certain accuracy.
Hence, we  need to specify the concept of an approximate solution for \eqref{eq:separable_convex}.
\begin{definition}\label{de:approx_sol}
Given a target accuracy $\varepsilon \geq 0$, a point $\tilde{\xb}^{\star}\in\Xc$ is said to be an $\varepsilon$-solution of \eqref{eq:separable_convex} if $ \abs{f(\tilde{\xb}^{\star}) - f^{\star}} \leq \varepsilon$ and $\Vert\Ab\tilde{\xb}^{\star} - \bb\Vert_2 \leq \varepsilon$.
\end{definition}
Here, we assume in Definition \ref{de:approx_sol} that $\tilde{\xb}^{\star}\in\Xc$, i.e., $\tilde{\xb}^{\star}$ is exactly feasible to $\Xc$. This requirement is reasonable in practice since $\Xc$ is usually a ``simple'' set where the projection onto $\Xc$ can be computed exactly. 
Moreover, we can use different accuracy levels for the absolute value of the primal objective residual $\abs{f(\tilde{\xb}^{\star}) - f^{\star}}$ and the primal feasibility gap $\Vert\Ab\tilde{\xb}^{\star} - \bb\Vert_2$ in Definition \ref{de:approx_sol}.

%%%%%%%%%%%%%%%%%%%%%%%%%%%%%%%%%%%%
%% 2.4. Variational inequality formulation and its gap function.
%%%%%%%%%%%%%%%%%%%%%%%%%%%%%%%%%%%%
\subsection{Mixed-variational inequality formulation and gap function}
Let $\wb := (\xb, \yb) \equiv (\xb^T, \yb^T)^T \in\R^n\times \R^m$ be the primal-dual variable and $F(\wb) := \begin{pmatrix}\Ab^T\yb\\ \bb - \Ab\xb\end{pmatrix}$ be a partial Karush-Kuhn-Tucker  mapping. 
Then, the optimality condition of \eqref{eq:separable_convex} becomes
\begin{equation}\label{eq:VIP}
f(\xb) - f(\xb^{\star}) + F(\wb^{\star})^T(\wb - \wb^{\star}) \geq 0, ~~ \forall \wb \in\Xc\times \R^m,
\end{equation}
which is known as a \textit{mixed-variational inequality} \cite{Facchinei2003}.
If we define 
\begin{equation}\label{eq:gap_func}
G(\wb^{\star}) := \max_{\wb\in\Wc :=\Xc\times \R^m}\set{f(\xb^{\star}) - f(\xb) + F(\wb^{\star})^T(\wb^{\star} - \wb)},
\end{equation}
then $G$ is known as the Auslender gap function of \eqref{eq:VIP} \cite{Auslender1976}.

Let $\Wc := \Xc \times \R^m$. Then, by the definition of $F$, we can see that 
\begin{equation*}
G(\wb^{\star}) = \displaystyle\max_{(\xb, \yb)\in\Wc}\set{f(\xb^{\star})  + \yb^T(\Ab\xb^{\star} - \bb) - f(\xb) - (\Ab\xb - \bb)^T\yb^{\star}} = f(\xb^{\star}) - g(\yb^{\star}) \geq 0.
\end{equation*} 
It is clear that $G(\wb^{\star}) = 0$ if and only if $\wb^{\star} := (\xb^{\star}, \yb^{\star}) \in \Wc^{\star} := \Xc^{\star}\times\Yc^{\star}$, which is indeed the \textit{strong duality} property.

\section{Primal-dual smoothing techniques}\label{sec:smoothing}
% studies the concept of smoothable functions . % for our primal-dual framework.
This section shows how to use augmented Lagrangian functions and Bregman distances as a principled smoothing technique \cite{Nesterov2004,Beck2012} within our primal-dual framework. We can then obtain different algorithmic variants by simply choosing an appropriate prox-center at each iteration.  %and allow it to vary with respect to the iterations of the algorithms.

%%%%%%%%%%%%%%%%%%%%%%%%%%%%%%%%%
%% 3.1. Smoothing function.
%%%%%%%%%%%%%%%%%%%%%%%%%%%%%%%%%
\subsection{Dual function is a smoothable function}
The dual function $g$ defined by \eqref{eq:dual_func} is convex but in general nonsmooth. 
We  approximate this function by a smoothed function $g_{\gamma}$ defined as:
\begin{equation}\label{eq:g_gamma}
g_{\gamma}(\yb) := \min_{\xb\in\Xc}\set{ f(\xb) + \yb^T(\Ab\xb - \bb) + \gamma d_b(\Sb\xb, \Sb\xb_c) },
\end{equation}
where $d_b$ is a given Bregman distance with the strong convexity parameter $\sigma_{d} > 0$, 
$\xb_c\in\Xc$ is the prox-center of $d_b$, $\Sb$ is a given consistent projection matrix and $\gamma > 0$ is a [primal] \textit{smoothness} parameter.
The following definition characterizes approximation properties of the smoothed function $g_{\gamma}$.

% Definition 3.1.
\begin{definition}[\cite{Beck2012}]\label{def:smooth_func}
The dual function $g$ defined by \eqref{eq:dual_func} is called a $(\gamma, D, \bar{L}^g)$-smoothable function if there exist positive numbers $\gamma$, $D$ and $\bar{L}^g$ and a concave and smooth function $g_{\gamma} : \dom{g} \to \R\cup\set{+\infty}$ so that:
\begin{equation}\label{eq:smoothed_gap_of_g}
g_{\gamma}(\yb) - \gamma D \leq g(\yb)  \leq g_{\gamma}(\yb), ~~\forall \yb\in\dom{g}.
\end{equation}
In addition, $\nabla{g}_{\gamma}(\cdot)$ is Lipschitz continuous with a Lipschitz constant $L^g_{\gamma} := \gamma^{-1}\bar{L}^g$.
\Eproof
\end{definition}
We call $g_{\gamma}$ the $(\gamma, D, \bar{L}^g)$-smoothed function of $g$ or simply the smoothed function of $g$ when these parameters are specified.
We note that $g_{\gamma}$ defined by \eqref{eq:g_gamma} is not necessarily Lipschitz gradient for an arbitrary choice of $\Sb$ and $\xb_c$. We consider two cases as follows.

%%%%%%%%%%%%%%%%%%%%%%%%%%%%%%%%%
%% 3.1. Smoothing via augmented Lagrangian function.
%%%%%%%%%%%%%%%%%%%%%%%%%%%%%%%%%
\subsubsection{Smoothing via augmented Lagrangian}
Let us choose $d_b(\ub,\ub_c) := (1/2)\norm{\ub-\ub_c}^2_2$, $\Sb \equiv \Ab$ and $\xb_c\in\Xc$ so that $\Ab\xb_c = \bb$. Then, we have trivially $d_b(\Sb\xb, \Sb\xb_c) :=  (1/2)\norm{\Ab\xb - \bb}_2^2$. As a result, the function $g_{\gamma}$ defined by \eqref{eq:g_gamma} becomes the augmented dual function, that is
\begin{equation}\label{eq:aug_dual_func}
\tilde{g}_{\gamma}(\yb) :=   \min_{\xb\in\Xc}\set{f(\xb) + \yb^T(\Ab{\xb} - \bb) + (\gamma/2)\norm{\Ab{\xb} - \bb}_2^2}.
\end{equation}
Here, $\mathcal{L}_{\gamma}(\xb, \yb) := f(\xb) + \yb^T(\Ab{\xb} - \bb) + (\gamma/2)\norm{\Ab{\xb} - \bb}_2^2$ is exactly the augmented Lagrangian of \eqref{eq:separable_convex} associated with the linear constraint $\Ab\xb = \bb$.
We denote by $\tilde{\xb}^{\star}_{\gamma}(\yb)$ the solution of  \eqref{eq:aug_dual_func} and  $\dom{\tilde{g}_{\gamma}} := \set{\yb\in\R^m ~|~ \tilde{\xb}^{\star}_{\gamma}(\yb)~\textrm{exists}}$.
It is well-known that $\tilde{g}_{\gamma}$ is concave as well as smooth, and its gradient is Lipschitz continuous with a Lipschitz constant $L^{\tilde{g}}_{\gamma} := \gamma^{-1}$.
We refer to $\tilde{g}_{\gamma}$ as an augmented Lagrangian smoother (in short, \textit{AL smoother}) of $g$.
The following lemma shows that $\tilde{g}_{\gamma}$ is a smoothed function of $g$, whose proof can be found, e.g., in \cite{Bertsekas1996d}.

% Lemma 2.3.
\begin{lemma}\label{le:aug_Lag_smoothing}
For any $\gamma > 0$, $\tilde{g}_{\gamma}$  defined by \eqref{eq:aug_dual_func} is concave and smooth. Its gradient is given by $\nabla\tilde{g}_{\gamma}(\yb) = \Ab\tilde{\xb}^{\star}_{\gamma}(\yb) - \bb$ and satisfies:
\begin{equation}\label{eq:main_est1_for_aug_Lag_g}
\norm{\nabla{\tilde{g}_{\gamma}}(\yb) - \nabla{\tilde{g}_{\gamma}}(\hat{\yb})}_2 \leq L^{\tilde{g}}_{\gamma}\norm{\yb - \hat{\yb}}_2, ~~\forall \yb, \hat{\yb}\in \dom{\tilde{g}_{\gamma}},
\end{equation}
where $L^{\tilde{g}}_{\gamma} := \gamma^{-1} > 0$. 

Consequently, $\tilde{g}_{\gamma}$ is a ($\gamma, D^{\Ab}_{\Xc}, \bar{L}^{\tilde{g}}$)-smoothed function of $g$ in the sense of Definition \ref{def:smooth_func}, i.e., $\tilde{g}_{\gamma}(\yb) - \gamma D_{\Xc}^{\Ab} \leq g(\yb) \leq \tilde{g}_{\gamma}(\yb)$, where $D^{\Ab}_{\Xc} := (1/2)\sup_{\xb\in\Xc}\norm{\Ab\xb - \bb}_2^2$ and $\bar{L}^{\tilde{g}} := 1$.
\end{lemma}
% End of Lemma 2.1.

%%%%%%%%%%%%%%%%%%%%%%%%%%%%%%%%%%%
%% 3.2. Smoothing via proximity funciton.
%%%%%%%%%%%%%%%%%%%%%%%%%%%%%%%%%%%
\subsubsection{Smoothing via Bregman distances}\label{subsec:smooth_prox_func}
If we choose $\Sb := \mathbb{I}$ to be the identity matrix of $\mathbb{R}^n$, then the smoothed function $g_{\gamma}$ defined by \eqref{eq:g_gamma} becomes
\begin{equation}\label{eq:smoothed_g}
\hat{g}_{\gamma}(\yb) := \min_{\xb\in\Xc} \set{f(\xb) + \yb^T(\Ab\xb - \bb) +  \gamma d_b(\xb, \xb_c)}.
\end{equation}
Let us denote by $\hat{\xb}_{\gamma}^{\star}(\yb)$ the solution of \eqref{eq:smoothed_g}, which always exists. 
We refer to $\hat{g}_{\gamma}$ as a Bregman distance smoother (shortly, \textit{BD smoother}) of $g$.
The following lemma summarizes the properties of $\hat{g}_{\gamma}$ (see, e.g., \cite{Nesterov2005c,TranDinh2012a}):

%% Lemma 2.1.
\begin{lemma}\label{le:smoothed_gap}
The function $\hat{g}_{\gamma}$ defined by \eqref{eq:smoothed_g} satisfies:
\begin{equation}\label{eq:smoothed_gap}
\hat{g}_{\gamma}(\yb) - \gamma D^{\mathbb{I}}_{\Xc} \leq \hat{g}_{\gamma}(\yb) - \gamma d_b(\xb^{\star}(\yb), \xb_c) \leq g(\yb) \leq g_{\gamma}(\yb), ~\forall \yb\in \R^m,
\end{equation}
where $D^{\mathbb{I}}_{\Xc}$ is the prox-diameter of $\Xc$ with respect to $d_b$ and $\xb^{\star}(\yb)$ is the solution of \eqref{eq:dual_func}.

Moreover, $\hat{g}_{\gamma}$  is concave and smooth. Its gradient is given by $\nabla{\hat{g}_{\gamma}}(\yb) := \Ab\hat{\xb}^{\star}_{\gamma}(\yb) - \bb$ for all $\yb\in\R^m$, and satisfies 
\begin{equation}\label{eq:main_est1_for_g}
\big\Vert\nabla{\hat{g}_{\gamma}}(\yb) - \nabla{\hat{g}_{\gamma}}(\hat{\yb})\big\Vert_2 \leq L^{\hat{g}}_{\gamma}\norm{\yb - \hat{\yb}}_2, ~~\forall \yb, \hat{\yb}\in\R^m,
\end{equation}
for $L^{\hat{g}}_{\gamma} := \frac{\norm{\Ab}^2_2}{\gamma\sigma_d}$.
Consequently, $\hat{g}_{\gamma}$ is a $(\gamma, D_{\Xc}^{\mathbb{I}}, \bar{L}^{\hat{g}})$-smoothed function of $g$, where $\bar{L}^{\hat{g}} := \frac{\norm{\Ab}_2^2}{\sigma_d}$ and $\sigma_d$ is the strong convexity parameter of $d_b$.
\end{lemma}

We note that if $\Xc$ is bounded and $f$ is continuous (or $\Xc \subset \mathrm{relint}(\dom{f})$), then $\xb^{\star}(\yb)$ always exists for any $\yb\in\R^m$. 
In this case, the prox-diameter $D_{\Xc}^{\mathbb{I}}$ of $\Xc$ is finite. Consequently, \eqref{eq:smoothed_gap} holds for all $\yb\in\R^m$.

%%%%%%%%%%%%%%%%%%%%%%%%%%%%%%%%%%%
%% 3.4. Partly primal smoothing technique.
%%%%%%%%%%%%%%%%%%%%%%%%%%%%%%%%%%%
\subsection{Smoothed gap function}\label{subsec:excessive_gap}
As we observe from the previous section,  the optimality condition of \eqref{eq:separable_convex} can be represented as a variational inequality of the form \eqref{eq:VIP}.
By using Auslender's gap function $G(\cdot)$ defined by \eqref{eq:gap_func}, we can show that $\wb^{\star}\in\Wc^{\star}$ is a primal-dual solution to \eqref{eq:separable_convex} and \eqref{eq:dual_prob}.
Since the gap function $G(\cdot)$ is generally nonsmooth, we smooth it by adding the following smoothing function:
\begin{equation}\label{eq:d_gammabeta}
d_{\gamma\beta}(\wb) \equiv d_{\gamma\beta}(\xb,\yb) := \gamma d_b(\Sb\xb,\Sb\xb_c)  + (\beta/2)\norm{\yb}_2^2,
\end{equation}
where $d_b$ is a given Bregman distance, $\Sb$ is a projection matrix and $\gamma$ and $\beta$ are two positive \textit{smoothness} parameters.

\begin{remark}\label{re:different_choice_of_dual_smoother}
For simplicity of our analysis, we use a simple quadratic prox-function $(\beta/2)\Vert\yb\Vert_2^2$ in \eqref{eq:d_gammabeta} for the dual variable $\yb$. 
However, we can replace this term by $\beta d_{b_{\yb}}(\yb, \yb_c)$, where $d_{b_\yb}$ is a given Bregman distance and $\yb_c$ is a given point in $\R^m$.
However, depending on the choice of $d_{b_{\yb}}$, the dual variable $\yb^{*}_{\beta}(\cdot)$ may no longer have a closed form expression. However, the overall practical performance may be improved. 
\end{remark}

The smoothed gap function for $G$ is then defined as follows:
\begin{equation}\label{eq:smoothed_gap_func}
G_{\gamma\beta}(\bar{\wb}) := \max_{\wb\in \Xc\times\R^m}\set{ f(\bar{\xb}) - f(\xb) + F(\bar{\wb})^T(\bar{\wb} - \wb) - d_{\gamma\beta}(\wb) },
\end{equation}
where $F$ is defined in \eqref{eq:VIP}. The function $G_{\gamma\beta}$ can be considered as Fukushima's gap function \cite{Fukushima1992a} for the variational inequality problem \eqref{eq:VIP}.
We can see that $G_{\gamma\beta}(\bar{\wb}) \to G_{00}(\bar{\wb}) \equiv G(\bar{\wb})$ as $\gamma$ and $\beta \to 0^{+}$ simultaneously.

It is clear that the maximization problem \eqref{eq:smoothed_gap_func} is a convex optimization problem. We denote by $\wb^{\star}_{\gamma\beta}(\bar{\wb}) := (\xb^{\star}_{\gamma}(\bar{\yb}), \yb^{\star}_{\beta}(\bar{\xb}))$ the solution of this problem. 
Then, by using the optimality condition of \eqref{eq:smoothed_gap_func} we can easily check that $\xb^{\star}_{\gamma}(\bar{\yb})$ is the optimal solution to \eqref{eq:g_gamma} at $\yb := \bar{\yb}$, while $\yb^{\star}_{\beta}(\bar{\xb})$ can be computed explicitly as 
\begin{equation}\label{eq:y_beta}
\yb^{\star}_{\beta}(\bar{\xb}) := \beta^{-1}(\Ab\bar{\xb} - \bb).
\end{equation}
Our goal is to generate two sequences $\set{\bar{\wb}^k}_{k\geq 0}\subseteq \Wc$ and $\set{(\gamma_k,\beta_k)}_{k\geq 0}\in \R^2_{++}$ so that $\set{G_{\gamma_k\beta_k}(\bar{\wb}^k)}_{k\geq 0}$ becomes firmly contractive. We formally encode this idea using the following definition.

% Definition 2.1.
\begin{definition}[Model-based Excessive Gap]\label{de:decreasing_sequence}
Given $\bar{\wb}^k := (\bar{\xb}^k, \bar{\yb}^k) \in\Wc$ and $(\gamma_k,\beta_k) > 0$, a new point $\bar{\wb}^{k+1} := (\bar{\xb}^{k+1}, \bar{\yb}^{k+1}) \in \Wc$ and $(\gamma_{k+1}, \beta_{k+1}) > 0$ so that $\gamma_{k+1}\beta_{k+1} < \gamma_k\beta_k$ is said to be \textit{firmly contractive} $($w.r.t. $G_{\gamma\beta}$ defined by \eqref{eq:smoothed_gap_func}$)$ if:
\begin{equation}\label{eq:descent_gap}
G_{k+1}(\bar{\wb}^{k+1}) \leq (1-\tau_k)G_{k}(\bar{\wb}^{k}) - \psi_{k},
\end{equation}
where $G_k(\cdot) := G_{\gamma_k\beta_k}(\cdot)$, $\tau_k\in [0, 1)$ and $\psi_k \in \R$ are two given parameters. 
\Eproof
\end{definition}
Here, the parameter $\tau_k$ and the decay term $\psi_k$ will be specified accordantly with different algorithmic schemes. 

In the context of excessive gap technique introduced by Nesterov, the smoothed gap function $G_{\mu_1\mu_2}(\bar{\wb})$ measures the excessive gap $f_{\mu_2}(\bar{\xb}) - \phi_{\mu_1}(\bar{\yb})$  in \cite[cf., (2.5) and (2.9)]{Nesterov2005d}). 
Hence,  we will call  $G_{\gamma\beta}(\bar{\wb})$ Nesterov's smoothed gap function customized for the constrained convex problem \eqref{eq:separable_convex}.
We note that the excessive gap condition $f_{\mu_2}(\bar{\xb}) \leq \phi_{\mu_1}(\bar{\yb})$ in \cite[(3.2)]{Nesterov2005d} only requires $G_{\mu_1\mu_2}(\bar{\wb}) \leq 0$.
In our case, we structure this condition using the basic model in \eqref{eq:descent_gap} so that we can manipulate $\tau_k$ and the new parameter $\psi_k$ simultaneously to analyze the convergence of our algorithms.

In the sequel, we often assume that the second parameter $\psi_k$ is nonnegative, which allows us to estimate the convergence rate of $\set{ G_k(\bar{\wb}^k)}_{k\geq 0}$. 
However, the following remark shows that  the sequence $\set{ G_k(\bar{\wb}^k)}_{k\geq 0}$ can still converge to $0^{+}$ even if $\psi_k$ is positive. However, we find the ensuing convergence analysis to be difficult.
% Remark 1.
\begin{remark}\label{re:gap_decreasing}
Let $\set{\tau_k}_{k\geq 0} \subseteq (0, 1)$ and $\set{\psi_k}_{k\geq 0}$ be sequences in Definition \ref{de:decreasing_sequence}.  If
\begin{equation}\label{eq:condition}
\lim_{k\to\infty}\tau_k = 0, ~~\sum_{k=0}^{\infty}\tau_k =+\infty, ~~\textrm{and}~ \sum_{k=0}^{\infty}\psi_k < +\infty,
\end{equation}
then the sequence $\set{ G_k(\bar{\wb}^k)}_{k\geq 0}$ converges to $0^{+}$. 
\Eproof
\end{remark}

From Definition \ref{de:decreasing_sequence},  if $\set{\bar{\wb}^k}_{k\geq 0}\subseteq \Wc$ and $\set{(\gamma_k,\beta_k)}_{k\geq 0}\in \R^2_{++}$ satisfy the condition \eqref{eq:descent_gap}, then we have $G_k(\bar{\wb}^k) \leq \omega_kG_0(\bar{\wb}^0) - \Psi_k$ by induction, where 
\begin{equation}\label{eq:omega_Psi}
\omega_k := \prod_{j=0}^{k-1}(1-\tau_j)  ~~(k\geq 1)~~\textrm{and}~~\Psi_k :=  \psi_k + \sum_{j=0}^{k}\prod_{i=j}^{k}(1-\tau_i)\psi_{j-1}~~(k\geq 0). 
\end{equation}
Consequently, the rate of convergence of $\set{G_k(\bar{\wb}^k)}_{k\geq 0}$ depends on the rate of $\set{\tau_k}_{k\geq 0}$ and $\set{\psi_k}_{k\geq 0}$.

The next lemma shows the relation between problem \eqref{eq:separable_convex} and its smoothed function $g_{\gamma}$ and $g$.
The proof of this lemma can be found in the appendix.

% Lemma 3.6.
\begin{lemma}\label{le:excessive_gap_aug_Lag_func}
Let $g_{\gamma}$ be defined by \eqref{eq:g_gamma} and $G_{\gamma\beta}$ defined by \eqref{eq:smoothed_gap_func}. Also, let $\set{\bar{\wb}^k}_{k\geq 0} \subset \Wc$ and $\set{(\gamma_k,\beta_k)}_{k\geq 0}\in\R^2_{++}$ be the sequences satisfying Definition \ref{de:decreasing_sequence}. Then we have
\begin{equation}\label{eq:main_bound_gap}
f(\bar{\xb}^k) - g_{\gamma_k}(\bar{\yb}^k) \leq \omega_kG_{0}(\bar{\wb}^0) - \Psi_k - (1/(2\beta_k))\norm{\Ab\bar{\xb}^k - \bb}_2^2.
\end{equation}
In addition, we also have the following bound:
\begin{align}
- \norm{\yb^{\star}}_2\Vert \Ab\bar{\xb}^k - \bb\Vert_2 \leq f(\bar{\xb}^k) - f^{\star} &\leq f(\bar{\xb}^k) - g(\bar{\yb}^k)  \leq S_k\label{eq:main_bound_gap_2a},\\
\Vert\Ab\bar{\xb}^k - \bb\Vert_2 &\leq \beta_k\Big[\norm{\yb^{\star}}_2 + \sqrt{\norm{\yb^{\star}}_2^2 + 2\beta_k^{-1}S_k}\Big] \label{eq:main_bound_gap_2b}
\end{align}
where $S_k := \omega_kG_0(\bar{\wb}^0) + \gamma_kD_{\Xc}^{\Sb} - \Psi_k$, provided that $\beta_k\norm{\yb^{\star}}_2^2 + 2S_k \geq 0$.
\end{lemma}

From Lemma \ref{le:excessive_gap_aug_Lag_func} we can see that if $G_0(\bar{\wb}^0) \leq \Psi_k$, then  the primal objective residual $\abs{f(\bar{\xb}^k) - f^{\star}}$ and the primal feasibility gap $\Vert\Ab\bar{\xb}^k - \bb\Vert_2$ of \eqref{eq:separable_convex} are bounded by
\begin{equation}\label{eq:main_bound_gap_3}
\left\{\begin{array}{ll}
&\abs{f(\bar{\xb}^k) - f^{\star}} \leq \max\Big\{\gamma_kD_{\Xc}^{\Sb}, \big[2\beta_kD_{{\Yc}^{\star}} + \sqrt{2\gamma_k\beta_kD_{\Xc}^{\Sb}}\big]D_{{\Yc}^{\star}}\Big\},\\
&\Vert\Ab\bar{\xb}^k - \bb\Vert_2 \leq 2\beta_kD_{{\Yc}^{\star}} + \sqrt{2\gamma_k\beta_kD_{\Xc}^{\Sb}},
\end{array}\right.
\end{equation}
where $D_{{\Yc}^{\star}} := \min\set{ \norm{\yb^{\star}}_2 ~|~ \yb^{\star}\in\Yc^{\star}}$, which is the norm of a minimum norm dual solution.
The estimate \eqref{eq:main_bound_gap_3} hints that we can derive algorithms based on $\set{(\gamma_k,\beta_k)}$ whose convergence rate depends directly on how we update the sequence $\set{(\gamma_k,\beta_k)}_{k\geq 0}$.
\section{The main algorithmic framework}\label{sec:algorithm}
The key objective in this section is to design a primal-dual update template from $\bar{\wb}^k \in\Wc$ and $(\gamma_k, \beta_k)\in\R_{++}^2$ to $\bar{\wb}^{k+1} \in\Wc$ and $(\gamma_{k+1}, \beta_{k+1})\in\R_{++}^2$ so that the conditions in Definition \ref{de:decreasing_sequence} hold.
We develop two distinct schemes to update $\bar{\wb}^k$ and $(\gamma_k, \beta_k)$ in the following two subsections.

%% 4.1.1. The iterative scheme with two primal steps
\subsection{An iteration scheme with two primal steps}
Since the objective function is not necessary smooth, we consider the following mapping under Assumption \ref{as:A1}:
\begin{equation}\label{eq:prox_oper2}
\mathrm{prox}_{\Sb f}(\hat{\xb}, \hat{\yb}; \beta) := \mathrm{arg}\!\min_{\xb\in\Xc}\set{ f(\xb) + \hat{\yb}^T\Ab(\xb - \hat{\xb}) + (\bar{L}^g/(2\beta))\norm{\Sb(\xb - \hat{\xb})}_2^2},
\end{equation}
where $\beta > 0$ and $\Sb$ is a projection matrix that satisfies the following condition:
\begin{equation}\label{eq:condition_dp}
\norm{\Ab\xb - \bb}_2^2 \leq \norm{\Ab\hat{\xb} - \bb}_2^2 + 2(\Ab\hat{\xb}-\bb)^T\Ab(\xb -\hat{\xb}) + \bar{L}^g\norm{\Sb(\xb -\hat{\xb})}_2^2, ~~\forall \xb,\hat{\xb}\in\Xc.
\end{equation}
An obvious choice of $\Sb$ is either $\Sb \equiv \Ab$ and $\bar{L}^g = 1$ or $\Sb \equiv \mathbb{I}$ and $\bar{L}^g = \norm{\Ab}_2^2$. Since $\Ab$ is known, both are feasible. Alternatively, local variable metrics can be used here, which might lead to different adaptation and computation tradeoffs in optimization. 

Now, given $\bar{\wb}^k := (\bar{\xb}^k,\bar{\yb}^k) \in\Wc$ and $(\gamma_k,\beta_k)\in\R^2_{++}$, we compute $\xb^{\star}_{\gamma}(\bar{\yb}^k)$ the solution of the minimization problem in \eqref{eq:g_gamma} and $\yb^{\star}_{\beta}(\bar{\xb}^k)$ by \eqref{eq:y_beta}. 
Then, we update the point $\bar{\wb}^{k+1} := (\bar{\xb}^{k+1}, \bar{\yb}^{k+1})$ and $(\gamma_{k+1}, \beta_{k+1})$ based on the following scheme:
\begin{equation}\label{eq:pd_scheme_2p}
\left\{\begin{array}{ll}
\hat{\xb}^k &:= (1-\tau_k)\bar{\xb}^k + \tau_k\xb^{\star}_{\gamma_k}(\bar{\yb}^k),\\
\hat{\yb}^k &:= \beta_{k+1}^{-1}(\Ab\hat{\xb}^k - \bb),\\
\bar{\xb}^{k+1} &:= \mathrm{prox}_{\Sb f}(\hat{\xb}^k, \hat{\yb}^k; \beta_{k+1}),\\
\bar{\yb}^{k+1} &:= (1-\tau_k)\bar{\yb}^k + \tau_k\hat{\yb}^k,
\end{array}\right.\tag{\textrm{2P1D}}
\end{equation}
where $\tau_k \in (0, 1]$ and $(\beta_{k+1}, \gamma_{k+1})$ is updated as
\begin{equation}\label{eq:update_param1}
\beta_{k+1}  = (1-\tau_k)\beta_k ~~\textrm{and} ~~ \gamma_{k+1} = (1- c_k\tau_k)\gamma_k,
\end{equation}
for some $c_k \in (-1, 1]$, which will be specified later.
It is important to note that if $f$ is nonsmooth, solving problem \eqref{eq:prox_oper2} requires the same cost as solving \eqref{eq:g_gamma}. Therefore, we can refer to \eqref{eq:pd_scheme_2p} as a \textit{primal-dual scheme with two primal steps}. 

% Remark 4.1.
\begin{remark}
If $f$ is $L_f$-Lipschitz gradient, then we can replace $f(\xb)$ in the proximal step at the third line of \eqref{eq:pd_scheme_2p} by its linearization, which leads to the following gradient step:
\begin{equation*}
\mathrm{grad}_{\Sb f}(\hat{\xb}, \hat{\yb}; \beta) \!:=\! \mathrm{arg}\!\!\min_{\xb\in\Xc}\!\!\!\set{\! (\nabla{f}(\hat{\xb}) + \Ab^T\yb)^T\!\!(\xb \!-\! \hat{\xb}) \!+\! (L_f/2)\Vert\xb \!-\! \hat{\xb}\Vert_2^2 \!+\! (2\beta)^{-1}\!\norm{\Sb(\xb \!-\! \hat{\xb})}_2^2\!}.
\end{equation*}
In particular, when $f$ is $p$-decomposable as in \eqref{eq:separable_structure} and if $f_i$ is Lipschitz gradient for some $i=1,\dots, p$, then we can use the gradient step for such a $f_i$ \cite{TranDinh2012a}.
\end{remark}

The following lemma provides conditions such that $(\bar{\xb}^{k+1}, \bar{\yb}^{k+1})$ updated by \eqref{eq:pd_scheme_2p} satisfies Definition \ref{de:decreasing_sequence}, whose proof  is deferred to the appendix.

%% Lemma 4.1.
\begin{lemma}\label{le:maintain_excessive_gap1}
Let $(\bar{\xb}^{k+1}, \bar{\yb}^{k+1})$ and $(\gamma_{k+1}, \beta_{k+1})$ be updated as \eqref{eq:pd_scheme_2p} and \eqref{eq:update_param1}. If $\Sb$ satisfies \eqref{eq:condition_dp} and $\tau_k$ is chosen such that
\begin{equation}\label{eq:pd_condition1}
\beta_{k+1}\gamma_{k+1} \geq \bar{L}^g\tau^2_k,
\end{equation}
then  $(\bar{\xb}^{k+1}, \bar{\yb}^{k+1})\in \Wc$ and satisfies Definition \ref{de:decreasing_sequence}, i.e., $G_{k+1}(\bar{\wb}^{k+1}) \leq (1-\tau_k)G_{k}(\bar{\wb}^{k}) - \psi_{k}$ for  $\psi_{k} := \frac{\tau_k^2}{2\beta_{k+1}}\norm{\Ab\xb^{\star}_{\gamma_k}(\bar{\yb}^k) - \bb}_2^2 \geq 0$.
\end{lemma}

%% 4.1.2. The iterative scheme with one primal step.
\subsection{An iteration scheme with two dual steps}
Alternatively to \eqref{eq:pd_scheme_2p}, we can switch from two primal steps to two dual steps.  
In this case, the new point $(\bar{\xb}^{k+1}, \bar{\yb}^{k+1})$ is updated as follows:
\begin{equation}\label{eq:pd_scheme_2d}
\left\{\begin{array}{ll}
\hat{\yb}^k        &  := (1-\tau_k)\bar{\yb}^k + \tau_k\yb^{\star}_{\beta_k}(\bar{\xb}^k),\\
\bar{\xb}^{k+1} & := (1 \!-\! \tau_k)\bar{\xb}^k + \tau_k\xb^{\star}_{\gamma_{k+1}}(\hat{\yb}^k),\\
\bar{\yb}^{k + 1} & := \hat{\yb}^k + \frac{\gamma_{k+1}}{\bar{L}^g}\big(\Ab\xb^{\star}_{\gamma_{k+1}}(\hat{\yb}^k) - \bb\big),
\end{array}\right.\tag{\textrm{1P2D}}
\end{equation}
where $\tau_k \in (0, 1)$ and the parameters $\beta_{k+1}$ and $\gamma_{k+1}$ are updated as \eqref{eq:update_param1}.
We refer to \eqref{eq:pd_scheme_2d} as a \textit{primal-dual scheme with two dual steps}.

The following lemma shows that $(\bar{\xb}^{k+1}, \bar{\yb}^{k+1})$ updated by \eqref{eq:pd_scheme_2d} maintains \eqref{eq:descent_gap}, whose proof can also be found in the appendix.

% Lemma 4.2.
\begin{lemma}\label{le:maintain_excessive_gap2}
Let $(\bar{\xb}^{k+1}, \bar{\yb}^{k+1})$ and $(\gamma_{k+1}, \beta_{k+1})$ be updated by \eqref{eq:pd_scheme_2d} and \eqref{eq:update_param1}, respectively.
If $\tau_k$ is chosen such that
\begin{equation}\label{eq:pd_condition2}
\beta_{k+1}\gamma_{k+1} \geq \bar{L}^g\tau^2_k,
\end{equation}
then  $(\bar{\xb}^{k+1}, \bar{\yb}^{k+1})\in \Wc$ and satisfies $G_{k+1}(\bar{\wb}^{k+1}) \leq (1-\tau_k)G_{k}(\bar{\wb}^{k}) - \psi_{k}$ for  
\begin{equation*}
\psi_k :=  \tau_k(1-\tau_k)\gamma_{k}\big[d_b(\Sb\xb^{\star}_{\gamma_{k+1}}(\hat{\yb}^k), \Sb\xb_c) - c_kd_b(\Sb\xb^{\star}_{\gamma_{k\!+\!1}}(\bar{\yb}^k),\Sb\xb_c)\big]\ge 0.
\end{equation*}
\end{lemma}

%% 4.3. Find a starting point.
\subsection{Finding a starting point}
In principle, we can start our algorithm at any point $(\bar{\xb}^0, \bar{\yb}^0)\in\Wc$.
However, we can find a point $\bar{\wb}^0 := (\bar{\xb}^0, \bar{\yb}^0)\in\Wc$ such that $G_{\gamma_0\beta_0}(\bar{\wb}^0) \leq 0$.
The following lemma shows how to compute such a point, whose proof can be found in the appendix.

% Lemma 2.1.
\begin{lemma}\label{le:starting_point}
Given $\xb^0_c\in\Xc$, the point  $\bar{\wb}^0 := (\bar{\xb}^0, \bar{\yb}^0)\in\Wc$ computed by
\begin{equation}\label{eq:starting_point1}
\left\{\begin{array}{ll}
&\bar{\xb}^0  = \xb^{\star}_{\gamma_0}(0^m), \\
&\bar{\yb}^0 := \beta_0^{-1}\left(\Ab\bar{\xb}^0 - \bb\right).
\end{array}\right.
\end{equation}
satisfies $G_{\gamma_0\beta_0}(\bar{\wb}^0) \leq -\gamma_0d_b(\Sb\bar{\xb}^0, \Sb\xb_c) \leq 0$ provided that $\beta_0\gamma_0 \geq \bar{L}^g$. 

Alternatively, the point  $\bar{\wb}^0 := (\bar{\xb}^0, \bar{\yb}^0) \in \Wc$ generated by 
\begin{equation}\label{eq:starting_point2}
\left\{\begin{array}{ll}
\bar{\yb}^0 &:=  \beta_0^{-1}\left(\Ab\xb_c - \bb\right),\\
\bar{\xb}^0  & := \mathrm{prox}_{\Sb f}(\xb_c, \bar{\yb}^0; \beta_0),
\end{array}\right.
\end{equation}
also satisfies $G_{\gamma_0\beta_0}(\bar{\wb}^0) \leq -\gamma_0d_b(\Sb\bar{\xb}^0, \Sb\xb_c) \leq 0$ provided that $\beta_0\gamma_0 \geq \bar{L}^g$. 
\end{lemma}
% End of the lemma.

%% 4.4. Update the parameters.
\subsection{Updating step-size parameter}
It remains to derive an update rule for the step-size $\tau_k$ in both scheme \eqref{eq:pd_scheme_2p} and \eqref{eq:pd_scheme_2d}.
The update rule is derived by using the same condition in both Lemma \ref{le:maintain_excessive_gap1} and Lemma \ref{le:maintain_excessive_gap2}.

Since $\tau_k$ satisfies $\beta_{k+1}\gamma_{k+1} \geq \bar{L}^{g}\tau_k^2$, $\tau_{k+1}$ also satisfies the same condition, i.e., $\beta_{k+2}\gamma_{k+2}\geq \bar{L}^g\tau_{k+1}^2$.
In addition, by \eqref{eq:update_param1}, we have $\beta_{k+2} := (1 - \tau_{k+1})\beta_{k+1}$ and $\gamma_{k+2} := (1 - c_{k+1}\tau_{k+1})\gamma_{k+1}$.
These conditions lead to $\tau_{k+1}^2 \leq (1 - \tau_{k+1})(1 - c_{k+1}\tau_{k+1})\tau_k^2$. Since we want to maximize the value of $\tau_{k+1}$, we take the equality, i.e., $\tau_{k+1}^2 = (1 - \tau_{k+1})(1 - c_{k+1}\tau_{k+1})\tau_k^2$.
The last condition leads to
\begin{equation}\label{eq:tau_update_rule}
a_{k+1} := \big(1 + c_{k+1} + \sqrt{4a_k^2 + (1-c_{k+1})^2}\big)/2, ~~\textrm{and}~~\tau_k := a_k^{-1}.
\end{equation}
In addition, from Lemma \ref{le:starting_point}, we have $\beta_0\gamma_0 \geq \bar{L}^g$. Let us choose $\beta_0 := \gamma_0^{-1}\bar{L}_g$. We need to choose $\tau_0 \in (0, 1]$ such that $\gamma_1\beta_1 = (1-\tau_0)(1-c_0\tau_0)\beta_0\gamma_0 \geq \bar{L}^g\tau_0^2$. 
Therefore, we get
\begin{equation}\label{eq:tau_update_rule0}
a_0 := \Big(1 + c_0 + \sqrt{4(1-c_0)+(1+c_0)^2}\Big)/2, ~\textrm{and}~\tau_0 := a_0^{-1}.
\end{equation}
The following Lemma shows the convergence rate of $a_k$, $\beta_k$ and $\beta_k\gamma_k$. 
The proof of this lemma can be found in the appendix.

% Lemma 3.2.
\begin{lemma}\label{le:update}
Let $s_k := \sum_{i=1}^k c_i$.  Then, the sequence $\set{a_k}$ updated by \eqref{eq:tau_update_rule} with $a_0$ given by \eqref{eq:tau_update_rule0} satisfies
\begin{equation}\label{eq:tau_rate}
(k + a_0 + s_k)/2 \leq a_k \leq  k + a_0.
\end{equation}
Consequently, the sequences $\set{\beta_k}$ and $\set{\gamma_k}$ updated by \eqref{eq:update_param1} satisfy
\begin{equation}\label{eq:beta_gamma_rate}
\frac{\bar{L}^g}{(k + a_0)^2} \leq \gamma_{k+1}\beta_{k+1}  \leq \frac{4\bar{L}^g}{(k + a_0 + s_{k})^2},
\end{equation}
where $\bar{L}^g$ is given in Definition \ref{def:smooth_func}.
Moreover, we also have
\begin{equation}\label{eq:beta_rate1}
\left\{\begin{array}{lllll}
&\frac{\beta_0}{(k+2)^2} \leq &\beta_{k+1}  &\leq \frac{4\beta_0}{(k+1)^2},  ~~&\textrm{if}~c_k = 0,\\
& &\beta_{k+1}  &= \frac{\beta_0}{k+2},  ~~&\textrm{if}~c_k = 1.
\end{array}\right.
\end{equation}
\end{lemma}

\subsection{A primal-dual algorithmic template}\label{subsec:algorithm}
Now, we combine all ingredients presented in the previous subsection to obtain the template for solving \eqref{eq:separable_convex} shown in Algorithm \ref{alg:pd_alg}.

%%%%%%%%%%%%%%%%%%%%%%%%%%%%%%%%%%%%%%%%%%%%%%%%
%+ Algorithm 4.1.
%%%%%%%%%%%%%%%%%%%%%%%%%%%%%%%%%%%%%%%%%%%%%%%%
\begin{algorithm}[!ht]\caption{(\textit{Primal-dual template using model-based excessive gap technique})}\label{alg:pd_alg}
\begin{algorithmic}[1]
\Statex{\hskip-6ex}\textbf{Inputs:}  $\gamma_0 > 0$, $c_0 \in (-1, 1]$, and a smoother (AL or BD). 
\Statex{\hskip-6ex}\textbf{Initialization:} 
\State $a_0 := \big(1 + c_0 + [4(1-c_0) + (1+c_0)^2]^{1/2}\big)/2$ and $\tau_0 := a_0^{-1}$.
\State Use $\bar{L}^g := 1$ for  AL smoother and $\bar{L}^g := \sigma_d^{-1}\norm{\Ab}_2^2$ for BD smoother.
\State $\beta_0 := \bar{L}^g/\gamma_0$.
\State  Compute $(\bar{\xb}^0, \bar{\yb}^0)$ by either \eqref{eq:starting_point1} or \eqref{eq:starting_point2}.
\Statex{\hskip-6ex}\textbf{For}~{$k=0$ {\bfseries to} $k_{\max}$}
\State If \textbf{stopping\_criterion}, terminate.
\State Given $(\bar{\xb}^k, \bar{\yb}^k)$, update $(\bar{\xb}^{k+1}, \bar{\yb}^{k+1})$ by either \eqref{eq:pd_scheme_2p} or \eqref{eq:pd_scheme_2d}.
\State $\beta_{k+1} := (1-\tau_k)\beta_k$ and update $\gamma_{k+1} := (1 - c_k\tau_k)\gamma_k$.
\State Update $c_{k+1}$ from $c_k$ if necessary.
\State Update $a_{k+1} := \big(1 + c_{k+1} + [4a_k^2 + (1-c_{k+1})^2]^{1/2}\big)/2$ and set $\tau_{k+1} := a_{k+1}^{-1}$.
\Statex{\hskip-6ex}\textbf{End For}
\end{algorithmic}
\end{algorithm}
%\vskip-0.5cm
%%%%%%%%%%%%%%%%%%%%%%%%%%%%%%%%%%%%%%%%%%%%%%%%

The main step of Algorithm \ref{alg:pd_alg} is Step 5, where we need to update $(\bar{\xb}^{k+1}, \bar{\yb}^{k+1})$ based on either \eqref{eq:pd_scheme_2p} or  \eqref{eq:pd_scheme_2d}.
If we use  \eqref{eq:pd_scheme_2p}, then $\gamma_{k}$ can be updated as $\gamma_{k+1} := (1-\tau_k)\gamma_k$, i.e., $c_k = 1$. 
We can also fix $\gamma_k = \gamma_0 > 0$ for all the iterations $k \geq 0$, i.e., $c_k = 0$.
It is important to note that Step 5 and Step 6 are mixed. Depending on the use of either \eqref{eq:pd_scheme_2p} or  \eqref{eq:pd_scheme_2d}, the corresponding parameter $\beta_k$ or $\gamma_k$ is updated before Step 5.
If we choose $c_k < 0$, then $\set{\gamma_k}$ is increasing. Since the rate of $\beta_k\gamma_k$ is fixed at $\mathcal{O}(1/k^2)$ due to \eqref{eq:beta_gamma_rate}, if we decrease the rate of $\set{\gamma_k}$ (i.e., increase $\gamma_k$), then $\set{\beta_k}$ converges faster than the $\mathcal{O}(1/k^2)$ rate.
We will discuss the stopping condition at Step 4 later.
We note that we can also alternate between \eqref{eq:pd_scheme_2p} and \eqref{eq:pd_scheme_2d} in Algorithm \ref{alg:pd_alg}. However, it is not clear whether this strategy would yield any numerical advantage.

%% 2.3. Worst-case complexity estimates.
\subsection{Convergence analysis}
Under Assumption \ref{as:A1}, the dual solution set $\mathcal{Y}^{\star}$ is nonempty. Recall that $D_{\mathcal{Y}^{\star}} := \displaystyle\min_{\yb^{\star}\in\mathcal{Y}^{\star}}\norm{\yb^{\star}}_2 < +\infty$ is the norm of a minimum norm dual solution. 
The following theorem shows the convergence of  Algorithm  \ref{alg:pd_alg}.

% Theorem 3.1.
\begin{theorem}\label{th:convergence_A1}
Let $\set{(\bar{\xb}^k, \bar{\yb}^k)}_{k\geq 0}$ be the sequence generated by Algorithm \ref{alg:pd_alg} after $k \geq 1$ iterations.
Then, if $g_{\gamma} \equiv \tilde{g}_{\gamma}$, i.e., using augmented Lagrangian smoother $\tilde{g}_{\gamma}$, then:
\begin{itemize}
\item[$\mathrm{a})$]
If  $c_k := 0$ for all $k \geq 0$, $\gamma_0 := \bar{L}^{\tilde{g}} = 1$, then:
\begin{equation}\label{eq:convergence_A1}
\left\{\begin{array}{rcl}
 &\Vert\Ab\bar{\xb}^k \!-\! \bb\Vert_2 \leq &\frac{8D_{\mathcal{Y}^{\star}}}{(k+1)^2},\\
-\frac{1}{2}\Vert\Ab\bar{\xb}^k \!-\! \bb\Vert_2^2 \!-\! D_{\mathcal{Y}^{\star}}\Vert\Ab\bar{\xb}^k \!-\! \bb\Vert_2 \leq & f(\bar{\xb}^k) \!-f^{\star} \leq &0,
\end{array}\right.
\end{equation}
for all $k \geq 0$. Moreover, the spectral norm of $\Ab$ does not affect the bounds in \eqref{eq:convergence_A1}.
\end{itemize}
As a consequence, the worst-case analytical complexity of Algorithm \ref{alg:pd_alg} to achieve an $\varepsilon$-primal solution $\bar{\xb}^k$ for \eqref{eq:separable_convex} in the sense of Definition \ref{de:approx_sol} is $\mathcal{O}\left(\varepsilon^{-1/2}\right)$.

Alternatively, if $g_{\gamma} \equiv \hat{g}_{\gamma}$, i.e., using Bregman distance smoother $\tilde{g}_{\gamma}$, then:
\begin{itemize}
\item[$\mathrm{b})$]
If Algorithm \ref{alg:pd_alg} uses \eqref{eq:pd_scheme_2p}, $\gamma_0 :=\sqrt{\bar{L}^g}$ and $c_k := 1$ for all $k \geq 0$, then:
\begin{equation}\label{eq:convergence_A2a}
\left\{\begin{array}{rcl}
&\Vert\Ab\bar{\xb}^k \!-\! \bb\Vert_2 &\leq \frac{\sqrt{\bar{L}^g}\left(2D_{\mathcal{Y}^{\star}} + \sqrt{2D_{\Xc}^{\mathbb{I}}}\right)}{k+1},\\
-D_{\mathcal{Y}^{\star}}\Vert\Ab\bar{\xb}^k \!-\! \bb\Vert_2 \leq &f(\bar{\xb}^k) - f^{\star} &\leq \frac{\sqrt{\bar{L}^g}}{k + 1}D_{\Xc}^{\mathbb{I}}.
\end{array}\right.
\end{equation}

\item[$\mathrm{c})$]
If Algorithm \ref{alg:pd_alg} uses \eqref{eq:pd_scheme_2d}, $\gamma_0 := \frac{2\sqrt{2\bar{L}_g}}{K+1}$ and $c_k := 0$ for all $k= 0, \ldots, K$, then:
\begin{equation}\label{eq:convergence_A2b}
\left\{\begin{array}{rcl}
&\Vert\Ab\bar{\xb}^K \!-\! \bb\Vert_2 &\leq \frac{2\sqrt{2\bar{L}_g}(D_{\mathcal{Y}^{\star}} + \sqrt{D_{\Xc}^{\mathbb{I}}})}{(K+1)},\\
-D_{\mathcal{Y}^{\star}}\Vert\Ab\bar{\xb}^K \!-\! \bb\Vert_2 \leq& f(\bar{\xb}^K) - f^{\star} &\leq \frac{2\sqrt{2\bar{L}^g}}{(K + 1)}D_{\Xc}^{\mathbb{I}}.
\end{array}\right.
\end{equation}
\end{itemize}
As a consequence, the worst-case analytical complexity of Algorithm \ref{alg:pd_alg} to achieve an $\varepsilon$-primal solution $\bar{\xb}^k$ for \eqref{eq:separable_convex} in the sense of Definition \ref{de:approx_sol} is $\mathcal{O}\left(\varepsilon^{-1}\right)$. 
\end{theorem}

We note that the choice of $\gamma_0$ in Theorem \ref{th:convergence_A1} trades-off the primal objective residual and the primal feasibility gap. 
Indeed, smaller $\gamma_0$ leads to smaller  $\vert f(\bar{\xb}^k) - f^{\star}\vert$.

We chose the $(\mathrm{1P2D})$ scheme above due to its close relationship to some well-known primal dual methods we describe below. Unfortunately, the $(\mathrm{1P2D})$ scheme has the drawback of fixing the total number of iterations \emph{a priori}, which the $(\mathrm{2P1D})$ scheme can avoid at the expense of more proximal operator calculations.

\section{Instances of Algorithm \ref{alg:pd_alg}}\label{sec:variants}
This section specifies Algorithm \ref{alg:pd_alg} under different assumptions  to obtain specific instances of this algorithm for solving \eqref{eq:separable_convex}.

%% 6.1. Strong convexity assumption.
\subsection{Strong convexity assumption}
If the objective function $f$ of \eqref{eq:separable_convex} is strongly convex with a convexity parameter $\sigma_f > 0$. Then it is well-known that  (see, e.g., \cite{Nesterov2005c}) the dual function $g(\cdot)$ defined by \eqref{eq:dual_func} is smooth and Lipschitz gradient with a Lipschitz constant $L^g_f := \frac{\norm{\Ab}_2^2}{\sigma_f}$.
In this case, we modify accordingly both schemes \eqref{eq:pd_scheme_2p} and \eqref{eq:pd_scheme_2d} as follows:
\begin{equation*}
(\mathrm{2P1D}_{\sigma})
\left\{\begin{array}{ll}
\hat{\xb}^k        &  {\!\!\!\!\!}:= (1-\tau_k)\bar{\xb}^k + \tau_k\xb^{\star}(\bar{\yb}^k),\\
\bar{\xb}^{k\!+\!1}  &{\!\!\!\!\!}:=  \mathrm{prox}_{\mathbb{I}f}(\hat{\xb}^k, \beta_{k}^{-1}\!(\Ab\hat{\xb}^k \!-\! \bb); \beta_{k}),\\
\bar{\yb}^{k \!+\! 1} &{\!\!\!\!\!} := (1-\tau_k)\bar{\yb}^k + \frac{\tau_k}{\beta_k}\big(\Ab\hat{\xb}^k - \bb\big).
\end{array}\right.
~~~~~(\mathrm{1P2D}_{\sigma})
\left\{\begin{array}{ll}
\hat{\yb}^k          &{\!\!\!\!\!}  := (1 \!-\! \tau_k)\bar{\yb}^k \!+\! \tau_k\yb^{\star}_{\beta_k}(\bar{\xb}^k),\\
\bar{\xb}^{k \!+\! 1}   &{\!\!\!\!\!} := (1 \!-\! \tau_k)\bar{\xb}^k \!+\! \tau_k\xb^{\star}(\hat{\yb}^k),\\
\bar{\yb}^{k \!+\! 1} & {\!\!\!\!\!} := \hat{\yb}^k \!+\! \frac{1}{L^g_f}\big(\Ab\xb^{\star}(\hat{\yb}^k) \!-\! \bb\big).
\end{array}\right. 
\end{equation*}
While the scheme $(\mathrm{1P2D}_{\sigma})$ remains similarly to \eqref{eq:pd_scheme_2d}, the parameter $\beta_k$ in $(\mathrm{2P1D}_{\sigma})$ has not updated yet as in \eqref{eq:pd_scheme_2p}.

The starting point  $\bar{\wb}^0 := (\xb^{\star}(0^m), \bar{\yb}^0)\in\Wc$ for Algorithm \ref{alg:pd_alg} with respect to this variant can be computed as $\bar{\yb}^0 := (L_f^g)^{-1}(\Ab\xb^{\star}(0^m) - \bb)$ and $\xb^{\star}(\yb)$ is the unique solution of the minimization in \eqref{eq:dual_func}.
The parameters $\beta_k$ and $\tau_k$ are updated as follows:
\begin{equation}\label{eq:update_tau_beta}
\beta_{k+1} := (1-\tau_k)\beta_k, ~~\tau_{k+1} := (\tau_k/2)[(\tau_k^2 + 4)^{1/2} - \tau_k], ~~~k\geq 0,
\end{equation}
where $\beta_0 := L_f^g$ and $\tau_0 := (\sqrt{5}-1)/2$.
The following corollary shows the convergence of both schemes, whose proof is in the appendix.

% Corollary 6.1.
\begin{corollary}\label{co:strong_convex_convergence}
Assume that $f$ of \eqref{eq:separable_convex} is $\sigma_f$-strongly convex.
Let $\set{(\bar{\xb}^k, \bar{\yb}^k)}_{k\geq 0}$ be a sequence generated by either $(\mathrm{2P1D}_{\sigma})$ or $(\mathrm{1P2D}_{\sigma})$ using the update rule \eqref{eq:update_tau_beta}. Then
\begin{equation}\label{eq:strong_cvx_main_estimate}
\left\{\begin{array}{rcl}
&\Vert\Ab\bar{\xb}^k - \bb\Vert_2 \leq &\frac{4\norm{\Ab}_2^2}{(k+2)^2\sigma_f}D_{\Yc^{\star}},\\
-D_{\Yc^{\star}}\Vert\Ab\bar{\xb}^k - \bb\Vert_2 \leq &f(\bar{\xb}^k) - f^{\star} \leq & 0, \\
&\Vert\bar{\xb}^k - \xb^{\star}\Vert_2 \leq &\frac{4\norm{\Ab}_2}{(k+2)\sigma_f}D_{\Yc^{\star}},
\end{array}\right.
\end{equation}
where $D_{\Yc^{\star}}$ is defined in Theorem \ref{th:convergence_A1} and $\xb^{\star}\in\Xc^{\star}$.

As a consequence, the worst-case analytical complexity for finding an $\varepsilon$-primal solution $\bar{\xb}^k$ of \eqref{eq:separable_convex} in the sense of Definition \ref{de:approx_sol} is $\mathcal{O}(1/\sqrt{\varepsilon})$.
\end{corollary}

\begin{remark}
 The bounds in \eqref{eq:strong_cvx_main_estimate} do not depend on the prox-diameter $D_{\Xc}^{\mathbb{I}}$ of the feasible set $\Xc$. 
 Hence, the boundedness of $\Xc$ is no longer required.
\end{remark}

\begin{remark}
Convergence of the objective indeed depends on the absolute value of the primal residual, i.e.,  $\vert f(\bar{\xb}^k) - f^{\star}\vert \leq \frac{4\Vert\Ab\Vert_2^2}{(k+2)^2\sigma_f}D_{\Yc^{\star}}^2$. 
\end{remark}

%%%%%%%%%%%%%%%%%%%%%%%%%%%%%%%%%%%%%
%%% 6.2. Lipschitz gradient Assumption.
%%%%%%%%%%%%%%%%%%%%%%%%%%%%%%%%%%%%%
\subsection{Lipschitz gradient assumption}
The aim of this subsection is to develop a variant of Algorithm \ref{alg:pd_alg} using  \eqref{eq:pd_scheme_2d} without fixed the accuracy as stated in Theorem \ref{th:convergence_A1}(c).
However, this variant is only limited to problems of the form \eqref{eq:separable_convex} that satisfy the following technical assumption:

%% Assumption A.3.
\begin{assumption}\label{ass:A3}
The following conditions hold:
\begin{itemize}
\item[(a)] The objective function $f$ and the feasible set $\Xc$ of \eqref{eq:separable_convex} are separable as in \eqref{eq:separable_structure}.
\item[(b)] The last term $f_p$ is $L_{f_p}$-Lipschitz gradient and the smallest eigenvalue $\lambda_{\min}(\Ab_p^T\Ab_p)$ of matrix $\Ab_p$ is positive.
\item[(c)] The Bregman distance $d(\Sb\xb,\Sb\xb_c)$ is chosen as $d(\Sb\xb,\Sb\xb_c) := \sum_{i=1}^{p}d_i(\Sb_i\xb_i,\Sb_i\xb_i^c)$, 
where $\Sb_p\equiv\mathbb{I}$ and $d_p(\cdot, \xb_p^c)$ is smooth and $\nabla{d_p}(\cdot, \xb_p^c)$ is $1$-Lipschitz continuous.
\item[(d)] The last term $g_{\gamma}^p$ of the smoothed dual function $g_{\gamma}$ defined by  \eqref{eq:g_gamma} satisfies 
\begin{equation}\label{eq:g_gamma^p}
g_{\gamma}^p(\yb) = \min_{\xb_p\in\R^{n_p}}\set{ f_p(\xb_p) + \yb^T\Ab_p\xb_p + (\gamma/2)d_p(\xb_p, \xb_p^c)}.
\end{equation}
That is the primal constraint on the last component is not active. 
\end{itemize}
\end{assumption}
Under Assumption A.\ref{ass:A3}, we can write the function $g_{\gamma}$ defined by \eqref{eq:g_gamma} as $g_{\gamma}(\yb) := \sum_{i=1}^pg_{\gamma}^i(\yb) - \bb^T\yb$, where
\begin{equation*}
g_{\gamma}^i(\yb) := \min_{\xb_i\in\Xc_i}\set{ f_i(\xb_i) + \yb^T\Ab_i\xb_i + (\gamma/2)d_i(\Sb_i\xb_i, \Sb_i\xb_i^c)}, ~~i=1,\dots, p.
\end{equation*}
A simple example for $d_p$ is $d_p(\xb_p) := (1/2)\Vert\xb_p - \xb_p^c\Vert_2^2$.
The last condition in Assumption A.\ref{ass:A3} shows that the solution $\xb^{\star}_{p,\gamma}(\yb)$ of the minimization problem in $g_{\gamma}^p$ must be attained in $\mathrm{relint}(\X_p)$.
This condition is not too restrictive, since we only require it for the last component $g_{\gamma}^p$. It is automatically fulfilled if $f_p$ is strongly convex and $\xb_p^c\in \mathrm{relint}(\X_p)$.
Now, we show that the function $g_{\gamma}^p$ is strongly concave in the following lemma, whose proof can be found in the appendix.

% Lemma 2.1.
\begin{lemma}\label{le:strongly_concave}
Under Assumption A.\ref{ass:A3}, the function $g_{\gamma}^p$ defined by \eqref{eq:g_gamma^p} is strongly concave with the parameter $\sigma_{g_{\gamma}^p} := (L_{f_p}+\gamma)^{-1}\lambda_{\min}(\Ab_p^T\Ab_p) > 0$. Consequently, the function $g_{\gamma}$ defined by \eqref{eq:g_gamma} is also strongly convex with the same parameter $\sigma_{g_{\gamma}^p}$.
\end{lemma}
Using the result of Lemma \ref{le:strongly_concave}, we can update $\gamma_k$ and $\beta_k$ in the scheme \eqref{eq:pd_scheme_2d} as 
\begin{equation}\label{eq:update_tau_beta2}
\gamma_{k+1} := \big(1 -  \tau_k/(1+\tau_k)\big)\gamma_k,~~\beta_{k+1} := (1-\tau_k)\beta_k ~~\textrm{and}~~\tau_k := (k+1)^{-1} ~~\forall k\geq 0,
\end{equation}
where $\beta_0 = \gamma_0 := \sqrt{\bar{L}^g}$.
In this case, we have $\gamma_{k+1}\beta_{k+1} \geq \bar{L}^{g}\tau_k^2$ for $k\geq 0$.
The following corollary shows the convergence of this variant.

% Corollary 6.2.
\begin{corollary}\label{co:Lipschitz_convergence}
Under Assumption A.\ref{ass:A3}, let $\set{(\bar{\xb}^k, \bar{\yb}^k)}_{k\geq 0}$ be a sequence generated by \eqref{eq:pd_scheme_2d} using the update rule \eqref{eq:update_tau_beta2}. Then
\begin{equation}\label{eq:Lipschitz_main_estimate}
\left\{\begin{array}{rcl}
&\Vert\Ab\bar{\xb}^k \!-\! \bb\Vert_2 &\leq \frac{2\sqrt{2\bar{L}^g}\big(D_{\mathcal{Y}^{\star}} + \sqrt{D_{\Xc}^{\Sb}}\big)}{k+1},\\
-D_{\mathcal{Y}^{\star}}\Vert\Ab\bar{\xb}^k \!-\! \bb\Vert_2 \leq&  f(\bar{\xb}^k) - f^{\star} &\leq \frac{2\sqrt{2\bar{L}^g}}{k + 1}D_{\Xc}^{\Sb},\\
\end{array}\right.
\end{equation}
where $D_{\Yc^{\star}}$ and $\bar{L}^g$ are defined in Theorem \ref{th:convergence_A1}.
\end{corollary}

\begin{remark}
Corollary \ref{co:Lipschitz_convergence} shows that, for certain subclass of problems \eqref{eq:separable_convex} satisfying Assumption A.\ref{ass:A3},  it allows us to simultaneously update both parameters $\gamma_k$ and $\beta_k$ instead of fixing a priori $\gamma_0$ as in Theorem \ref{th:convergence_A1}(c).
\end{remark}

%\begin{remark}
%If $f$ is simultaneously $\sigma_f$-strongly convex and Assumption A.\ref{ass:A3} is fulfilled, then it is easy to show that the convergence of $\Vert\Ab\bar{\xb}^k \!-\! \bb\Vert_2$ and $\vert f(\bar{\xb}^k) - f^{\star} \vert$ is linear, which we omit the proof details here.
%\end{remark}

%%% 5.3. Inexact solution of the augmented Lagrangian smoother
\subsection{Inexact solution of the augmented Lagrangian smoother}
In the augmented Lagrangian smoothing method, solving the minimization problem \eqref{eq:aug_dual_func} exactly can be impracticable. However, we can often solve this subproblem up to a given accuracy $\delta >0$, i.e.,
\begin{equation}\label{eq:approx_sol}
\tilde{\xb}^{\delta}_{\gamma}(\yb) := \delta\textrm{-}\arg\min_{\xb\in\Xc}\set{\mathcal{L}_{\gamma}(\xb, \yb) := f(\xb) + \yb^T(\Ab{\xb} - \bb) + (\gamma/2)\norm{\Ab{\xb} - \bb}_2^2},
\end{equation}
in the following sense:
\begin{equation}\label{eq:aug_Lagrangian_inexact}
\mathcal{L}_{\gamma}(\tilde{\xb}^{\delta}_{\gamma}(\yb), \yb) - \mathcal{L}_{\gamma}(\tilde{\xb}^{\star}_{\gamma}(\yb), \yb) \leq \gamma\delta^2/2,
\end{equation}
where $\tilde{\xb}^{\star}_{\gamma}(\yb)$ is an exact solution of \eqref{eq:aug_dual_func}.

The condition $\tilde{\xb}^{\delta}_{\gamma}(\yb) \in \Xc$  is reasonable in practice since the feasible set $\Xc$ can be assumed to be ``simple'' so that the computation of the projection onto $\Xc$ can be carried out exactly. 
In addition, there exist several convex optimization algorithms (e.g., Nesterov's accelerated algorithms \cite{Nesterov2004}) for computing $\tilde{\xb}^{\delta}_{\gamma}(\yb)$ that satisfy  \eqref{eq:aug_Lagrangian_inexact}.

By the definition of $\mathcal{L}_{\gamma}$, we can easily show that 
\begin{equation*}
\mathcal{L}_{\gamma}(\tilde{\xb}^{\delta}_{\gamma}(\yb), \yb) - \mathcal{L}_{\gamma}(\tilde{\xb}^{\star}_{\gamma}(\yb), \yb) \geq (\gamma/2)\Vert\Ab(\tilde{\xb}^{\delta}_{\gamma}(\yb) - \tilde{\xb}^{\star}_{\gamma}(\yb))\Vert_2^2,
\end{equation*}
which leads to $\Vert\Ab(\tilde{\xb}^{\delta}_{\gamma}(\yb) - \tilde{\xb}^{\star}_{\gamma}(\yb))\Vert_2 \leq \delta$. 
Now, if we define $\nabla{\tilde{g}}_{\gamma}^{\delta}(\yb) := \Ab\tilde{\xb}^{\delta}_{\gamma}(\yb) - \bb$ an approximation for the gradient $\nabla{\tilde{g}_{\gamma}}(\yb)$, then \eqref{eq:aug_Lagrangian_inexact} and the last inequality implies
\begin{equation}\label{eq:aug_Lagrangian_sol_inexact}
\Vert\nabla{\tilde{g}}^{\delta}_{\gamma}(\yb) - \nabla{\tilde{g}}_{\gamma}(\yb)\Vert_2 \leq \delta.
\end{equation}
In addition, we also denote by $\tilde{g}^{\delta}_{\gamma}(\yb) := \tilde{\mathcal{L}}_{\gamma}(\tilde{\xb}^{\delta}_{\gamma}(\yb), \yb)$ as an approximation to $\tilde{g}_{\gamma}(\yb)$.

Instead of using the true solution $\xb^{\star}_{\gamma}(\yb)$ in the schemes \eqref{eq:pd_scheme_2p} and \eqref{eq:pd_scheme_2d}, we use the approximate solutions $\tilde{\xb}^{\delta}_{\gamma}(\yb)$ to obtain the following inexact iterative schemes:
\begin{equation}\label{eq:inexact_pd_scheme}
\boxed{
\left\{\begin{array}{ll}
&~~~~~~~~~(\textrm{i2P1D})\\
\hat{\xb}^k &{\!\!\!\!}        :=   (1 - \tau_k)\bar{\xb}^k + \tau_k\tilde{\xb}^{\delta_k}_{\gamma_k}(\bar{\yb}^k),\\
\hat{\yb}^k  &{\!\!\!\!}       :=   \beta_{k+1}^{-1}(\Ab\hat{\xb}^k - \bb),\\
\bar{\xb}^{k+1}  &{\!\!\!\!}:=  \widetilde{\mathrm{prox}}^{\delta_k}_{\Ab f}(\hat{\xb}^k, \hat{\yb}^k; \beta_{k+1}),\\
\bar{\yb}^{k+1}   &{\!\!\!\!}:=  (1 - \tau_k)\bar{\yb}^k + \tau_k\hat{\yb}^k.
\end{array}\right.
}~~~~~
\boxed{
\left\{\begin{array}{ll}
&~~~~~~~~~(\textrm{i1P2D})\\
\bar{\yb}^{\star}_k  &{\!\!\!}      :=  \beta_k^{-1}(\Ab\bar{\xb}^k - \bb),\\ 
\hat{\yb}^k        &{\!\!\!}      :=  (1-\tau_k)\bar{\yb}^k + \tau_k\bar{\yb}^{\star}_k,\\
\bar{\xb}^{k\!+\!1} &{\!\!\!}  :=  (1 \!-\! \tau_k)\bar{\xb}^k \!+\! \tau_k\tilde{\xb}^{\delta_k}_{\gamma_{k}}\!(\hat{\yb}^k),\\
\bar{\yb}^{k\!+\!1} &{\!\!\!}  :=  \hat{\yb}^k \!+\! \gamma_k\big(\Ab\tilde{\xb}^{\delta_k}_{\gamma_{k}}\!(\hat{\yb}^k) \!-\! \bb\big).
\end{array}\right.
}
\end{equation}
Here, the inexact proximal operator $\widetilde{\mathrm{prox}}^{\delta}_{\Ab f}$ is defined as:
\begin{equation}\label{eq:prox_oper2_inexact}
\widetilde{\mathrm{prox}}^{\delta}_{\Ab f}(\bar{\xb}, \hat{\yb}; \beta) :=\delta\textrm{-}\mathrm{arg}\!\min_{\xb\in\Xc}\set{\mathcal{H}_{\beta}(\xb; \hat{\yb}, \bar{\xb}) \!:=\! f(\xb) \!+\! \hat{\yb}^T\Ab(\xb \!-\! \bar{\xb}) \!+\! \frac{\bar{L}^g}{2\beta}\Vert\Ab(\xb \!-\! \bar{\xb})\Vert^2_2},
\end{equation}
where $\Ab$ and $\delta \geq 0$ are given and the inexactness is also defined as in \eqref{eq:aug_Lagrangian_inexact}.

The starting point $\bar{\wb}^0 := (\bar{\xb}^0, \bar{\yb}^0)\in\Wc$ can be computed from one of the following formulations:
\begin{equation}\label{eq:starting_point1_inexact}
\left\{\begin{array}{ll}
\bar{\xb}^0  &:= \tilde{\xb}^{\delta_0}_{\gamma_0}(0^m), \\
\bar{\yb}^0 &:=  \beta_0^{-1}\left(\Ab\bar{\xb}^0 - \bb\right),
\end{array}\right.
~~\textrm{or}~~
\left\{\begin{array}{ll}
\bar{\yb}^0 &:= \beta_0^{-1}\left(\Ab\xb_c - \bb\right),\\
\bar{\xb}^0  & :=  \widetilde{\mathrm{prox}}^{\delta_0}_{\Ab f}(\xb_c, \bar{\yb}^0;\beta_0).
\end{array}\right.
\end{equation}

The following theorem shows the convergence of the inexact variant of Algorithm \ref{alg:pd_alg} using scheme \eqref{eq:inexact_pd_scheme}, called  $(\mathrm{i1P2D})$, whose proof can be found in the appendix. Analogously, we can also prove the same result as  in Theorem \ref{th:inexact_convergence} for the $(\mathrm{i2P1D})$ scheme but we omit the laborious details.

% Theorem .
\begin{theorem}\label{th:inexact_convergence}
Let $\set{(\bar{\xb}^k, \bar{\yb}^k)}_{k\geq 0}$ be the sequence generated by Algorithm \ref{alg:pd_alg}  using $(\mathrm{i1P2D})$ in \eqref{eq:inexact_pd_scheme} and the first initial point $(\bar{\xb}^0, \bar{\yb}^0)$ in \eqref{eq:starting_point1_inexact}.
Then, if  $\gamma_0 = \bar{L}^{\tilde{g}} = 1$, $c_k := 0$ and $q_k\delta_k \leq q_{k-1}\delta_{k-1}$ for all $k \geq 0$ and $q_k := (1-\tau_k)\tau_k\norm{\bar{\yb}^k \!-\! \bar{\yb}^{\star}_k}_2 + (D_{\Xc}^{\Ab} + 1)/2$ then:
\begin{equation}\label{eq:convergence_inexact}
\left\{\begin{array}{rcl}
&{\!\!\!\!\!} \Vert\Ab\bar{\xb}^k \!-\! \bb\Vert_2 &\leq \frac{4}{(k+1)^2}\Big( 2D_{\Yc^{\star}} + \sqrt{\frac{14q_0\delta_0}{(k+1)^2}}\Big),\\
{\!\!\!\!\!}-(1/2)\Vert\Ab\bar{\xb}^k \!-\! \bb\Vert_2^2 \!-\! D_{\Yc^{\star}} \Vert\Ab\bar{\xb}^k \!-\! \bb\Vert_2 \leq & f(\bar{\xb}^k) - f^{\star}  &\leq  7q_0\delta_0.
\end{array}\right.
\end{equation}
As a consequence, if $\delta_0 = \mathcal{O}\left(\frac{q_0}{k^2}\right)$, then  the worst-case analytical complexity of Algorithm \ref{alg:pd_alg} to achieve an $\varepsilon$-primal solution $\bar{\xb}^k$ of \eqref{eq:separable_convex} in the sense of Definition \ref{de:approx_sol} is $\mathcal{O}\left(\varepsilon^{-1/2}\right)$.
\end{theorem}

Theorem \ref{th:inexact_convergence} shows that the primal feasibility gap $\Vert\Ab\bar{\xb}^k - \bb\Vert_2$ converges to $0^{+}$ at the rate $\mathcal{O}(1/k^2)$, while the objective residual $\abs{f(\bar{\xb}^k) - f^{\star}}$ depends on the numerical accuracy $\delta_0$ of \eqref{eq:approx_sol} at the initial iteration $k = 0$.
If $\delta_0$ is not sufficiently small, we only obtain a sub-optimal solution of \eqref{eq:separable_convex}. 
Practically, we can solve \eqref{eq:approx_sol} at $k=0$ with relatively high accuracy and use a warm-start strategy to significantly reduce the computational burden of the subsequent iterations. 

%%%%%%%%%%%%%%%%%%%%%%%%%%%%%%%%%%%%%%%%%%%%
%% 7. Connection to existing methods}
%%%%%%%%%%%%%%%%%%%%%%%%%%%%%%%%%%%%%%%%%%%%
\section{Explicit connections to existing methods}\label{sec:existing_methods}\label{sec: explicit}
To better differentiate our contributions, it is important to make explicit comparisons of Algorithm \ref{alg:pd_alg} with the  dual fast gradient methods, alternating direction methods of multipliers (ADMM) and proximal-based decomposition methods here. 
%
%It is important to  concrete connections to existing primal-dual methods much more concrete. 
%We only discuss the connection between  and the following related methods: .
\subsection{Connections to the fast gradient methods}
Dual fast gradient methods were studied in, e.g., \cite{Beck2014,Necoara2008,Nedelcu2014,Polyak2013}. 
The main idea is to use either the strong convexity of the objective \cite{Beck2014,Polyak2013} or smoothing technique via  prox-functions \cite{Necoara2008} or augmented Lagrangian function \cite{Nedelcu2014}, which leads to the Lipschitz continuity of the gradient of the dual function. Then, Nesterov's fast gradient method \cite{Nesterov2004} is applied to solve the smoothed dual problem.

In this paper, we also smooth the dual function by using either augmented Lagrangian function or Bregman distances to obtain a smoothed dual function with Lipschitz continuous gradient. 
In order to obtain both primal objective residual and primal feasibility gap simultaneously, we exploit the concepts of excessive gap technique introduced by Nesterov \cite{Nesterov2005d} and Auslander's gap function \cite{Auslender1976} to build a primal-dual sequence $\set{(\bar{\xb}^k, \bar{\yb}^k)}_{k\geq 0}$ that converges to the primal-dual optimal solution $(\xb^{\star}, \yb^{\star})$ of \eqref{eq:separable_convex}.
In \cite{Necoara2008,Polyak2013} the authors only proved the convergence results in terms of the dual objective values $g$, which is different from Theorem \ref{th:convergence_A1}, where we both have the convergence rate guarantee both on the primal objective residual and the primal feasibility gap.
In \cite{Nedelcu2014} the authors characterized the convergence rate of an inexact augmented Lagrangian method both in the primal objective values and the primal feasibility gaps. 
However, the approach  is directly based on Nesterov's accelerated scheme for the dual problem and the convergence results are presented in an ergodic sense.
In \cite{Beck2014} the authors considered a special case of \eqref{eq:separable_convex}, where the objective function is strongly convex as in Corollary \ref{co:strong_convex_convergence}. They also characterized the feasibility gap. However, the convergence rate of this quantity drops to $\mathcal{O}(1/k)$ instead of the better $\mathcal{O}(1/k^2)$ rate established by our Corollary \ref{co:strong_convex_convergence}.

We close this discussion by showing that our results in Corollary \ref{co:strong_convex_convergence} can be applied to non-strongly convex problems of the form \eqref{eq:separable_convex}.
We process this procedure as follows. 
Assume that $f$ of \eqref{eq:separable_convex} is not strongly convex, we consider the function $f_{\sigma}(\xb) := f(\xb) + (\sigma_f/2)\Vert\xb - \xb_c\Vert_2^2$, where $\sigma_f > 0$ and $\xb_c\in\Xc$. Then, the function $f_{\sigma}$ is strongly convex with the parameter $\sigma_f > 0$.
Next, we apply either $(\mathrm{2P1D}_{\sigma})$ or $(\mathrm{1P2D}_{\sigma})$ to solve \eqref{eq:separable_convex} with $f$ is substituted by $f_{\sigma}$. 
In this case, Corollary \ref{co:strong_convex_convergence} is still valid. Moreover, we have $f_{\sigma}(\bar{\xb}^k) = f(\bar{\xb}^k) + (\sigma_f/2)\Vert\bar{\xb}^k - \xb_c\Vert_2^2$ and $f^{\star}_{\sigma} = f^{\star} + (\sigma_f/2)\Vert\xb^{\star} - \xb_c\Vert_2^2$, which imply 
\begin{equation*}
\vert f(\bar{\xb}^k) - f^{\star}\vert \leq \vert f_{\sigma}(\bar{\xb}^k) - f_{\sigma}^{\star}\vert + 2\sigma_f D_{\Xc}^{\mathbb{I}},
\end{equation*}
where $D_{\Xc}^{\mathbb{I}} := \max_{\xb\in\Xc}(1/2)\Vert\xb - \xb_c\Vert_2^2$.
Combining this estimate and Corollary \ref{co:strong_convex_convergence} we obtain $\vert f(\bar{\xb}^k) - f^{\star}\vert \leq \frac{4\Vert\Ab\Vert_2^2}{\sigma_f(k+2)^2}D_{\mathcal{Y}^{\star}}^2 + 2\sigma_fD_{\Xc}^{\mathbb{I}}$. Hence, if we choose $\sigma_f := \frac{\sqrt{2}\Vert\Ab\Vert_2D_{\mathcal{Y}^{\star}}}{(k+2)\sqrt{D_{\Xc}^{\mathbb{I}}}}$ then we obtain the worst-case analytical complexity of this algorithm as 
\begin{equation*}
\vert f(\bar{\xb}^k) - f^{\star}\vert \leq \frac{2\sqrt{2}\Vert\Ab\Vert_2D_{\mathcal{Y}^{\star}}(D_{\Xc}^{\mathbb{I}})^{1/2}}{(k+2)}~~\textrm{and}~~\Vert \Ab\bar{\xb}^k -\bb\Vert_2 \leq \frac{2\sqrt{2}\Vert\Ab\Vert_2(D_{\Xc}^{\mathbb{I}})^{1/2}}{(k+2)}.
\end{equation*}
Comparing this complexity and Theorem \ref{th:convergence_A1}, we conclude that depending on the values of $D_{\mathcal{Y}}^{\star}$ and $D_{\Xc}^{\mathbb{I}}$ we can use choose an appropriate variant of Algorithm \ref{alg:pd_alg} for solving the given problem. However, note that we do not generally have access to $D_{\mathcal{Y}^{\star}}$, hence we can instead use the standard (1P2D) or (2D1P) schemes which do not require the knowledge of the smoothing parameter.

%% 7.1. Connection to alternating direction methods of multipliers.
\subsection{Connection to ADMMs}
Several algorithms based on method of multipliers such as alternating minimization algorithm (AMA) \cite{Tseng1991a}, alternating direction method of multipliers (ADMM) \cite{Bertsekas1989b} and alternating linearization methods (ALM) \cite{Goldfarb2012} have been developed in the literature. Such methods aim at solving instances of \eqref{eq:separable_convex} when $f$ and $\Xc$ are separable with $p = 2$. In this case, the primal step is computed by solving two subproblems w.r.t. $\xb_{[1]}$ and $\xb_{[2]}$ alternatively. 

Let $f(\xb) := f_1(\xb_1) + f_2(\xb_2)$ and $\Xc := \Xc_1\times \Xc_2$. The standard ADMM algorithm \cite{Tseng1991a} can be presented as follows:
\begin{equation}\label{eq:alternating_method}
\left\{\begin{array}{ll}
\xb_1^{k+1} &{\!\!\!\!\!\!\!} := \mathrm{arg}\!\!\!\displaystyle\min_{\xb_1\in\Xc_1}\set{f_1(\xb_1)  \!+\! \frac{\eta_k}{2}\norm{\Ab_1\xb_1 \!+\! \Ab_2\xb^k_2 \!-\! \bb + \eta_k^{-1}\yb^k}_2^2},\\
\xb^{k+1}_2 &{\!\!\!\!\!\!\!} := \mathrm{arg}\!\!\displaystyle\min_{\xb_2\in\Xc_2}\set{f_2(\xb_2)  \!+\! \frac{\eta_k}{2}\norm{\Ab_1\xb_1^{k\!+\!1} + \Ab_2\xb_2 - \bb + \eta_k^{-1}\yb^k}_2^2},\\
\yb^{k+1} &{\!\!\!\!\!\!\!}  := \yb^k + \eta_k(\Ab_1\xb_1^{k+1} + \Ab_2\xb_2^{k+1} - \bb),
\end{array}\right.
\end{equation}
where $\eta_k > 0$ is a penalty parameter. 
ADMM works very well in practice and has been widely used in many disciplines. 
When $f_i$ is tractably proximal and $\Ab_i^T\Ab_i$  is diagonalizable, the solutions $\xb_i$ can be computed efficiently ($i = 1,2$). In the opposite case, computing $\xb_1^{k+1}$ and $\xb_2^{k+1}$ may require an iterative algorithm.

Let us modify the \eqref{eq:pd_scheme_2d} scheme by using the primal step as in \eqref{eq:alternating_method} to obtain:
\begin{equation}\label{eq:ADMM_variants}
\left\{\begin{array}{ll}
\hat{\yb}^k & {\!\!\!\!\!\!\!} := (1-\tau_k)\bar{\yb}^k + \tau_k\beta_k^{-1}(\Ab\bar{\xb}^k - \bb),\\
\xb^{k+1}_1 &{\!\!\!\!\!\!\!} := \mathrm{arg}\displaystyle\min_{\xb_1\in\Xc_1}\set{f_1^k(\xb_1) \ \!+\! \frac{\rho_k}{2}\norm{\Ab_1\xb_1 \!+\! \Ab_2\xb^k_2 \!-\! \bb + \rho_k^{-1}\yb^k}_2^2},\\
\xb^{k+1}_2 &{\!\!\!\!\!\!\!} := \mathrm{arg}\displaystyle\min_{\xb_2\in\Xc_2}\set{f_2(\xb_2)  \!+\! \frac{\eta_k}{2}\norm{\Ab_1\xb_1^{k\!+\!1} + \Ab_2\xb_2 - \bb  + \eta_k^{-1}\yb^k}_2^2},\\
\bar{\xb}_{k+1} &{\!\!\!\!\!\!\!} := (1-\tau_k)\bar{\xb}^k + \tau_k\xb^{k+1},\\
\bar{\yb}^{k+1} &{\!\!\!\!\!\!\!} := \hat{\yb}^k + \eta_k(\Ab_1\xb_1^{k+1} + \Ab_2\xb^{k+1}_2 - \bb),
\end{array}\right.
\end{equation}
where $f_1^k(\cdot) := f_1(\cdot) + (\gamma_{k+1}/2)\norm{\Ab_1(\xb_1-\bar{\xb}^c_1)}^2$ for a fixed $\bar{\xb}_1^c\in\Xc_1$.
It is trivial that if $\tau_k = 0$, $\gamma_{k+1} = 0$ and $\rho_k = \eta_k$, then \eqref{eq:ADMM_variants} coincides with the standard ADMM scheme \eqref{eq:alternating_method}.

As indicated in \cite{Tran-Dinh2015}, the parameters $\tau_k$, $\gamma_k$, $\beta_k$, $\rho_k$ and $\eta_k$ are updated by:
\begin{equation}\label{eq:update_admm_params}
\begin{array}{llll}
&\tau_k  := \frac{3}{k\!+\!4}, ~&\gamma_k := \frac{2\gamma_0}{k\!+\!2}, ~&\beta_k := \frac{9(k+3)}{\gamma_0(k\!+\!1)(k\!+\!7)},\\
&\rho_k  := \frac{3\gamma_0}{(k\!+\!3)(k\!+\!4)}, ~&\eta_k := \frac{\gamma_0}{k+3}, &
\end{array}
\end{equation}
where $\gamma_0 > 0$ is chosen arbitrarily to trade off the primal objective residual $\vert f(\bar{\xb}^k) - f^{\star}\vert$ and the primal feasibility gap $\norm{\Ab\bar{\xb}^k - \bb}$.
The starting point $\bar{\wb}^0 := [\bar{\xb}^0, \bar{\yb}^0]$ can be computed from the steps 2, 3 and 5 of \eqref{eq:ADMM_variants} by choosing $\hat{\yb}^0 = \boldsymbol{0}^m$.

The following corollary shows the convergence of the $\mathrm{PADMM}$ scheme \eqref{eq:ADMM_variants}-\eqref{eq:update_admm_params}, whose proof can be found in \cite{Tran-Dinh2015}.
%% Corollary 7.3.
\begin{corollary}\label{co:PADMM}
Let $\set{(\bar{\xb}^k, \bar{\yb}^k)}_{k\geq 0}$ be a sequence generated by Algorithm \ref{alg:pd_alg} using the $\mathrm{ADMM}$ scheme \eqref{eq:ADMM_variants}-\eqref{eq:update_admm_params}. 
If  $\gamma_0 := 3$, then:
\begin{equation}\label{eq:convergence_A2d}
\left\{\begin{array}{ll}
\abs{f(\bar{\xb}^k) - f^{*}} &\leq \frac{6D_{\max}}{k+2},\\
\Vert\Ab\bar{\xb}^k \!-\! \bb\Vert_2 &\leq \frac{6\left[D^{*}_{\mathcal{Y}} + \sqrt{D_{\Xc_1}^1 + 4D^{\Ab}_{\Xc}}\right]}{k+2}.
\end{array}\right.
\end{equation}
where $D_{\Xc_1}^1 := (1/2)\max\big\{ \norm{\Ab_1(\xb_1 - \bar{\xb}_1^c)}_2^2 : \xb_1 \in\Xc_1 \big\}$, $D^{\Ab}_{\Xc} := (1/2)\max\big\{ \norm{\Ab\xb - \bb}^2 : \xb\in\Xc \big\}$, and $D_{\max} := \max\big\{D_{\Xc_1}^1 + 3D_{\Xc}^{\Ab}, D_{\Yc^{*}}(D^{*}_{\mathcal{Y}} + \sqrt{D_{\Xc_1}^1 + 4D_{\Xc}^{\Ab}})\big\}$.
As a consequence, the worst-case complexity of Algorithm \ref{alg:pd_alg} to achieve an $\varepsilon$-primal-dual solution $(\bar{\xb}^k, \bar{\yb}^k)$ is $\mathcal{O}\left(\varepsilon^{-1}\right)$.
\end{corollary}

If $f_1$ and $f_2$ has a tractable proximal operator $\mathrm{prox}_{\lambda f}$ then instead of solving two minimization problems in \eqref{eq:alternating_method}, we can linearize the quadratic term to obtain a preconditioning $\mathrm{ADMM}$ ($\mathrm{PADMM}$) as considered in \cite{Chambolle2011}. In this case, the primal step \eqref{eq:alternating_method} becomes
\begin{equation}\label{eq:linearized_alternating_method}
\left\{\begin{array}{ll}
\xb_1^{k+1} &{\!\!\!\!\!\!\!} := \mathrm{arg}\displaystyle\min_{\xb_1\in\Xc_1}\Big\{f_1(\xb_1) + \frac{\kappa_k}{2\alpha_{1k}}\norm{\xb_1 - (\mathbf{g}_1^k + \kappa_k^{-1}(\Ab_1^T\yb^k)}_2^2 \Big\},\\
\xb_2^{k+1} &{\!\!\!\!\!\!\!} := \mathrm{arg}\displaystyle\min_{\xb_2\in\Xc_2}\Big\{ f_2(\xb_2)  + \frac{\eta_k}{2\alpha_{2k}}\norm{\xb_2- (\mathbf{g}_2^k + \eta_k^{-1}(\Ab_2^T\yb^k))}_2^2 \Big\},
\end{array}\right.
\end{equation}
where $\kappa_k := \gamma_{k+1}+\rho_k$.
$\mathbf{g}_1^k := \xb_1^k - \alpha_{1k}\Ab_1^T(\Ab_1\xb_1^k + \Ab_2\xb_2^k - \bb) - \alpha_{0k}\Ab_1^T\Ab_1(\xb_1^k - \bar{\xb}_1^c)$, $\mathbf{g}_2^k := \xb_2^k - \alpha_{2k}\Ab_2^T(\Ab_1\xb_1^{k+1} + \Ab_2\xb_2^k - \bb)$, and step size $\alpha_{1k}$, $\alpha_{0k}$ and $\alpha_{2k}$ are chosen from gradient methods \cite{Yang2011}.

In \cite{He2012a,He2012} the authors proved the convergence of the standard $\mathrm{ADMM}$ algorithm at the rate of $\mathcal{O}(1/k)$ but in the sense of Auslender's gap function and requires the boundedness of both the primal and dual feasible sets.
In \cite{Ouyang2014} the authors considered other variant of ADMM, which requires the Lipschitz gradient assumption and still obtained the $\mathcal{O}(1/k)$ convergence rate both on the objective values $f(\xb^k) - f^{*}$ and the feasibility gap.
Other variants of ADMM can be found, e.g., in \cite{Deng2012,Ouyang2013,Wang2013a} and the references quoted therein, which were applied to stochastic cases or using different set of assumptions.

\subsection{Connections to proximal-based decomposition method}
If we set $\xb^k_c \equiv \bar{\xb}^{k-1}$ for $k\geq 1$ in our \eqref{eq:pd_scheme_2d} scheme, then the resulting scheme closely relates to the proximal-based decomposition method ($\mathrm{PBDM}$) studied in \cite{Chen1994,Shefi2014}.
Indeed, the main steps of $\mathrm{PBDM}$ can be expressed as follows:
\begin{equation}\label{eq:proximal_based_alg}
\left\{\begin{array}{ll}
\hat{\yb}^k & {\!\!} := \bar{\yb}^k + \gamma_k^{-1}(\Ab\bar{\xb}^k - \bb),\\
\xb^{k+1}_1 &{\!\!} := \mathrm{arg}\displaystyle\min_{\xb_1\in\Xc_1}\set{f_1(\xb_1) \!+\! (\hat{\yb}^k)^T\Ab_1\xb_1 + (\gamma_k/2)\Vert\xb_1 - \xb^k_1\Vert^2_2},\\
\xb^{k+1}_2 &{\!\!} := \mathrm{arg}\displaystyle\min_{\xb_2\in\Xc_2}\set{f_2(\xb_2) \!+\! (\hat{\yb}^k)^T\Ab_2\xb_2 + (\gamma_k/2)\Vert\xb_2 - \xb^k_2\Vert_2^2},\\
\bar{\yb}^{k+1} &{\!\!} := \bar{\yb}^k + \gamma^{-1}_k(\Ab_1\xb_1^{k+1} + \Ab_2\xb^{k+1}_2 - \bb).
\end{array}\right.
\end{equation}
Clearly, this method looks very similar to \eqref{eq:pd_scheme_2d}, where it has two dual steps and one primal step. 
Here, \eqref{eq:proximal_based_alg} uses only one parameter $\gamma_k$, $d_b(\xb,\hat{\xb}) := \frac{1}{2}\norm{\xb - \hat{\xb}}_2^2$ the Euclidian distance and $\Sb \equiv\mathbb{I}$.
In \cite{Shefi2014} the authors prove the convergence of the scheme \eqref{eq:proximal_based_alg} in a joint criterion $f(\tilde{\xb}^k) - f^{\star} + r\Vert\Ab\tilde{\xb}^k - \bb\Vert_2 \leq \frac{1}{k+1}\left[c_1 + c_2\max_{\Vert\yb\Vert\leq r}\Vert\yb - \yb^0\Vert^2_2\right]$, where $\tilde{\xb}^k := \frac{1}{k}\sum_{j=0}^{k-1}\xb^{j+1}$ and $c_1$, $c_2$ and $r$ are given constants.
This result is very similar to the ones in \cite{He2012a,He2012} for ADMM, which combines the primal objective residual and the primal feasibility gap. 
However, since \eqref{eq:separable_convex} is constrained,  $f(\xb^k) - f^{\star}$ may take an arbitrarily negative value. Hence, the joint criterion does not imply the approximation of the primal objective residual and the primal  feasibility gap separately. 
Moreover, as indicated in \cite{Goldstein2013}, convergence guarantee in a joint criterion is not sufficient to ensure that primal-dual methods work well in practice. 
It is important to control algorithmic parameters to trade-off between the objective residual and the feasibility of the problem. 
In our case, we prove a separated criterion on the objective residual and the primal feasibility, which allows one to control the parameters in order to trade-off these quantities. 
At the same time, our methods still exploit $p$-decomposability with parallel updates in the primal steps \eqref{eq:g_gamma} and \eqref{eq:prox_oper2}.

%%% 8. The implementation issues
\section{Implementation enhancements}\label{sec:impl_issues}
We discuss in this section how to enhance the practical performance of  Algorithm \ref{alg:pd_alg}. 
We observe that at least three steps in Algorithm \ref{alg:pd_alg} can be modified to enhance its practical performance: the choice of $\xb_c^k$, the update rule for parameters as well as the parallel and distributed implementation choices.

%% 8.1. The choice of the center point xc.
\subsection{The choice of proximal-point $\xb^k_c$ and Bregman distances} 
In \eqref{eq:pd_scheme_2p} and \eqref{eq:pd_scheme_2d}, we can adaptively choose the center point $\xb^k_c$ of the Bregman distance at each iteration. We propose two options:

\begin{itemize}
\item \textit{Proximal-point:} We can choose $\xb_c^k := \xb^{\star}_{\gamma_{k-1}}(\hat{\yb}^{k-1})$ for $k\geq 1$ in \eqref{eq:g_gamma}. 
This makes Algorithm \ref{alg:pd_alg} similar to the proximal-based decomposition algorithm in \cite{Chen1994}, which employs the proximal term $d_b(\cdot, \hat{\xb}^{\star}_{k-1})$ with the Bregman distance $d_b$.

\item\textit{ADMM variant:} If we choose $d_b$ to be the Euclidean distance, $\Sb$ and $\xb_c$ such that $d_b(\Sb\xb,\Sb\xb_c) := (1/2)\big[\Vert\Ab_1\xb_1 + \Ab_2(\xb_{\gamma_{k-1}}^{\star}(\bar{\yb}^{k-1}))_2 - \bb\Vert_2^2 + \Vert\Ab_1(\xb_{\gamma_k}^{\star}(\bar{\yb}^k))_1 + \Ab_2\xb_2 - \bb\Vert_2^2\big]$, then  \eqref{eq:pd_scheme_2d} becomes a new variant of $\mathrm{ADMM}$ as discussed in \eqref{eq:ADMM_variants}.
However, the convergence guarantee of this variant as well as the case where the center point $\xb_c^k$ changes remains unknown.

\item\textit{Preconditioned ADMM variant:} 
We can choose $\xb_c^k := (\mathbf{g}_1^k, \mathbf{g}_2^k)$, where  $\mathbf{g}_1^k$ and $\mathbf{g}_2^k$ are given in \eqref{eq:convergence_A2d}.
The step-size $\alpha_{1k}$ and $\alpha_{2k}$ can be taken as $\alpha_{1k} := \norm{\Ab_1}_2^{-2}$ and $\alpha_{2k} :=  \norm{\Ab_2}_2^{-2}$ or computed from the exact line-search rule. 
In this case, \eqref{eq:pd_scheme_2d} becomes a new variant of the preconditioned $\mathrm{ADMM}$ algorithm in \cite{Chambolle2011}.
\end{itemize}
In addition to the choice of $\xb_c^k$, we can also choose an appropriate prox-function $b_{\Xc}$ for the feasible set $\Xc$ in order to define the Bregman distance $d_b$. 
For instance, if $\Xc$ is a standard simplex, i.e., $\Xc := \big\{\xb\in\R^n_{+}~:~ \sum_{j=1}^n\xb_j = 1\big\}$, then the entropy prox-function $b_{\Xc}(\xb) := \xb^T\ln(\xb) + n$ becomes an appropriate choice.

%% 8.2. Guidance on tuning the parameters.
\subsection{Guidance on tuning the parameters}\label{subsec:tuning_params}
Since Algorithm \ref{alg:pd_alg} generates a sequence $\set{\bar{\wb}^k}_{k\geq 0}$ that decreases the smoothed gap function $G_{\gamma_k\beta_k}(\bar{\wb}^k)$ as required in Definition \ref{de:decreasing_sequence}. The actual decrease on the objective residual  is $f(\bar{\xb}^k) - f^{\star} \leq \gamma_k(D_{\Xc}^{\Sb} - \Psi_k/\gamma_k)$.
In practice, $D_k := D_{\Xc}^{\Sb} - \Psi_k/\gamma_k$ can be dramatically smaller than $D_{\Xc}^{\Sb}$ in the early iterations. This implies that increasing $\gamma_k$ in the early iterations might improve practical performance.  

Our strategy is based on the following observations.
 If $\gamma_k$ increases, then $\tau_k$ also increases. Consequently, $\beta_k$ decreases. 
 Since $\beta_k$ measures the primal feasibility gap $\mathcal{F}_k := \norm{\Ab\bar{\xb}^k - \bb}_2$ due to Lemma \ref{le:excessive_gap_aug_Lag_func}, 
 we only increase $\gamma_k$ if the feasibility gap $\mathcal{F}_k$ is relatively high. For instance, in the case  $\xb_c^k := (\mathbf{g}_1^k, \mathbf{g}_2^k)$, we can compute the dual feasibility gap as $\mathcal{H}_k := \gamma_k\norm{\Ab_1^T\Ab_2((\hat{\xb}^{\star}_{k+1})_2 - (\hat{\xb}^{\star}_k)_2)}$. Then, if $\mathcal{F}_k \geq s\mathcal{H}_k$ for some $s > 0$, we increase $\gamma_{k+1} := (1-c\tau_k)\gamma_k$ for some $c < 0$. In our implementation, we suggest the value $c=1.05\tau_k^{-1}$ as a default option. 

We can also decrease the parameter $\gamma_k$ in \eqref{eq:pd_scheme_2d} by $\gamma_{k+1} := (1-c_k\tau_k)\gamma_k$, where  $c_k := d_b(\Sb\xb_{\gamma_k}^{\star}(\hat{\yb}^k),\Sb\xb_c)/D^{\Sb}_{\Xc} \in [0, 1]$ after updating the vector $(\bar{\xb}^{k+1},\bar{\yb}^{k+1})$ in \eqref{eq:pd_scheme_2d} if we know a priori an upper bound estimate for $D_{\Xc}^{\Sb}$.

% 7.2. Parallel and distributed implementation
\subsection{Parallel and distributed implementation}\label{subsec:parallel_impl}
\begin{figure*}[!ht]
%\vskip-0.5cm
\begin{center}
\includegraphics[scale=0.62]{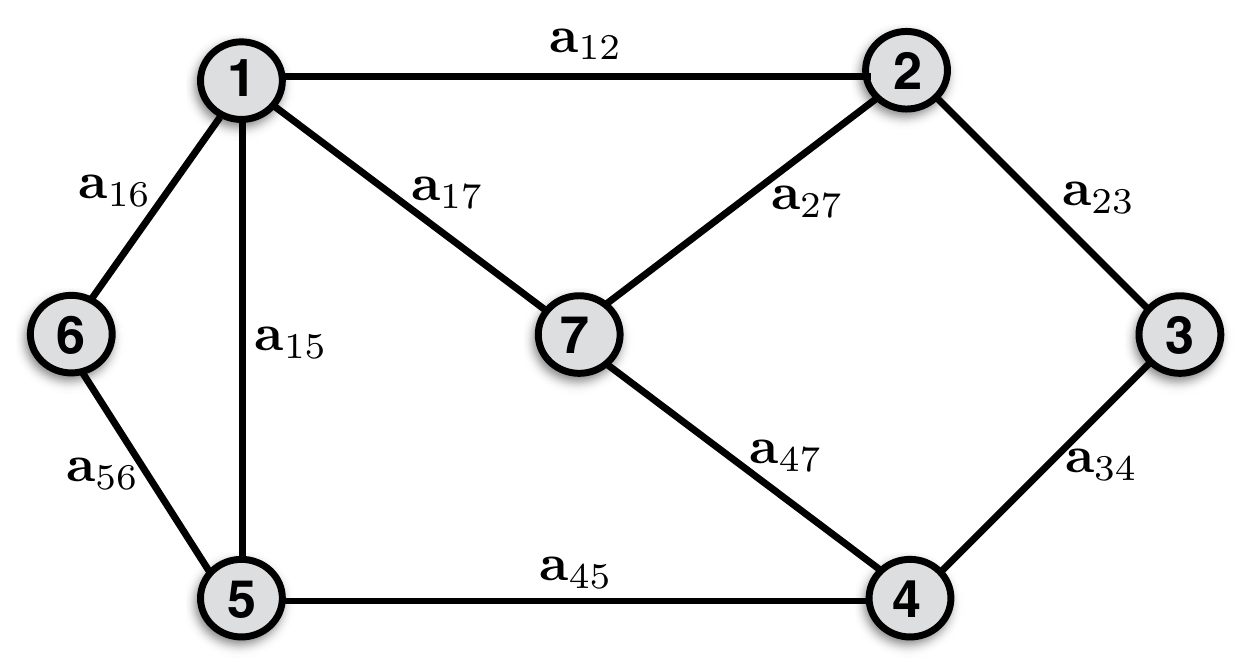}
\caption{A graph representing the structure of problem \eqref{eq:separable_convex} in a separable case.}\label{fig:graph}
\end{center}
\vskip-0.6cm
\end{figure*}

Suppose that $f$ and $\Xc$ are both separable as defined in \eqref{eq:separable_structure}, where each objective component $f_i$ and feasible set $\Xc_i$ correspond to the subsystem $i$ ($i=1,\dots, p$) of a large-scale network represented by a graph as illustrated in Figure \ref{fig:graph}.
The variable $\xb_i$ represents the unknown parameters of the subsystem $i$, and $\xb_i\in\Xc_i$ is its local constraint. Each subsystem $i$ communicates with its neighbors $j$ by asking the information from them via communication links $(i,j)$. Let  $\mathbf{a}_{ij}$ be the information the subsystem $i$ requests from its neighbor $j$ extracted from the neighbor's variable $\xb_j$. In this case, the information requested from all neighbors needs to be constrained by $\bb_i$, which leads to $\sum_{j\in \mathcal{N}_i}\mathbf{a}_{ij}\xb_j = \bb_i$, where $\mathcal{N}_i$ denotes all the neighbors of the subsystem $i$, for $i=1,\dots, p$. We note that each subsystem $i$ can have more than one links, the number of links leads to the number of coupling constraints.
To this end, one can reformulate a convex optimization problem over this network into a constrained problem of the form \eqref{eq:separable_convex} with separable objective, coupling constraints and separable local constraints.

Now, we assume that each $\Xc_i$ engages to a Bregman distance $d_{\Xc_i}$ with the convexity parameter $\sigma_i > 0$. 
We also choose either $\Sb := \norm{\Ab}_2\mathbb{I}$, $\Sb := \mathrm{diag}(A_1,\cdots, A_p)$ or $\Sb := \mathrm{diag}(\norm{A_1}_2\mathbb{I}_1,\dots, \norm{A_p}_2\mathbb{I}_p)$.
In this case, the Bregman distance of $\Xc$ becomes $d_{\Xc}(\xb, \xb_c) := \sum_{i=1}^pd_{\Xc_i}(\xb_i, \xb_{ic})$, where the strong convexity parameter of $d_{\Xc}$ is $\sigma_{d} := \min_{1\leq i\leq p}\sigma_i$.
The main step of Algorithm \ref{alg:pd_alg} is Step 5, where we need to perform the primal dual scheme \eqref{eq:pd_scheme_2p} or \eqref{eq:pd_scheme_2d}. 
We show how to implement these steps in a parallel and distributed manner based on the graph structure shown in Figure \ref{fig:graph}.

\paragraph{\textbf{Computation:}}
\vskip-0.08cm
The primal step in \eqref{eq:pd_scheme_2p} or \eqref{eq:pd_scheme_2d} requires to solve \eqref{eq:g_gamma} and \eqref{eq:prox_oper2}.
By the separability of $f$ and $\Xc$, \eqref{eq:g_gamma} can be solved \textbf{in parallel}. 
More precisely, each subsystem $i$ needs to estimate its local variable $\xb^k_i$ independently by solving a subproblem of the form:
\begin{equation*}
\xb_i^{k} := \mathrm{arg}\min_{\xb_i\in\Xc_i}\Big\{ f_i(\xb_i) + (\yb^k)^T\Ab_i\xb_i + \gamma d_{\Xc_i}(\xb_i, \xb_{ic}^k)\Big\}, ~~i=1,\dots, p,
\end{equation*}
where $\yb^k$ is a local copy of the Lagrange multiplier at the iteration $k$ for the subsystem $i$.

The dual step is updated as $\yb^{k+1} := \yb^k + \rho_k(\Ab\xb^k - \bb)$, where $\rho_k > 0$ is a given step size.
Here, each subsystem $i$ updates its local copy of the multiplier
\begin{equation*}
\yb_i^{k+1} := \yb_i^k + \rho_k\big(\sum_{j\in\mathcal{N}_i}\mathbf{a}_{ij}\xb_{ij} - \bb_i\big), ~~i=1\dots, p,
\end{equation*}
and sends this sub-vector to its neighbors to compute $\Ab_i^T\yb^{k+1}$ for the next iteration.

\paragraph{\textbf{Communication:}}
\vskip-0.08cm
At each iteration $k$, each subsystem $i$ requests the information from its neighbors to form the feasibility gap $\sum_{j\in\mathcal{N}_i}\mathbf{a}_{ij}\xb_{ij} - \bb_i$ and then updates $\yb^{k+1}_i$. This multiplier sub-vector is  then sent to the subsystem's  neighbors.

\paragraph{\textbf{Memory storage:}} 
Along with the local variable $\xb_i$ and the feasible set $\Xb_i$, each subsystem $i$ needs to store a copy of the dual variable $\yb^k$ and a part of  coefficient matrix $\Ab$ that represents the links to its neighbors, i.e., $\mathbf{a}_{ij}$ for $j\in\mathcal{N}_i$.

\paragraph{\textbf{Consensus and asynchronous operation:}} 
Note that our feasibility guarantees can be used to show the ``consensus'' of the distributed system with a corresponding rate when the communication graph is known \cite{Boyd2011}. 
Intriguingly, given that algorithms are tolerant to approximate proximal operators, we might expect them to also tolerate small levels of asynchronousity. 
Theoretical characterization of this important variant is left for future work. 

%% Extend to the inequality constraints and nonlinear constraints.
\subsection{Extension to inequality constraints}\label{sec:extensions}
The theory presented in the previous sections can be extended to solve convex optimization problems with linear inequality constraints of the form:
\begin{equation}\label{eq:ineq_cvx_prob}
f^{\star} := \displaystyle\min_{\xb\in\R^n} \set{ f(\xb) ~:~ \Ab\xb \leq \bb, ~\xb\in\Xc},
\end{equation}
where $f$, $\Xc$, $\Ab$ and $\bb$ are defined as in \eqref{eq:separable_convex}.

 A simple way to process \eqref{eq:ineq_cvx_prob} is using a slack variable $\sb \in\R^{m}_{+}$ such that $\Ab\xb + \sb = \bb$ and $\zb = (\xb, \sb)$ as the new variable. 
 Then we can transform \eqref{eq:ineq_cvx_prob} into \eqref{eq:separable_convex} with respect to the new variable $\zb$.
 
We can also process \eqref{eq:ineq_cvx_prob} by modifying the dual steps for updating $\hat{\yb}^k$, $\yb_{\beta_k}^{*}(\bar{\xb}^k)$ and $\bar{\yb}^{k+1}$ in both schemes \eqref{eq:pd_scheme_2p} and \eqref{eq:pd_scheme_2d}.
More precisely, we update these vectors as follows:
\begin{equation*}
\hat{\yb}^k := \big[\beta_{k+1}^{-1}(\Ab\bar{\xb}^k - \bb)\big]_{+}, ~\yb_{\beta_k}^{*}(\bar{\xb}^k) := \big[\beta_k^{-1}(\Ab\bar{\xb}^k - \bb)\big]_{+},
\end{equation*}
and 
\begin{equation*}
\bar{\yb}^{k+1} = \big[\hat{\yb}^k + (\gamma_{k+1}/\bar{L}^g)\big(\Ab\xb^{\star}_{\gamma_{k+1}}(\hat{\yb}^k) - \bb\big) \big]_{+},
\end{equation*}
where $[\cdot]_{+} := \max\set{0, \cdot}$. Indeed, the conclusion of Theorem \ref{th:convergence_A1} remains valid for this new variant for solving \eqref{eq:ineq_cvx_prob}.

\section{Numerical illustrations}\label{sec:num_results}
In this section, we present numerical simulations on several well-studied applications from machine learning, signal and image processing, and compressive sensing.
The numerical simulations  are performed using MATLAB R2012b, running on a Mac OS. i7 with 2.6Ghz and 16Gb RAM.  
We choose the Euclidean distance $d_b(\xb, \xb_c) := (1/2)\norm{\xb - \xb_c}^2$ in all test cases. 
We terminate Algorithm \ref{alg:pd_alg} if both primal feasibility gap 
\begin{equation*}
\mathcal{F}^r_k := \Vert\Ab\bar{\xb}^k - \bb\Vert_2/\max\set{1, \norm{\bb}_2} \leq \varepsilon_f, ~~\textrm{and}~~\Vert\bar{\xb}^{k+1} - \bar{\xb}^k\Vert_2/\max\big\{1, \Vert\bar{\xb}^k\Vert_2\big\} \leq\varepsilon_x,
\end{equation*}
for given default  tolerances $\varepsilon_f = 10^{-6}$ and $\varepsilon_x = 10^{-6}$ unless stated otherwise. 

%% 5.1. Convergence vs. worst-case complexity bound.
\subsection{Actual performance vs. theoretical bounds}
We demonstrate the empirical performance of the four variants of Algorithm \ref{alg:pd_alg} with respect to its theoretical bounds via a  basic non-overlapping sparse-group basis pursuit problem: 
\begin{equation}\label{eq:group_basis_pursuit}
\min_{\xb\in [\mathbf{l},\mathbf{u}]\subseteq \R^n} \sum_{i=1}^{n_g}w_i\norm{\xb_{g_i}}_2, ~~\mathrm{s.t.}~ \Ab\xb = \bb, 
\end{equation}
where $[\mathbf{l}, \mathbf{u}]$ is a box constraint,  and $g_i$ and $w_i$'s are the group indices and weights, respectively.

In this test, we choose $\xb^c = \mathbf{0}\in [\mathbf{l}, \mathbf{u}]$ and $d_b(\xb, \xb_c) := (1/2)\norm{\xb - \xb^c}^2$. We then evaluate $D_{\Xc}$ numerically, given $\Xc := [\mathbf{l}, \mathbf{u}]$. We estimate $D_{\mathcal{Y}^{\star}}$ and $f^{\star}$ by solving \eqref{eq:group_basis_pursuit} with an interior-point solver (SDPT3) \cite{Toh2010} up to accuracy $10^{-8}$. 
In the $(\mathrm{2P1D})$ scheme, we set $\gamma_0 = \beta_0 = \sqrt{\bar{L}_g}$, while, in the $(\mathrm{1P2D})$ scheme, we set $\gamma_0 := \frac{2\sqrt{2}\norm{\Ab}}{K+1}$ with $K:= 10^4$ and generate the theoretical bounds defined in Theorem \ref{th:convergence_A1}. 

We test the performance of the four variants using a synthetic sparse recovery problem, where $n = 1024$, $m = \lfloor n/3\rfloor = 341$, $n_g = \lfloor n/8\rfloor = 128$, and $\xb^{\natural}$ is a $\lfloor n_g/8\rfloor$-sparse vector. We set  $\mathbf{l} := \min(\xb^{\natural})$ and $\mathbf{u} := \max(\xb^{\natural})$. Matrix $\mathbf{A}$ are generated randomly from the iid standard Gaussian distribution and $\bb := \Ab\xb^{\natural}$. The group indices $g_i$ is also generated randomly for $i=1,\cdots, n_g$.

\paragraph{\textbf{Bregman smoothing case:}} Figure \ref{fig:theoretical_bound} shows the empirical performance of two variants: $(\mathrm{2P1D})$ and $(\mathrm{1P2D})$ of Algorithm \ref{alg:pd_alg}, where theoretical bounds are computed from Theorem \ref{th:convergence_A1}. 
\begin{figure*}[!ht]
\begin{center}
\includegraphics[width=12.5cm, height=8.3cm]{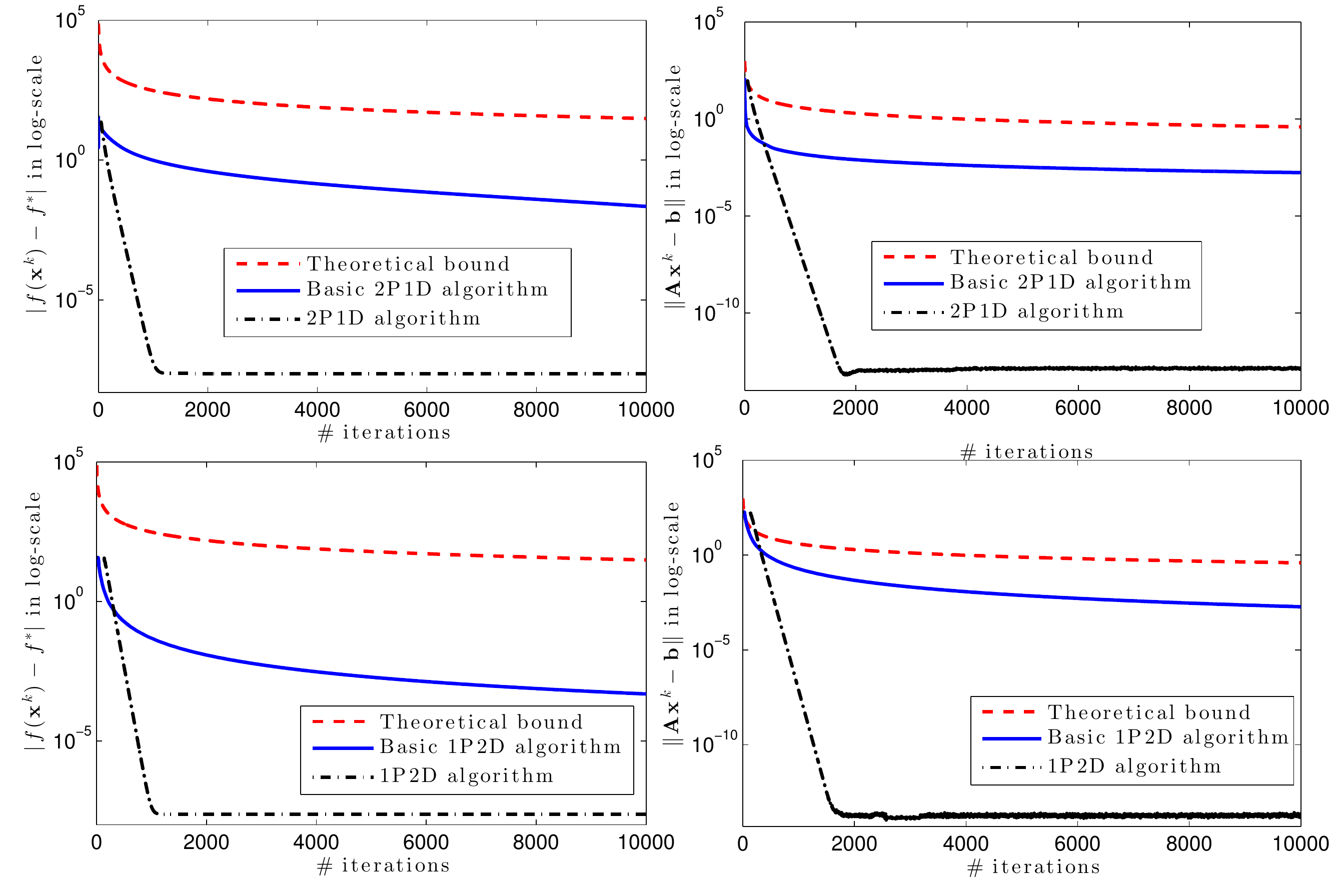}
\end{center}
\vskip-0.4cm
\caption{\footnotesize Actual performance vs. theoretical bounds of Algorithm \ref{alg:pd_alg} using Bregman smoother.}\label{fig:theoretical_bound}
%\vskip-0.4cm
\end{figure*} 
The basic algorithm refers to the case where $\xb_c$ is fixed and the parameters are not tuned. Hence, the iterations of the basic $(\mathrm{1P2D})$ use only 1 proximal calculation and applies $\Ab$ and $\Ab^T$ once each, and the iterations of the basic $(\mathrm{2P1D})$ use 2 proximal calculations and applies $\Ab$ twice and $\Ab^T$ once. In contrast, $(\mathrm{2P1D})$ and $(\mathrm{1P2D})$ variants whose iterations require one more application of $\Ab^T$ for adaptive parameter updates. 

It is clear from Figure~\ref{fig:theoretical_bound} that the empirical performance of the basic variants roughly follows the $\mathcal{O}(1/k)$ convergence rate both in terms of objective residual $\vert f(\bar{\xb}^k) - f^{\star}\vert$ and the feasibility gap $\norm{\Ab\bar{\xb}^k - \bb}_2$. The deviations from the bound are due to the increasing sparsity of the iterates, which improves empirical convergence.  With a kick-factor of $c_k=-0.02/\tau_k$ and adaptive proximal-center $\xb_c^k$ enhancements as suggested in Section \ref{sec:impl_issues}, the tuned $(\mathrm{2P1D})$ and $(\mathrm{1P2D})$ variants significantly outperform theoretical predictions.
Indeed, they approach the optimal solution up to $10^{-13}$ accuracy, i.e. $\norm{\bar{\xb}^k - \xb^{\star}} \leq 10^{-13}$ after only a few hundreds of iterations.

\paragraph{\textbf{Augmented Lagrangian smoothing case:}}
Similarly, Figure \ref{fig:theoretical_bound} illustrates the actual performance vs.\ the theoretical bounds $\mathcal{O}(1/k^2)$ by using augmented Lagrangian smoothing techniques. 
\begin{figure*}[!ht]
\begin{center}
\includegraphics[width=12.5cm, height=8.3cm]{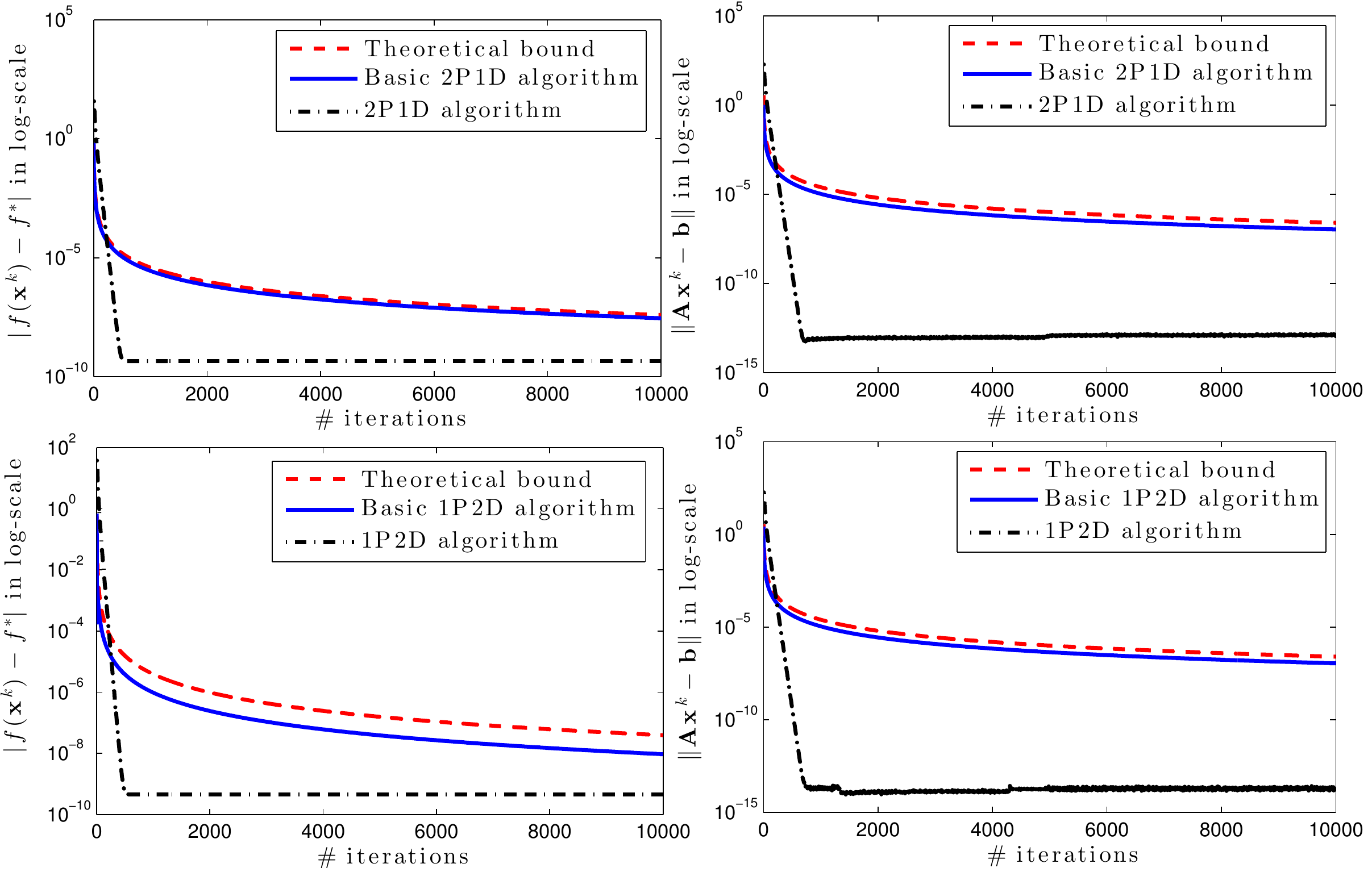}
\end{center}
\vskip-0.4cm
\caption{\footnotesize Actual performance vs. theoretical bounds of Algorithm \ref{alg:pd_alg} for augmented Lagrangian smoother.}\label{fig:theoretical_bound}
%\vskip-0.8cm
\end{figure*} 
Here, we solve the subproblems \eqref{eq:aug_dual_func} and \eqref{eq:prox_oper2_inexact} by using FISTA \cite{Beck2009}.
Since, we can not exactly estimate the true solution of the subproblems \eqref{eq:aug_dual_func} and \eqref{eq:prox_oper2_inexact}, we solve these problems up to at least the accuracy $\delta_0^2 = 10^{-8}$ as suggested by Theorem \ref{th:convergence_A1}.

In this case, the theoretical bounds and the actual performance of the basis variants are very close to each other both in terms of the objective residual $\vert f(\bar{\xb}^k) - f^{\star}\vert$ as well as the primal feasibility gap $\Vert\Ab\bar{\xb}^k - \bb\Vert_2$.
When the parameter $\gamma_k$ is updated, the algorithms exhibit a  better performance. 

\paragraph{\textbf{Strongly convex case:}} We demonstrate the theoretical bounds for the strongly convex case via the elastic net:
\begin{equation}\label{eq:elastic_net}
\min_{\xb}\set{ \Vert\xb\Vert_1 + (\sigma/2)\Vert\xb\Vert_2^2} ~~\mathrm{s.t.}~ \Ab\xb = \bb,
\end{equation}
where $\sigma > 0$ is a given constant, and other parameters are selected as in \eqref{eq:group_basis_pursuit}. 
The data of this test is also generated randomly as for \eqref{eq:group_basis_pursuit}, where $n := 2000$, $m = 700$ and $\xb^{\natural}$ is $100$-sparse.

\begin{figure*}[!ht]
\begin{center}
\vskip-0.4cm
\includegraphics[width=12.5cm, height=8.3cm]{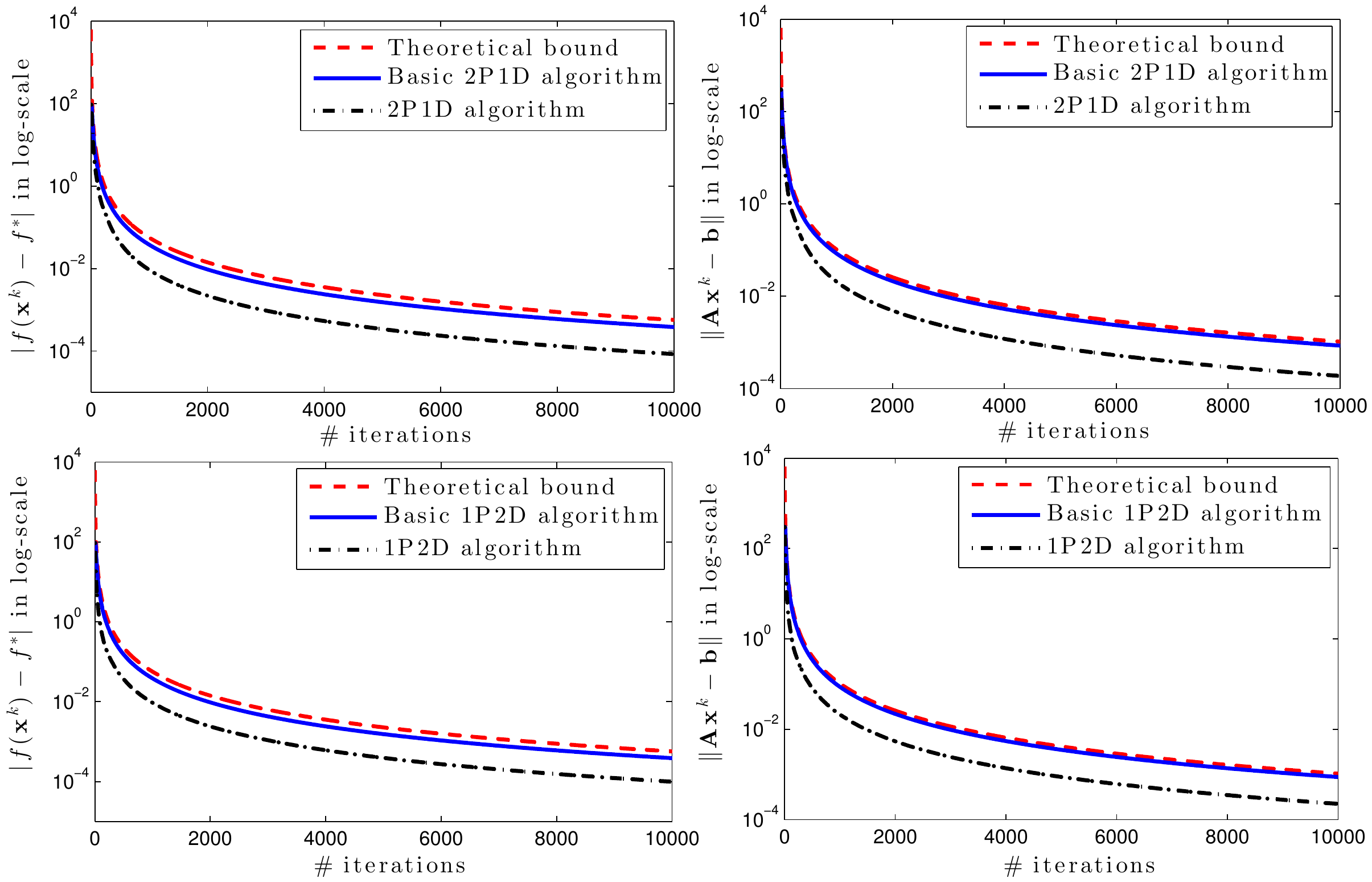}
\end{center}
\vskip-0.4cm
\caption{Actual performance vs.\ theoretical bounds for strongly convex case.}\label{fig:theoretical_bound3}
%\vskip-0.4cm
\end{figure*}

We test Algorithm \ref{alg:pd_alg} using both $(\mathrm{2P1D}_{\sigma})$ and $(\mathrm{1P2D}_{\sigma})$ to solve \eqref{eq:elastic_net} with $\sigma := 0.1$. The results are plotted in Figure \ref{fig:theoretical_bound3} for both $\vert f(\bar{\xb}^k) - f^{\star}\vert$ and $\Vert\Ab\bar{\xb}^k-\bb\Vert_2$, respectively after $K := 10^4$ iterations. The configuration of the basic variants are as before whereas the enhanced versions use a backtracking linesearch procedure to determine an approximation $L_k$ for the Lipschitz constant $L^g_f$. The iterates converge better than the theoretical rate (see the appendix). 

We obtain the final relative errors $(\vert f(\bar{\xb}^k) - f^{\star}\vert/\vert f^{\star}\vert, \Vert\Ab\bar{\xb}^k-\bb\Vert_2/\norm{\bb}_2)$ for both cases are $(4.0376, 2.8294)\times 10^{-6}$ and $(4.0744, 2.9064)\times 10^{-6}$, respectively. These values are $(0.8900, 0.6237)\times 10^{-6}$ and $(1.0462, 0.7400)\times 10^{-6}$, respectively, in the line-search variants, which are approximately $4$ times smaller than in the basic ones.
The relative recovery error $\Vert\bar{\xb}^k-\xb^{\star}\Vert_2/\norm{\xb^{\star}}_2$ is also $3.228\times 10^{-7}$ and $3.753\times 10^{-7}$, respectively. 
We also observe that after $642$ (reps., $691$) iterations, both algorithms reach the accuracy $\Vert\bar{\xb}^k-\xb^{\star}\Vert_2 \leq 10^{-2}$, and after $2034$ (resp., $2193$) iterations, which corresponds to an approximate relative error of $10^{-3}$. 

\begin{figure*}[!ht]
\begin{center}
%\vskip-0.5cm
%\includegraphics[scale=0.45]{figs/theo_bound_scvx_xb_20_06_2014}
\includegraphics[width=12.5cm, height=4.5cm]{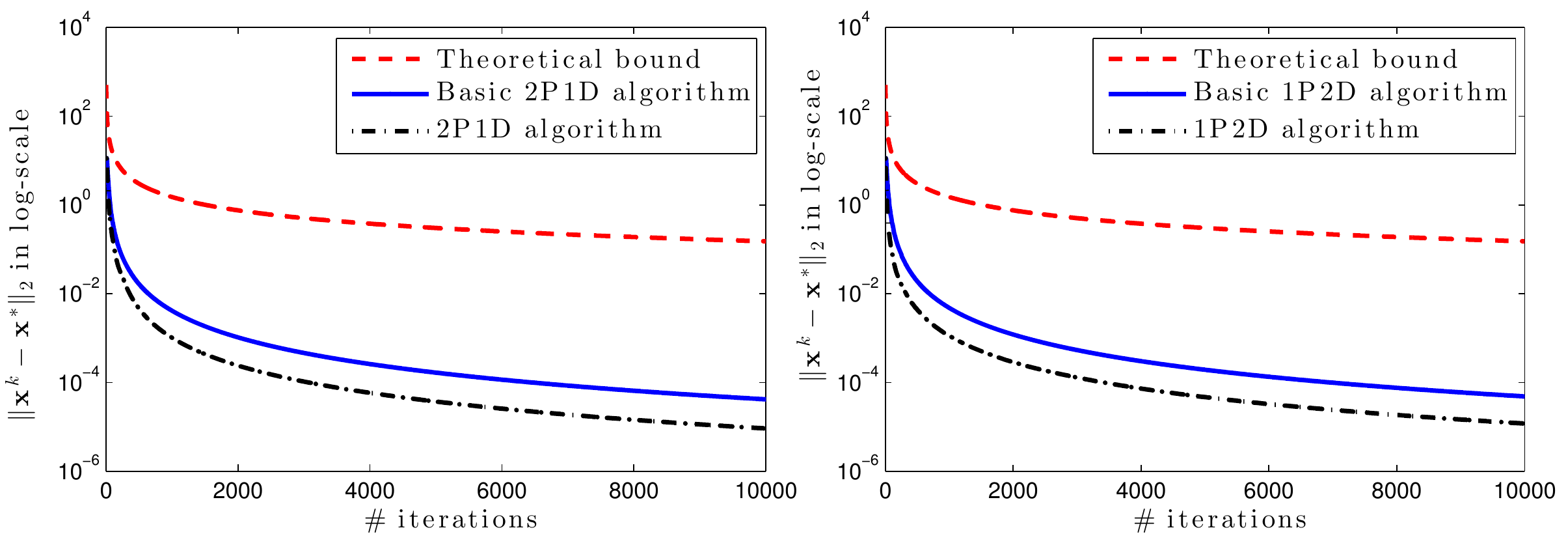}
\end{center}
\vskip-0.4cm
\caption{Actual performance vs.\ theoretical bound for strongly convex case (iterative sequence).}\label{fig:theoretical_bound3_xk}
%\vskip-0.4cm
\end{figure*}

We also compute the practical values of $\set{\Vert\bar{\xb}^k - \xb^{\star}\Vert_2}_{k\geq 0}$ and its theoretical bound shown in Corollary \ref{co:strong_convex_convergence} for \eqref{eq:elastic_net}.
The convergence of  $\set{\Vert\bar{\xb}^k - \xb^{\star}\Vert_2}_{k\geq 0}$ and its theoretical bound is plotted in Figure \ref{fig:theoretical_bound3_xk} for both algorithms: $(\mathrm{1P2D}_{\sigma})$ and $(\mathrm{2P1D}_{\sigma})$, and their line-search variants, respectively.

We can see that the theoretical bound given by Corollary \ref{co:strong_convex_convergence} is far from the actual performance. This is clearly observed due to a rough estimation of the upper bound. 
The line-search variants takes less iterations than the basic ones, but require additional computations for the line-search procedure, which makes them in the end slower.

\paragraph{\textbf{A new variant of preconditioned ADMM:}}
Finally, we verify the theoretical justification of the new PADMM variant given in \eqref{eq:linearized_alternating_method}. The same test can be done for  the new ADMM variant.

We use again the group basis pursuit problem \eqref{eq:group_basis_pursuit} by reformulating it into the following form:
\begin{equation}\label{eq:group_basis_pursuit2}
\min_{\xb\in [\mathbf{l},\mathbf{u}]\subseteq \R^n} \sum_{i=1}^{n_g}w_i\norm{\xb_{g_i}}_2 + \delta_{\set{0^m}}(\mathbf{r}), ~~\mathrm{s.t.}~ \Ab\xb + \mathbf{r} = \bb, ~\mathbf{r}\in [\underline{\mathbf{r}}, \bar{\mathbf{r}}],
\end{equation}
where $\delta_{\mathcal{S}}$ is the indicator function of the set $\mathcal{S}$, $\underline{\mathbf{r}}$ and  $\bar{\mathbf{r}}$ is computed from the bounds $\mathbf{l}$ and $\mathbf{u}$ of $\xb$ via the relation $\rb = -\Ab\xb - \bb$.

\begin{figure*}[!ht]
\begin{center}
%\vskip-0.4cm
%\includegraphics[scale=0.45]{figs/theo_padmm}
\includegraphics[width=12.5cm, height=4.5cm]{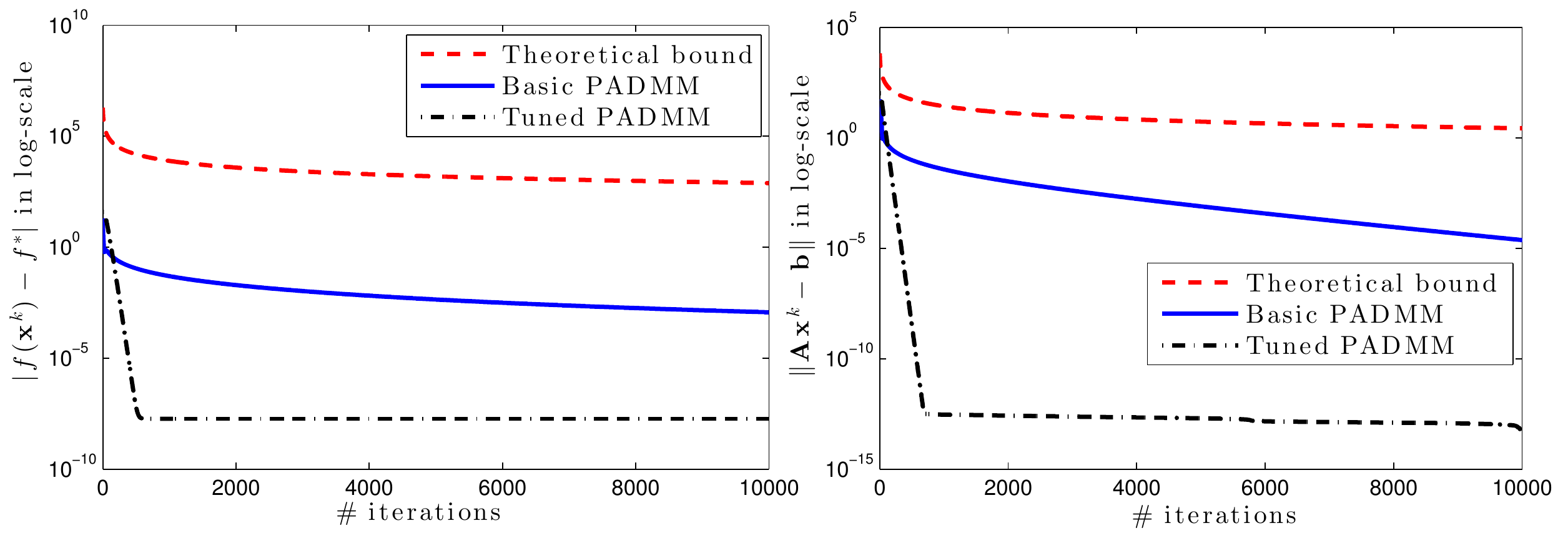}
\end{center}
\vskip-0.4cm
\caption{Actual performance vs.\ theoretical bound for PADMM.}\label{fig:theo_padmm}
%\vskip-0.2cm
\end{figure*}

We test the new variant of preconditioned ADMM (PADMM) and compare it with the tuned version, where we adaptively update the parameter $\gamma_k$ using the strategy in Section \ref{sec:impl_issues}.
In the basis PADMM variant, we fix $\gamma_0 := 2\sqrt{2}\Vert\Ab\Vert_2/(K+ 1)$, where $K = 10^4$.% as suggested by Corollary \ref{co:newPADMM}.
By using the same data as in the previous cases, we obtain the performance of this variant as shown in Figure \ref{fig:theo_padmm}.

As we can observe from Figure \ref{fig:theo_padmm} that, the basis PADMM variant relatively follows the curvature of the theoretical bounds, while the tuned variant reaches very high accuracy solution after few hundreds of iterations. This behavior is similar to the (1P2D) variant using Bregman smoother tested above.

%% 5.5. Performance with enhancements.
\subsection{Performance robustness.}
We demonstrate the performance robustness of our tuned $\mathrm{(1P2D)}$ variant by applying it to the following image deconvolution problem:
\begin{equation}\label{eq:image_deconv}
\min_{\xb: 0\le \xb\le 255} (1/2)\norm{\mathcal{B}(\xb) - \bb}_2^2 + \lambda \norm{\xb}_{\mathrm{TV}},
\end{equation}
where $\bb$ is a given blurry image with a known blur kernel $\mathcal{B}$, and  $\norm{\cdot}_{\mathrm{TV}}$ is the isotropic total variation norm and $\lambda > 0$ is a regularization parameter.  

As opposed to directly using the TV-norm proximal map, we simply use the linear mapping $\mathbf{D}$ of its norm operator $\Vert\xb\Vert_{\text{TV}} = \Vert\mathbf{D}\xb \Vert_1$ and introduce a slack variable $\rb = \mathbf{D}\xb$ to split \eqref{eq:image_deconv} into $\xb$ and $\rb$ variables with additional linear coupling constraint $\rb - \mathbf{D}\xb = 0$. Hence, we can reformulate  \eqref{eq:image_deconv} into \eqref{eq:separable_convex}, where $\rb \in \mathcal{R} := \set{\hat{\rb} ~|~ \hat{\rb} = \Db\xb, 0 \leq \xb \leq 255}$ is also bounded. 

We apply the $\mathrm{(1P2D)}$ variant of Algorithm \ref{alg:pd_alg} to solve the resulting problem and compare it with the ADMM solver implemented in \cite{Chan2011} since both algorithms have similar complexity per iteration.
We choose the center point as suggested in our practical enhancement guidelines, which leads to a \textit{new variant} of the standard ADMM method.
We test two cases:  without and with tuning based on our guidance. We choose the initial regularization parameters $\rho_0$ the same as the recent \emph{exact} ADMM solver suggests \cite{Chan2011}. 

Surprisingly, if we assume periodic boundary conditions for the TV-norm, then  \textrm{ADMM} can efficiently obtain accurate solutions to the subproblems in computing $\xb_1^k$ and $\xb_2^k$. The key idea is that the operator $ \mathbf{D}^T\mathbf{D} + \mathcal{B}^T\mathcal{B}$ is diagonalizable by the Fourier transform. Hence, the complexity per iteration in exact \textrm{ADMM} and $(\mathrm{1P2D})$ is approximately the same. Note however that our algorithm does not require periodic boundary conditions  to solve this class of problems, which may not be valid in other applications

\begin{figure*}[ht!]
\begin{center}
\centerline{\includegraphics[width = 15.5cm, height=4.0cm]{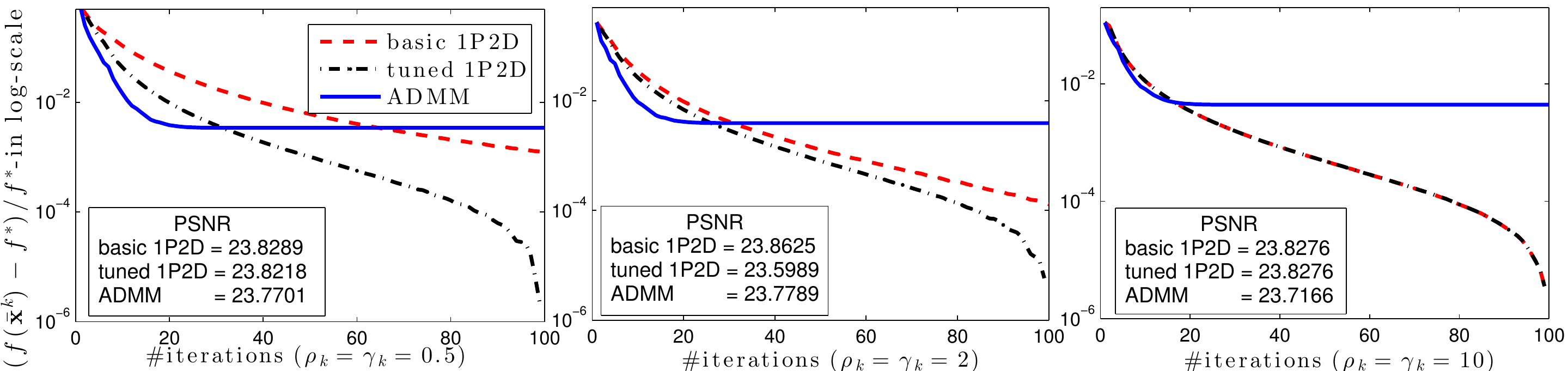}}
\caption{The performance of the augmented Lagrangian methods under different penalty parameters. }\label{fig:camera_man}
\end{center}
\vskip -0.4cm
\end{figure*}

Figure \ref{fig:camera_man}  illustrates the performance of $(\mathrm{1P2D})$ and  the \textrm{ADMM} code \cite{Chan2011} with different values of parameter $\gamma$ (resp., $\rho$ in the ADMM solver). Our test is based on the \texttt{camera\_man} image, with the regularization $\lambda = 0.01$ as done in \cite{Chan2011}. The suggested value for $\rho$ is $\rho = 2$ in \cite{Chan2011}. The exact ADMM  code  \cite{Chan2011} also uses a specific update rule for the penalty parameter, which is different from ours.  Figure \ref{fig:camera_man} shows the convergence of three algorithms wrt.\ three values of $\gamma$ (respectively, $\rho$) after $100$ iterations. We can see that ADMM decreases quickly first but then does not move, while $(\mathrm{1P2D})$ continues to descend on the objective function. 

We note that the ADMM solver is sensitive to the choice of $\rho$. For any value of $\rho$, if we run up to $200$ iterations then the exact  ADMM algorithm diverges, which is due to their aggressive update rule on the penalty parameter.  

%% 5.4. The effect of inexact computation.
\subsection{Inexact computations.}
In this test, we study the empirical impact of inexact proximal operator calculations to the performance of Algorithm \ref{alg:pd_alg}. Again, we choose the $(\mathrm{1P2D})$ variant, which has similar complexity per iteration as preconditioned ADMM \cite{Chambolle2011}. 
For this, we use a Schatten norm based regularizer on a Poisson likelihood data model:
\begin{equation}\label{eq:poisson}
\min_{\xb\in \mathcal{X}}(\mathcal{B}(\xb))^T\mathbf{1} - \sum_{i=1}^m\mathbf{c}_i\log((\mathcal{B}(\xb))_i + \bb) + \lambda \norm{\xb}_S,
\end{equation}
where $\mathcal{X} := \R^n_{+}$, $\mathbf{c}$ is a given photon count vector in $\mathbb{Z}^m$, $\bb$ is the background intensity, $\lambda > 0$ is a chosen regularization parameter, and $\mathcal{B}$ is a blur kernel. This likelihood model is quite common in scientific imaging problems.

The work in \cite{Lefkimmiatis2013b} proposed a norm based on exploiting self-similarities within the images via $\norm{\xb}_S := \norm{\mathrm{mat}(\mathcal{H}(\xb))}_{\star}$, which is the Schatten-norm of a matrix $\mathrm{mat}(\mathcal{H}(\xb))$ for a suitably chosen linear operator $\mathcal{H}$.
Since the proximal operator regarding the second term $f_2(\xb) := \lambda \norm{\xb}_S + \delta_{\mathcal{X}}(\xb)$, where $\delta_{\Xc}$ is the indicator  of $\Xc$, does not have a closed form, we need to iteratively compute it. 

The resulting inexact computation affects the performance of optimization algorithms. Here, we compare our new PADMM variant of Algorithm \ref{alg:pd_alg} (called tuned 1P2D) with  PADMM and  PADMM based on our tuning strategy in the enhancement paragraph as well as the exact ADMM solver provided by \cite{Lefkimmiatis2013b}. Here, the ADMM solver exploits boundary conditions and Fourier transform to invert $ I + \mathcal{B}^T\mathcal{B}$ for solving its subproblems. 
When $\bb$ is zero (i.e., there is no background), then the logarithmic term pose computational problems since its gradient is no longer Lipschitz. Fortunately, the proximal operator of the $\log$ function can be efficiently calculated. 
\begin{figure}[!ht]
\begin{center}
%\vskip-0.4cm
%\centerline{\includegraphics[scale=0.42]{figs/poisson_face_055_200}}
\centerline{\includegraphics[width = 15.5cm, height=4.0cm]{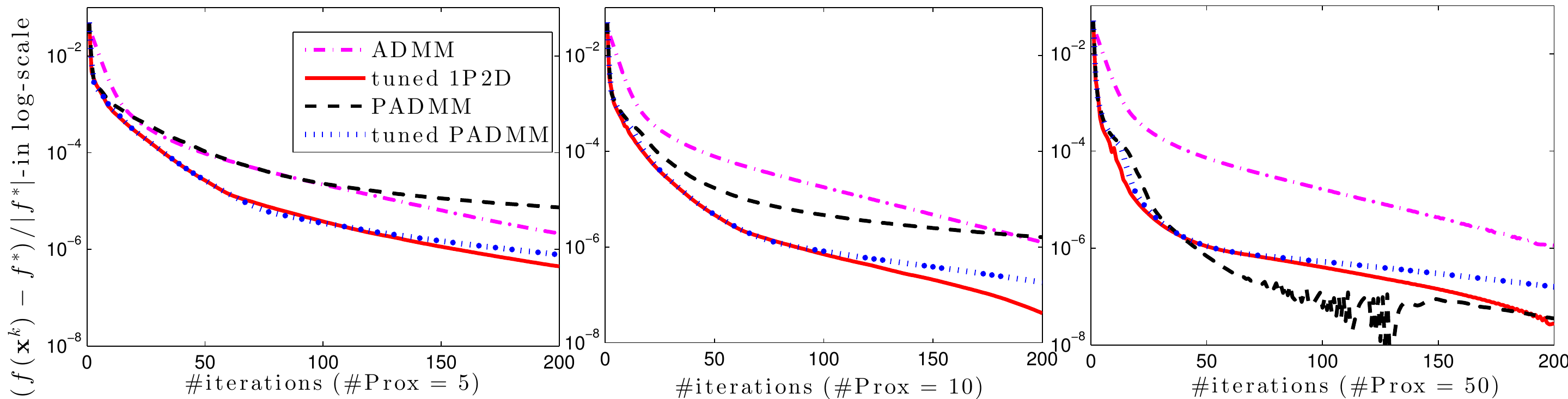}}
\caption{The performance of four algorithms on the Clown image \cite{Lefkimmiatis2013b}. }\label{fig:poisson1}
\vskip -0.4cm
\end{center}
\end{figure} 

We test these algorithms on the Clown image where we take the regularization parameter $\lambda = 0.055$ suggested in \cite{Lefkimmiatis2013b}. We use the \texttt{Denoise} solver in \cite{Lefkimmiatis2013b} to approximately compute the prox-operator of $f_2$ with inner iterations $\mathrm{nProx}$$=5,10,50$, where we can warm start each iteration using each algorithms current estimate. The exact ADMM solver is already implemented with penalty parameter updates. 

Figure \ref{fig:poisson1} illustrates that our tuned $(\mathrm{1P2D})$ solver and \textrm{PADMM} are quite robust to inexact prox calculations and outperform exact ADMM for a range of $\mathrm{nProx}$ values. Against intuition, we observe that PADMM exhibits numerical instability when $\mathrm{nProx}$ is highest. Overall, our algorithm provides the best time to reach an $\epsilon$-solution since doubling $\mathrm{nProx}$ roughly doubles the overall time. For instance, $\mathrm{nProx}=5$ and $200$ iterations roughly takes the same time as $\mathrm{nProx}=10$ and $100$ iterations, where our algorithm provides the best accuracy. 

In this setting, our solver and PADMM do not require periodic boundary conditions. When this assumption is removed, the subproblem are no longer dominated by just prox calculations. Then, we expect our algorithm obtain better timing performance due to its parallel updates. 

%% 5.2. Performance with enhancements.
\subsection{Additional comparisons with state-of-the-art.}
We compare our algorithms with existing state-of-the-art Matlab codes for solving five well-studied problems: standard basis pursuit, group-sparse basis pursuit, robust PCA, square-root LASSO and support vector machines with the Hinge loss.
While there are several software packages that can be used to solve these problems, we only select few of representatives which we find as the most efficient methods for corresponding problems.

%%%% 8.4.1. Standard basis pursuit.
\subsubsection{Standard basis pursuit.}
We consider the standard basis pursuit problem arising from compressive sensing \cite{Donoho2006}:
\begin{equation}\label{eq:standard_BP}
\min_{\xb\in\R^n} \norm{\xb}_1 ~~\mathrm{s.t.}~ \Ab\xb = \bb,
\end{equation}
where $\Ab\in\R^{m\times n}$ and $\bb\in\R^{m}$.

In this example, we compare our algorithms with \textrm{YALL1} \cite{Yang2011} and \textrm{SPGL$_1$} \cite{VanDenBerg2008} which are well-known solvers for the basis pursuit problem.
We use the data from the benchmark collection \texttt{Sparco} \cite{Sparco2007}. 
For \textrm{YALL1} and \textrm{SPGL$_1$}, we use the default settings and all the algorithms are terminated with the accuracy $10^{-6}$.
Within our methods, we run three algorithms: $(\mathrm{1P2D})$ via Bregman distance smoothing, inexact $\mathrm{1P2D(1)}$ (with only one FISTA iteration) and inexact $\mathrm{1P2D(5)}$ (with $5$ FISTA iterations) via augmented Lagrangian smoothing.
The two last algorithms are inexact variants of Algorithm \ref{alg:pd_alg} using the augmented Lagrangian smoother.
Table \ref{tbl:sparsco_benchmark} shows the problems selected from the \texttt{Sparco} test collection \cite{VanDenBerg2008}  that we use for our test. 
\begin{table*}[!ht]
\newcommand{\cell}[1]{{\!}#1{\!}}
\newcommand{\cellc}[1]{{\!}\textbf{\color{blue}#1}{\!}}
\begin{center}
%\vskip-0.3cm
\caption{The \texttt{Sparco} test problems used}\label{tbl:sparsco_benchmark}
\vskip0.25cm
\begin{tabular}{lccccl}\hline
Problems & ID & $m$ & $n$ & $\norm{\bb}_2$ & Operators \\ \hline
%\texttt{blocksig} & 2 & 1024 & 1024 & 7.9e+1 &  wavelet \\ \hline
%\texttt{cosspike} & 3 & 1024 & 1024 & 1.0e+2 & DCT \\ \hline
%\texttt{zsinspike} & 4 & 1024 & 2048 & 1.1e+2 & DCT \\ \hline
\texttt{gcosspike} & 5 & 300 & 2048 & 8.1e+1 & Gaussian ensemble, DCT \\ \hline
\texttt{p3poly} & 6 & 600 & 2560 & 2.2e+0 & Gaussian ensemble, wavelet \\ \hline
\texttt{sgnspike} & 7 & 600 & 2560 & 2.2e+0 & Gaussian ensemble \\ \hline
\texttt{zsgnspike} & 8 & 600 & 2560 & 2.9e+0 & Gaussian ensemble \\ \hline
\texttt{gausspike} & 11 & 256 & 1024 & 8.7e+1 & Gaussian ensemble \\ \hline
\texttt{srcsep1} & 401 & 29166 & 57344 & 2.2e+1 & windowed DCT \\ \hline
\texttt{srcsep2} & 402 & 29166 & 86016 & 2.3e+1 & windowed DCT \\ \hline
%\texttt{srcsep3} & 403 & 29166 & 86016 & 2.3e+1 & windowed DCT \\ \hline
\texttt{phantom1} & 501 & 629 & 4096 & 1.1e+1 & restricted FPT, wavelet \\ \hline
%\texttt{soccer1} & 601 & 3200 & 4096 & 5.5e+4 & binary ensemble, wavelet \\ \hline
%\texttt{yinyang} & 603 & 1024 & 4096 & 2.5e+1 & wavelet \\ \hline
\texttt{blurrycam} & 701 & 65536 & 65536 & 1.3e+2 & blurring, wavelet \\ \hline
\texttt{blurspike} & 702 & 16384 & 16384 & 2.2e+0 & blurring  \\ \hline
%\texttt{dcthdr} & 12 & 2000 & 8192 & 2.3e+3 & restricted DCT \\ \hline
%\texttt{finger} & 703 & 11013 & 125385 & 5.5e+4 & 2D curvelet \\ \hline
%\texttt{seismic} & 901 & 41472 & 480617 & 1.1e+2 & 2D curvelet \\ \hline
%\texttt{jitter} & 902 & 200 & 1000 & 4.7e+1 & DCT \\ \hline
%\texttt{spiketrn} & 903 & 1024 & 1024 & 5.7e+1 & 1D convolution \\ \hline
\end{tabular}
\end{center}
%\end{footnotesize}
\vskip-0.2cm
\end{table*}

%\textbf{QUOC: I would simply remove the coloring in the table}
The numerical results and performance information are reported in Table \ref{tbl:BP_test} for $10$ problems from Table \ref{tbl:sparsco_benchmark}.
Our algorithms and \textrm{YALL1} are still superior to  $\mathrm{SPGL}_1$ both in terms of number of iterations, matrix-vector multiplications and CPU time, while producing very similar final objective value $f(\xb^k)$ and the feasibility gap $\Vert\Ab\xb^k - \bb\Vert_2$.
\textrm{YALL1}  performs quite well compared to our methods in terms of timing. However, it fails for the last two problems (i.e., \texttt{blurrycam} and \texttt{bluspike}) due to their parameter update rules.
\begin{table*}[!ht]
%\vskip-0.2cm
\newcommand{\cell}[1]{{\!\!\!}#1{\!\!}}
\newcommand{\cellbf}[1]{{\!\!\!\!}\textbf{\color{blue}#1}{\!\!}}
%\begin{scriptsize}
%\begin{footnotesize}
\begin{center}
%\vskip-0.5cm
\caption{Comparison of the five algorithms: $(\mathrm{1P2D})$, \textrm{1P2D(1)}, \textrm{1P2D(5)}, \textrm{YALL1} and \textrm{SPGL$_1$}.}\label{tbl:BP_test}
\vskip0.25cm
\begin{tabular}{l|rrrrr|rrrrr}\hline
 & \cell{$\mathrm{1P2D}$} & \cell{\textrm{1P2D(1)}} & \cell{\textrm{1P2D(5)}} & \cell{\textrm{YALL1}} & \cell{$\mathrm{SPGL}_1$}  & \cell{$\mathrm{1P2D}$} & \cell{\textrm{1P2D(1)}} & \cell{\textrm{1P2D(5)}} & \cell{\textrm{YALL1}} & \cell{$\mathrm{SPGL}_1$} \\ \hline
\cell{Problems} & \multicolumn{5}{c}{\#Iterations} & \multicolumn{5}{c}{CPU time [s]} \\ \hline
%%%%
\cell{gcosspike} &\cell{330} & \cell{275} & \cell{274} & \cellbf{208} & \cell{1026} & \cell{0.87} & \cell{0.74} & \cell{2.16} & \cellbf{0.62} & \cell{2.53} \\ 
\cell{p3poly} &\cell{306} & \cell{100} & \cellbf{98} & \cell{252} & \cell{1775} & \cell{14.32} & \cellbf{4.81} & \cell{16.11} & \cell{13.34} & \cell{67.10} \\ 
\cell{sgnspike} &\cell{346} & \cell{157} & \cellbf{156} & \cell{178} & \cell{291} & \cell{0.96} & \cellbf{0.50} & \cell{1.48} & \cell{0.61} & \cell{1.01} \\ 
\cell{zsgnspike} &\cell{331} & \cell{307} & \cell{307} & \cellbf{152} & \cell{320} & \cell{1.59} & \cell{1.53} & \cell{4.65} & \cellbf{0.91} & \cell{1.87} \\ 
\cell{gausspike} &\cell{368} & \cell{320} & \cellbf{319} & \cell{170} & \cell{516} & \cell{0.31} & \cell{0.29} & \cell{0.67} & \cellbf{0.19} & \cell{0.52} \\ 
\cell{srcsep1} &\cell{380} & \cell{331} & \cellbf{330} & \cell{426} & \cell{1580} & \cell{22.65} & \cellbf{19.44} & \cell{67.88} & \cell{36.95} & \cell{119.60} \\ 
\cell{srcsep2} &\cell{376} & \cell{326} & \cellbf{325} & \cell{334} & \cell{1310} & \cell{34.64} & \cellbf{29.73} & \cell{102.88} & \cell{58.19} & \cell{155.12} \\ 
\cell{phantom1} &\cell{291} & \cell{285} & \cell{285} & \cellbf{166} & \cell{712} & \cell{1.16} & \cell{1.02} & \cell{2.70} & \cellbf{0.50} & \cell{2.42} \\ 
\cell{blurrycam} &\cell{1042} & \cell{3496} & \cellbf{569} & \cell{failed} & \cell{3629} & \cellbf{23.48} & \cell{72.97} & \cell{39.08} & \cell{failed} & \cell{152.96} \\ 
\cell{blurspike} &\cell{1255} & \cell{4191} & \cellbf{797} & \cell{failed} & \cell{2159} & \cellbf{5.83} & \cell{18.86} & \cell{10.14} & \cell{failed} & \cell{17.16} \\ 
%%%%%
\hline
\cell{Problems} & \multicolumn{5}{c}{\#$\Ab\xb$} & \multicolumn{5}{c}{ \#$\Ab^T\yb$} \\ \hline
%%%%
\cell{gcosspike} &\cell{332} & \cell{552} & \cell{1600} & \cellbf{312} & \cell{1815} & \cell{331} & \cellbf{276} & \cell{1325} & \cell{416} & \cell{1028} \\ 
\cell{p3poly} &\cell{308} & \cellbf{202} & \cell{590} & \cell{378} & \cell{3279} & \cell{307} & \cellbf{101} & \cell{491} & \cell{504} & \cell{1777} \\ 
\cell{sgnspike} &\cell{348} & \cell{316} & \cell{902} & \cellbf{267} & \cell{482} & \cell{347} & \cellbf{158} & \cell{745} & \cell{178} & \cell{293} \\ 
\cell{zsgnspike} &\cell{333} & \cell{616} & \cell{1766} & \cellbf{228} & \cell{557} & \cell{332} & \cell{308} & \cell{1458} & \cellbf{152} & \cell{322} \\ 
\cell{gausspike} &\cell{370} & \cell{642} & \cell{1840} & \cellbf{255} & \cell{858} & \cell{369} & \cell{321} & \cell{1520} & \cell{340} & \cell{518} \\ 
\cell{srcsep1} &\cellbf{382} & \cell{664} & \cell{1962} & \cell{639} & \cell{2639} & \cellbf{381} & \cell{332} & \cell{1631} & \cell{852} & \cell{1582} \\ 
\cell{srcsep2} &\cellbf{378} & \cell{654} & \cell{1922} & \cell{501} & \cell{2122} & \cell{377} & \cellbf{327} & \cell{1596} & \cell{668} & \cell{1312} \\ 
\cell{phantom1} &\cell{293} & \cell{572} & \cell{1687} & \cellbf{249} & \cell{1014} & \cell{292} & \cell{286} & \cell{1401} & \cellbf{166} & \cell{599} \\ 
\cell{blurrycam} &\cellbf{1044} & \cell{6994} & \cell{3420} & \cell{failed} & \cell{6800} & \cellbf{1043} & \cell{3497} & \cell{2850} & \cell{failed} & \cell{3631} \\ 
\cell{blurspike} &\cellbf{1257} & \cell{8384} & \cell{4180} & \cell{failed} & \cell{4127} & \cellbf{1256} & \cell{4192} & \cell{3382} & \cell{failed} & \cell{2161} \\ 
%%%%
\hline
\cell{Problems} & \multicolumn{5}{c}{The objective value $f(\xb^k)$} & \multicolumn{5}{c}{$\Vert \Ab\xb^k - \bb\Vert_2/\Vert \bb \Vert_2\times 10^5$} \\ \hline
%%%
\cell{gcosspike} &\cell{181.484} & \cell{183.050} & \cellbf{181.481} & \cell{181.483} & \cell{181.482} & \cell{0.096} & \cellbf{0.087} & \cell{0.091} & \cell{3.479} & \cell{0.187} \\ 
\cell{p3poly} &\cell{1748.023} & \cell{1838.254} & \cell{1747.982} & \cellbf{1747.954} & \cell{1748.363} & \cell{0.079} & \cell{0.079} & \cell{0.082} & \cell{1.374} & \cellbf{0.001} \\ 
\cell{sgnspike} &\cell{20.620} & \cell{20.620} & \cellbf{20.619} & \cell{20.621} & \cell{20.620} & \cell{0.211} & \cell{0.090} & \cellbf{0.090} & \cell{1.324} & \cell{9.963} \\ 
\cell{zsgnspike} &\cell{28.927} & \cell{28.927} & \cellbf{28.927} & \cell{28.928} & \cell{28.927} & \cell{0.349} & \cell{0.093} & \cellbf{0.092} & \cell{1.598} & \cell{7.169} \\ 
\cell{gausspike} &\cell{24.041} & \cell{24.041} & \cellbf{24.041} & \cell{24.041} & \cell{24.041} & \cell{0.152} & \cell{0.093} & \cellbf{0.092} & \cell{1.628} & \cell{0.112} \\ 
\cell{srcsep1} &\cell{1057.583} & \cell{1059.361} & \cellbf{1057.228} & \cell{1057.974} & \cell{1058.821} & \cell{0.123} & \cell{0.093} & \cellbf{0.091} & \cell{0.954} & \cell{0.647} \\ 
\cell{srcsep2} &\cell{1093.134} & \cell{1094.450} & \cellbf{1092.807} & \cell{1097.060} & \cell{1093.961} & \cell{0.118} & \cell{0.092} & \cellbf{0.094} & \cell{1.050} & \cell{0.446} \\ 
\cell{phantom1} &\cell{202.697} & \cell{202.828} & \cellbf{202.696} & \cell{202.856} & \cell{202.783} & \cell{0.572} & \cell{0.085} & \cellbf{0.085} & \cell{1.148} & \cell{2.412} \\ 
\cell{blurrycam} &\cellbf{10276.681} & \cell{10276.682} & \cell{10276.691} & \cell{failed} & \cell{10276.717} & \cell{0.125} & \cell{0.099} & \cell{0.097} & \cell{failed} & \cellbf{0.075} \\ 
\cell{blurspike} &\cell{576.482} & \cell{576.482} & \cell{576.482} & \cell{failed} & \cellbf{576.474} & \cell{0.125} & \cell{0.100} & \cellbf{0.100} & \cell{failed} & \cell{9.067} \\ 
\hline
\end{tabular}
\end{center}
%\end{scriptsize}
%\end{footnotesize}
\vskip-0.3cm
\end{table*}

We also note that within $s :=1$ to $5$ FISTA iterations, the inexact $\mathrm{(1P1D(s))}$ algorithms still perform well and produce more accurate solutions when the inner iteration number is increasing. 

%%%% 8.5.2. Sparse-group basis pursuit 
\subsubsection{Sparse-group basis pursuit.}
We consider again the sparse-group basis pursuit problem \eqref{eq:group_basis_pursuit}. In this case, we compare our algorithms and the group $\mathrm{YALL1}$ solver \cite{Yang2011}, which we find one of the most efficient algorithm for solving \eqref{eq:group_basis_pursuit}.  A further comparison with SPGL$_1$ can be found in \cite{Yang2011}.

One of the most common ways to compare the performance of different algorithms is using performance profile concept \cite{Dolan2002}.
In this example, we  benchmark seven algorithms with performance profiles.

Recall that a performance profile is built based on a set $\mathcal{S}$ of $n_s$
algorithms (solvers) and a collection $\mathcal{P}$ of $n_p$ problems. Suppose that we build a profile based on computational time (but the same concept can be used for different measurements). 
We denote by
$$T_{p,s} := \textit{computational time required to solve problem $p$ by solver $s$}.$$
We compare the performance of algorithm $s$ on problem $p$ with the best performance of any algorithm on this problem. That is, we compute the performance ratio
$r_{p,s} := \frac{T_{p,s}}{\min\{T_{p,\hat{s}} ~|~ \hat{s}\in \mathcal{S}\}}$.
Now, let 
$$
\tilde{\rho}_s(\tilde{\tau}) := (1/n_p)\mathrm{size}\left\{p\in\mathcal{P} ~|~ r_{p,s}\leq \tilde{\tau}\right\}~\textrm{for}~\tilde{\tau}\in\mathbb{R}_{+}.
$$ 
The function
$\tilde{\rho}_s:\mathbb{R}\to [0,1]$ is the probability for solver $s$ that a performance ratio is within a factor $\tilde{\tau}$ of the best possible ratio. We use the term
``performance profile'' for the distribution function $\tilde{\rho}_s$ of a performance metric.
We plotted the performance profiles in $\log$-scale, i.e. 
$$
\rho_s(\tau) := (1/n_p)\mathrm{size}\left\{p\in\mathcal{P} ~|~ \log_2(r_{p,s})\leq \tau := \log_2\tilde{\tau}\right\}.
$$

The data of this test is generated as follows. The problem size is set to $n := s\times 5120$, $m := \lfloor n/3\rfloor$ and $n_g := \lfloor m/4\rfloor$ for $s := 1,\cdots, 20$. Matrix $\Ab$ is drawn randomly from standard Gaussian distribution with $50\%$ correlated columns. Vector $\bb := \Ab\xb^{\star} + \sigma$, where $\xb^{\star}$ is a given test vector generated also randomly with the standard Gaussian distribution, and $\sigma$ is a Gaussian noise. 

\begin{figure*}[!ht]
%\vskip-0.4cm
\begin{center}
%\centerline{\includegraphics[scale=0.41]{figs/profile_gBP_ns7_03062014}}
\centerline{\includegraphics[width=15cm, height=9.5cm]{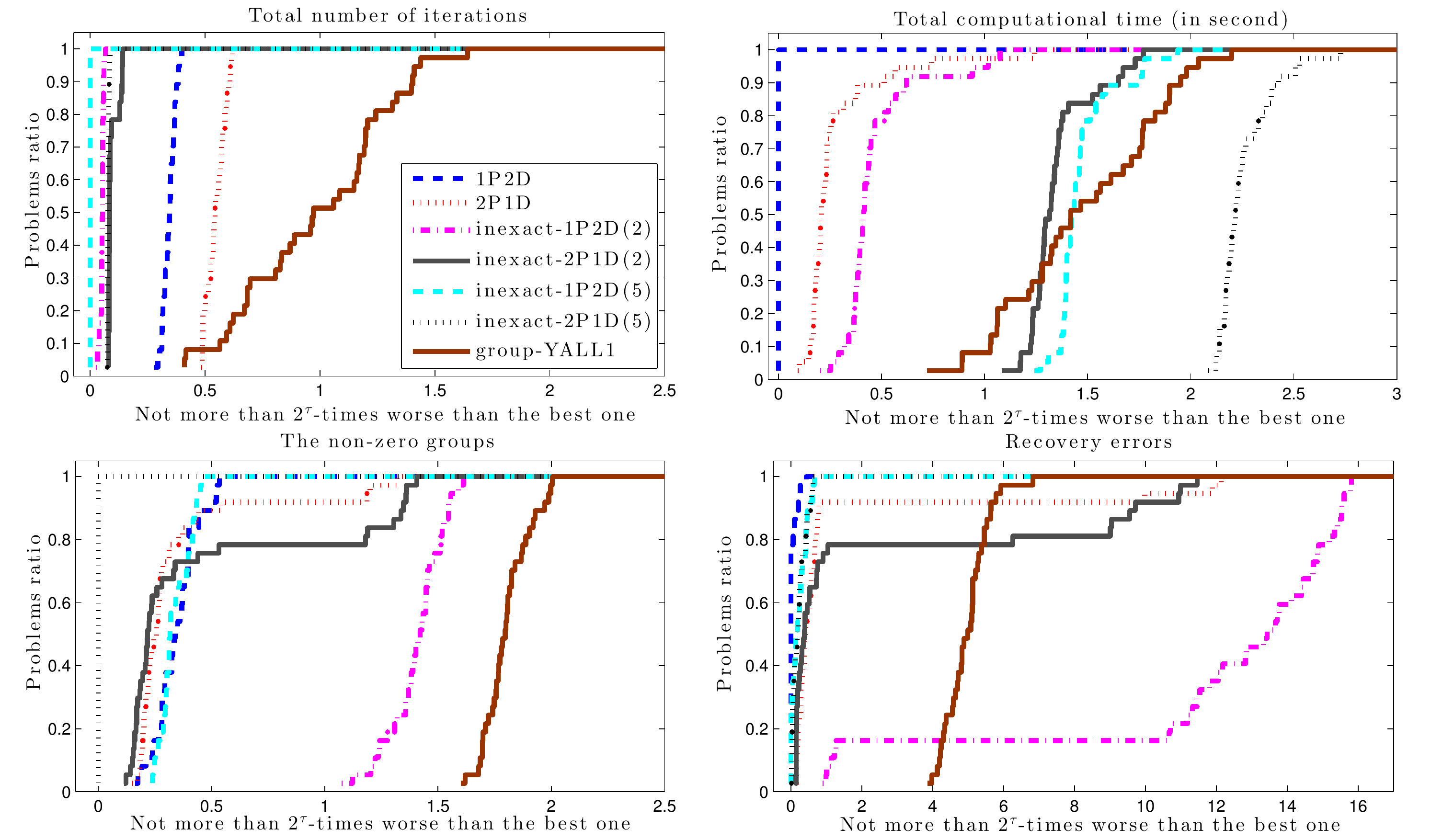}}
\caption{The performance profiles of $7$ algorithms on synthetic data without noise}\label{fig:profile_gBP_05}
\end{center}
\vskip-0.4cm
\end{figure*}

Figure \ref{fig:profile_gBP_05} shows the performance profile of $7$ algorithms: $6$ variants of Algorithm \ref{alg:pd_alg} and $\mathrm{group\_YALL1}$ \cite{Yang2011} in terms of iteration numbers, computational time (in second), the number of nonzero groups and the relative recovery errors $\Vert \xb^k - \xb^{\star}\Vert_2/\Vert\xb^{\star}\Vert_2$.
These performance profiles are built from $37$ problems for size $[m, n, n_g] = [1706, 5120, 427]$ to $[8533, 25600, 2133]$ without additive Gaussian noise.
The $y$-axis of these figures shows the problem ratio $\rho_s(\tau)$. If the problem ratio $\rho_s(\tau)$ is closer to $1$, then the corresponding algorithm has a better performance. 
The $x$-axis shows how many times ($2^{\tau}$) one algorithm is better than the others in $\log_2$-scale.

We can observe from the performance profiles in Figure \ref{fig:profile_gBP_05} for the noiseless case that:
The $\mathrm{(1P2D)}$ variant is the best one in terms of computational time while produces relatively good results (number of nonzero groups, solution recovery errors) compared to the rest. 
The inexact $\mathrm{(2P1D)}$ variant with $5$ FISTA iterations gives the best results  (number of nonzero groups, solution recovery errors)  but is slow due to two primal steps.
While the computational time of our algorithms  slightly increases with respect to the problem size, it increases linearly in $\mathrm{group\_YALL1}$ due to the solution of linear systems.  
The inexact $\mathrm{(2P1D)}$ is more robust to the FISTA iterations than the inexact $\mathrm{(1P2D)}$ one.

\begin{figure*}[!ht]
%\vskip-0.3cm
\begin{center}
%\centerline{\includegraphics[scale=0.41]{figs/profile_gBP_5n_03062014}}
\centerline{\includegraphics[width=15cm, height=9.5cm]{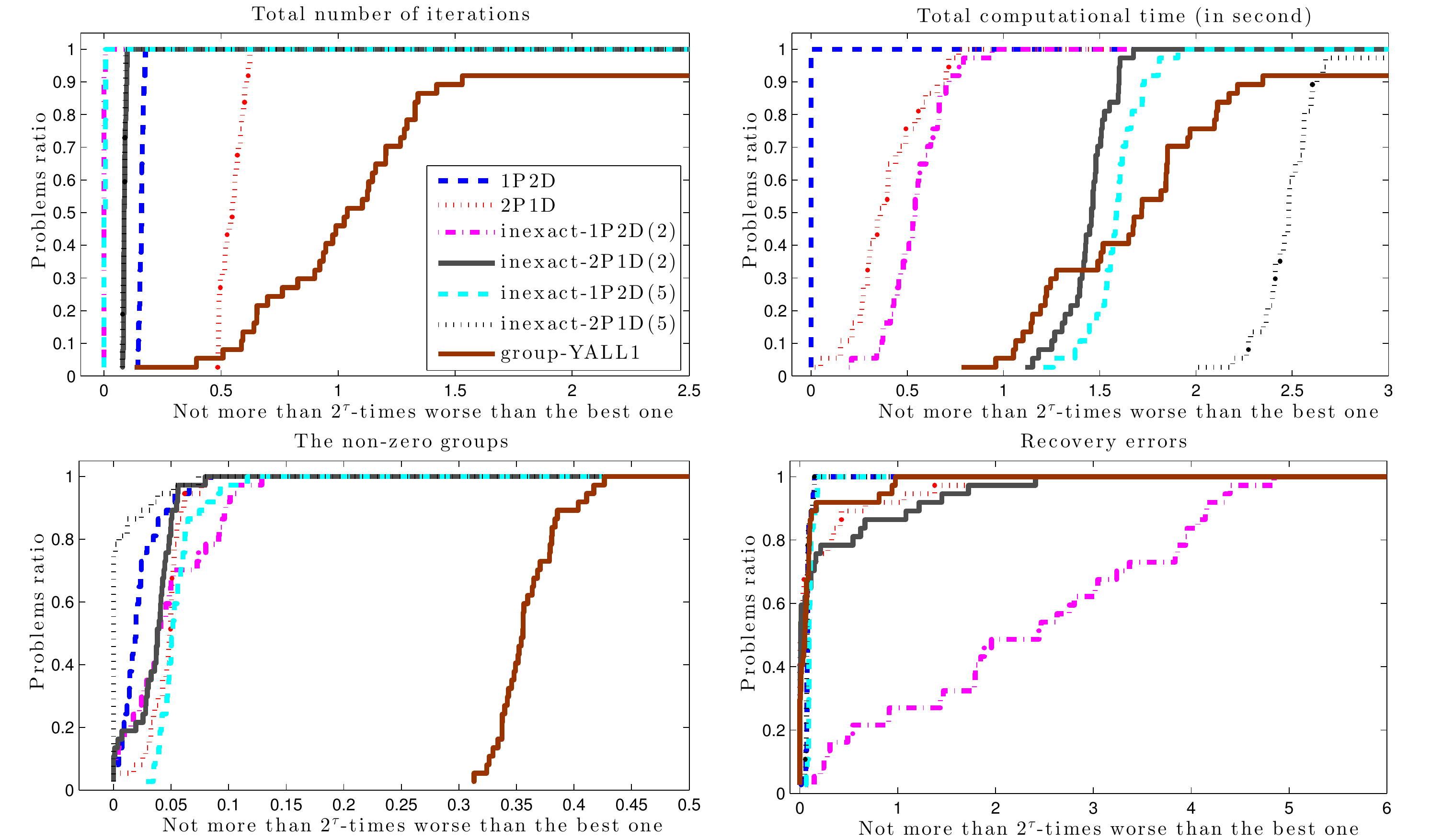}}
\caption{The performance profiles of $7$ algorithms on synthetic data with $5\%$ Gaussian noise}\label{fig:profile_gBP_05_5p}
\end{center}
\vskip -0.4cm
\end{figure*}

Figure \ref{fig:profile_gBP_05_5p} presents the performance profiles when we add $5\%$ Gaussian noise to the model. 
The performance of our algorithms basically remains the same as in the noiseless case, while the number of nonzero groups in $\mathrm{group\_YALL1}$ is increasing significantly compared to ours.
If we increase the noise level up to $10\%$, $\mathrm{group\_YALL1}$ starts oscillating and cannot converge to the solution with the desired accuracy. 
This happens due to the effect of the fixed penalty parameter in $\mathrm{group\_YALL1}$. We note that if we update this parameter, the linear system in $\mathrm{group\_YALL1}$ needs to be resolved, which slows down significantly the performance of the algorithm except some tricks are exploited.

%%%% 8.5.1. Robust PCA 
\subsubsection{Robust principle component analysis.}
We consider the following robust principle component analysis (RPCA) problem:
\begin{equation}\label{eq:robust_PCA}
\min_{\Xb,\Yb \in\mathbb{R}^{m\times n}} \set{ \lambda\norm{\vec{\Xb}}_1 + \norm{ \Yb}_{\star}: \Xb + \Yb = \mathbf{M}},
\end{equation}
where $\mathbf{M}\in\mathbb{R}^{m\times n}$ is a given matrix, $\norm{\cdot}_{\star}$ is the nuclear norm and $\lambda > 0$ is a regularization parameter.
As suggested in \cite{Candes2011}, we can choose $\lambda := \frac{c}{\sqrt{m}}$ to get a perfect recovery (i.e., with high probability), where $c > 0$ is a scaling constant.

In this example, we demonstrate our $(\mathrm{1P2D})$ algorithm on the video clip taken from a surveillance camera in a subway station, which is available at \url{http://perception.i2r.a-star.edu.sg/bk_model/bk_index.html}.
We crop $200$ gray frames from this video clip and preprocess it to obtain a $20'800\times 200$ matrix as an input $\mathbf{M}$. By tuning the regularization parameter $\lambda$, we pick the best possible value $\lambda := 0.01$.
We run our $(\mathrm{1P2D})$ algorithm and compare it with  three other open-source codes: exact ADMM, inexact ADMM \cite{Lin2009} and TFOCS \cite{Becker2011b}. 
All the algorithms are terminated with the same accuracy $10^{-3}$.

The results and performance of these algorithms are reported in Table \ref{tbl:robust_pca}, where $\#\mathrm{svd}$ is the number of SVDs required by the algorithms, $F(\Xb^k,\Yb^k) := \lambda\norm{\vec{\Xb}}_1 + \norm{ \Yb}_{\star}$.
\begin{table*}[!ht]
\newcommand{\cell}[1]{{\!}#1{\!}}
\newcommand{\cellbf}[1]{{\!}\textrm{\color{blue}#1}{\!}}
\begin{center}
%\vskip -0.4cm
\caption{The results and performance of four algorithms on the real-world data}\label{tbl:robust_pca}
\vskip0.25cm
\begin{tabular}{lrrrrr}\hline
Algorithms & $\#\mathrm{iterations}$ & $\#\mathrm{svd}$ & $F(\Xb^k, \Yb^k)$ & {\footnotesize\cell{$\frac{\norm{\Xb^k + \Yb^k - \mathbf{M}}_F}{\norm{\mathbf{M}}_F}$}} & Time[s] \\ \hline
\cell{$(\mathrm{1P2D})$} & \cell{13} & \cellbf{14} & \cellbf{547845.12485} & \cell{0.0004029} & \cell{10.53} \\ 
\cell{exactADMM} &  \cellbf{4} & \cell{662} & \cell{548333.09286} & \cellbf{0.0000676} & \cell{458.75} \\ 
\cell{inexactADMM} & \cell{19} & \cell{19} & \cell{548551.75715} & \cell{0.0004988} & \cellbf{9.33} \\ 
\cell{TFOCS} & \cell{38} & \cell{122} & \cell{566257.63794} & \cell{0.0008508} & \cell{111.89} \\ \hline
\end{tabular}
\end{center}
%\end{footnotesize}
%\vskip-0.4cm
\end{table*}
We can see from Table \ref{tbl:robust_pca}, $(\mathrm{1P2D})$ requires fewest SVD operations and has similar computational time as inexact ADMM, while reaches a better objective value $F(\Xb^k,\Yb^k)$ and the relative feasibility gap. The exact ADMM produces a better solution in terms of quality (lower relative feasibility gap) but requires too many SVDs.

\begin{figure*}[!ht]
\begin{center}
%\vskip-0.4cm
\centerline{\includegraphics[width=15cm, height=14cm]{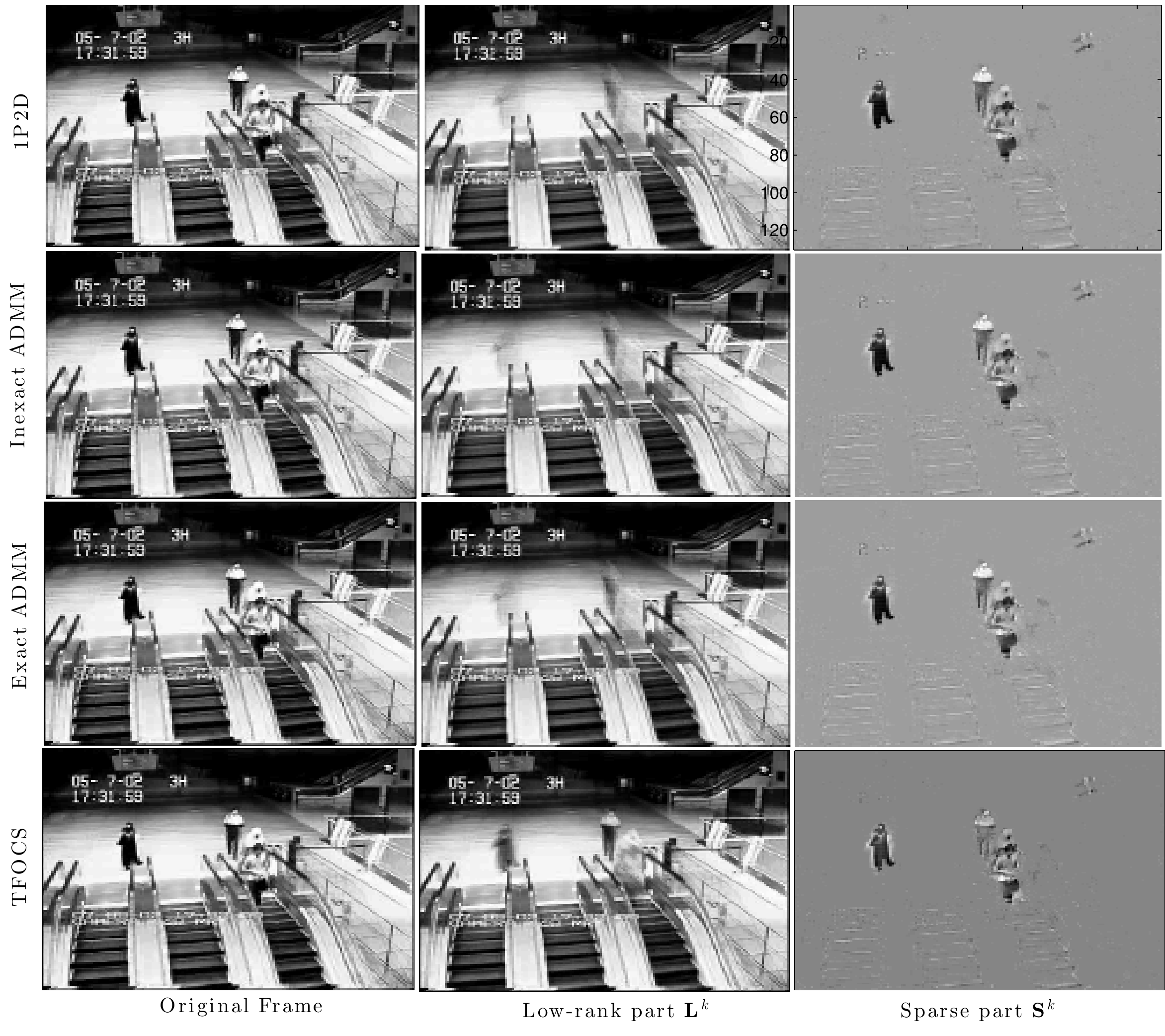}}
%\vskip-0.2cm
\caption{The results of four algorithms on the frame 25 of the video clip}\label{fig:rpca_25}
\end{center}
\vskip -0.4cm
\end{figure*}

The frame 25 of this video is plotted in Figure \ref{fig:rpca_25}, which illustrates how the output of the algorithms can be presented in object separation context.
We can see from this plot that the objects (humans) can be considered as sparse representation and are separated from the background. As can be observed from the second column in Figure \ref{fig:rpca_25}, $(\mathrm{1P2D})$ and ADMMs  give  a better low-rank image estimate as compared to TFOCS.

%% 5.3.1. Square-root LASSO.
\subsubsection{Square-root LASSO.}
Since the  $(\mathrm{1P2D})$ variant of Algorithm \ref{alg:pd_alg} has similar cost-per-iteration as \textrm{ADMM}, we compare this algorithm with the state-of-the-art solvers such as \textrm{TFOCS}, \textrm{ADMM} and \textrm{PADMM}.

For this purpose, we choose the square-root LASSO problem:
\begin{equation}\label{eq:heter_lasso}
\min_{\xb\in\R^n} \norm{\Ab\xb - \bb}_2 \!+\! \lambda\norm{\xb}_1,
\end{equation}
where  $\Ab\in\R^{m\times n}$, $\bb\in\R^m$ are given and $\lambda > 0$ is a regularization term. 
By introducing a new variable $\rb = \Ab\xb - \bb$, \eqref{eq:heter_lasso} can be reformulated in the form of \eqref{eq:separable_convex}:
\begin{equation}\label{eq:heter_lasso2}
\min_{\xb\in\R^n,\rb\in\R^m} \lambda\norm{\xb}_1 + \norm{\rb}_2, ~\mathrm{s.t.} ~\Ab\xb - \bb - \rb = 0.
\end{equation}
As shown in \cite{Belloni2011} that the regularization parameter $\lambda$ can be set at $\lambda = c\Phi^{-1}(1 - 0.5\alpha/n)$ for given $c>1$ and $\alpha \in (0, 1)$.
The suggested values for $c$ and $\alpha$ are $1.1$ and $0.05$, respectively. By choosing this value of $\lambda$, we can probably recover $\xb$ with probability $1-\alpha = 0.95$. 

We mimic the basis pursuit problem before and generate $5$ problems of size $(m, n, s) = i(350, 1000, 100)$, where $i = 1,\dots, 5$ and $s$ is the sparsity. 
We generate the matrix $\Ab$ randomly from Gaussian distribution with $0.5$ correlated columns. Vector $\bb$ is generated as $\bb := \Ab\xb^{\star} + \mathbf{n}$, and $\mathbf{n}$ is Gaussian noise with distribution $\mathcal{N}(0, 0.1)$.

We tune all the augmented Lagrangian algorithms: $(\mathrm{1P2D})$, the preconditioning ADMM (\textrm{PADMM}) and the exact ADMM (\textrm{ADMM}). In these algorithms, we use the same strategy to tune the smoothness parameter $\gamma_k$ and the penalty parameter $\rho_k$, as we observe this works best for three algorithms. The center point $\xb_c^k$ in Algorithm \ref{alg:pd_alg} is chosen as discussed in the enhancement paragraph. In stark contrast to the ADMM and PADMM, our  subproblems with respect to $\xb$ and $\rb$ are solved in \textbf{parallel}. Note that the \textrm{ADMM} requires one matrix inversion $\mathbf{I} + \Ab^T\Ab$. 

A Monte Carlo run of size $10$ shows that our algorithm is not only more accurate but is also faster (cf., Table \ref{tbl:sqrt_lasso}). 
\begin{table*}[!ht]
\newcommand{\cell}[1]{#1}
\newcommand{\cellc}[1]{\textrm{\color{blue}#1}}
\begin{center}
\caption{Performance comparison of Algorithm 1, PADMM, ADMM, and TFOCS.}\label{tbl:sqrt_lasso}
\vskip0.25cm
\begin{tabular}{ccc|rrrr}\hline
\multicolumn{3}{c|}{Size} & \multicolumn{4}{|c}{\# Iterations}  \\ \hline
$m$ & $n$ & $s$ & \cell{$(\mathrm{1P2D})$} & \cell{\textrm{PADMM}} &  \cell{\textrm{ADMM}} & \cell{\textrm{TFOCS}} \\ \hline
\cell{350} & \cell{1000} & \cell{100} & \cellc{ 1331} & \cell{ 1592} & \cell{ 3665} & \cell{ 5000}  \\
\cell{700} & \cell{2000} & \cell{200} & \cellc{ 1311} & \cell{ 1398} & \cell{ 2861} & \cell{ 5000}  \\
\cell{1050} & \cell{3000} & \cell{300} & \cellc{ 1307} & \cell{ 1335} & \cell{ 2797} & \cell{ 5000} \\
\cell{1400} & \cell{4000} & \cell{400} & \cellc{ 1318} & \cell{ 1330} & \cell{ 2631} & \cell{ 5000} \\
\cell{1750} & \cell{5000} & \cell{500} & \cellc{ 1316} & \cell{ 1322} & \cell{ 2594} & \cell{ 5000} \\ \hline
\multicolumn{3}{c|}{Size} & \multicolumn{4}{|c}{$\#\Ab\xb / \#\Ab^T\yb$} \\ \hline
$m$ & $n$ & $s$ & \cell{$(\mathrm{1P2D})$} & \cell{\textrm{PADMM}} &  \cell{\textrm{ADMM}} & \cell{\textrm{TFOCS}} \\ \hline
\cell{350} & \cell{1000} & \cell{100} &\cellc{ 1332/ 2661} & \cell{ 1593/ 3184} & \cell{ 3666/ 7330} & \cell{15996/ 5523} \\ 
\cell{700} & \cell{2000} & \cell{200} & \cellc{ 1312/ 2621} & \cell{ 1399/ 2796} & \cell{ 2862/ 5720} & \cell{16005/ 5548} \\ 
\cell{1050} & \cell{3000} & \cell{300} & \cellc{ 1308/ 2613} & \cell{ 1336/ 2670} & \cell{ 2798/ 5593} & \cell{15989/ 5826} \\ 
\cell{1400} & \cell{4000} & \cell{400} & \cellc{ 1319/ 2635} & \cell{ 1331/ 2659} & \cell{ 2632/ 5260} & \cell{16018/ 5801} \\ 
\cell{1750} & \cell{5000} & \cell{500} & \cellc{ 1317/ 2630} & \cell{ 1323/ 2644} & \cell{ 2595/ 5187} & \cell{16022/ 5790} \\ \hline
\multicolumn{3}{c|}{Size} & \multicolumn{4}{|r}{Objective values $f(\bar{\xb}^k)$~~~~~~~~~~~~~~~~~~}  \\ \hline
$m$ & $n$ & $s$ & \cell{$(\mathrm{1P2D})$} & \cell{\textrm{PADMM}} &  \cell{\textrm{ADMM}} & \cell{\textrm{TFOCS}} \\ \hline
\cell{350} & \cell{1000} & \cell{100} & \cellc{31.424461} & \cell{31.424537} & \cell{31.424762} & \cell{32.652869} \\
\cell{700} & \cell{2000} & \cell{200} & \cellc{74.917422} & \cell{74.917552} & \cell{74.919787} & \cell{77.039976} \\ 
\cell{1050} & \cell{3000} & \cell{300} & \cellc{120.904351} & \cell{120.904523} & \cell{120.909089} & \cell{123.684820} \\
\cell{1400} & \cell{4000} & \cell{400} & \cellc{150.458042} & \cell{150.458275} & \cell{150.465146} & \cell{156.510366} \\ 
\cell{1750} & \cell{5000} & \cell{500} & \cellc{192.030170} & \cell{192.030441} & \cell{192.040217} & \cell{201.906842} \\ \hline
\multicolumn{3}{c|}{Size} & \multicolumn{4}{|r}{Recovery errors $\norm{\bar{\xb}^k - \xb^{\star}}/\norm{\xb^{\star}}$~~~~~} \\ \hline
$m$ & $n$ & $s$ &  \cell{$(\mathrm{1P2D})$} &  \cell{\textrm{PDMM}} & \cell{\textrm{ADMM}} & \cell{\textrm{TFOCS}} \\ \hline
\cell{350} & \cell{1000} & \cell{100} & \cell{0.15120} & \cell{0.15122} & \cell{0.15180} & \cellc{0.14713} \\ 
\cell{700} & \cell{2000} & \cell{200} & \cell{0.04689} & \cell{0.04689} & \cell{0.04707} & \cellc{0.04447} \\ 
\cell{1050} & \cell{3000} & \cell{300} &\cell{0.03165} & \cell{0.03166} & \cell{0.03181} & \cellc{0.02947} \\ 
\cell{1400} & \cell{4000} & \cell{400} &\cellc{0.03013} & \cell{0.03014} & \cell{0.03025} & \cell{0.04040} \\ 
\cell{1750} & \cell{5000} & \cell{500} & \cellc{0.03802} & \cell{0.03803} & \cell{0.03824} & \cell{0.04973} \\ \hline
\end{tabular}
\end{center}
\vskip-0.4cm
\end{table*}
We count the number of matrix-vector multiplications both in $\Ab\xb$ and $\Ab^T\xb$ since these are more expensive than the prox operators. Since the iterative vector $\xb$ is sparse, the multiplication $\Ab^T\yb$ is  more expensive. As we can see through this example that \textrm{ADMM} requires more iterations than Algorithm \ref{alg:pd_alg} and \textrm{PADMM} while produces lower accurate solutions. At the same time, TFOCS is slowest and least accurate while sometimes obtaining better estimation error. 

\subsubsection{Binary linear support vector machine.}
This example is concerned the following binary linear support vector machine problem of the Hinge loss function:
\begin{equation}\label{eq:binary_svm}
\min_{\xb\in\R^n}\Big\{ F(\xb) := \sum_{j=1}^m\ell_j(y_j, \wb_j^T\xb - \bb_j) + g(\xb) \Big\},
\end{equation}
where $\ell_j(s, \tau)$ is the Hinge loss function given by $\ell_j(s,\tau) := \max\set{0, 1 - s\tau} = [1-s\tau]_{+}$, $\wb_j$ is the column of a given matrix $\mathbf{W}\in\R^{m\times n}$, $\bb \in\R^n$ is the bias vector, $\yb\in\set{-1,+1}^m$ is a classifier vector $g$ is a given regularization function, e.g., $g(\xb) := \frac{\lambda}{2}\norm{\xb}^2_2$ for $\ell_2$-regularizer or $g(\xb) := \lambda\norm{\xb}_1$ for $\ell_1$-regularizer ($\lambda > 0$ is a regularization parameter).

By introducing a slack variable $\rb = \mathbf{W}\xb - \bb$, we can write \eqref{eq:binary_svm} in terms of \eqref{eq:separable_convex} as
\begin{equation}\label{eq:binary_svm2}
\begin{array}{ll}
\displaystyle\min_{\xb\in\R^n, \rb\in\R^m} & \set{ \sum_{j=1}^m\ell_j(y_j, \rb_j) + g(\xb)} \\
\mathrm{s.t.} & \mathbf{W}\xb - \rb = \bb,
\end{array}
\end{equation}
Now, we can apply the $(\mathrm{1P2D})$ variant to solve this resulting problem. 
We test  this algorithm on \eqref{eq:binary_svm2} and compare it with \textrm{LibSVM} \cite{Chang2011}. 
We select only two problems from the \textrm{LibSVM} data set available at \url{http://www.csie.ntu.edu.tw/~cjlin/libsvmtools/datasets/}  for our test. The first problem is \texttt{a1a}, which has size $p = 119$ features and $N = 1605$ data points, while the second problem is \texttt{news20}, which has size $p = 1355191$ features and $N = 19996$ data points.

We compare two algorithms in terms of the final objective values $F(\xb^k)$, the classification accuracy $\mathrm{ca}_{\lambda} :=  1 - N^{-1}\sum_{j=1}^N\left[\mathrm{sign}(\mathbf{W}\xb^k - \rb) \neq \yb)\right]$ and the computational time.
The results of our test are reported in Table \ref{tbl:svm_results}.
\begin{table*}[!ht]
\newcommand{\cell}[1]{{\!\!\!}#1{\!\!}}
\newcommand{\cellbf}[1]{{\!\!\!}\textrm{#1}{\!\!}}
\begin{center}
\caption{The results of two algorithms on two real-world data problems}\label{tbl:svm_results}
\vskip0.25cm
%\begin{footnotesize}
\begin{tabular}{l|rrrrrrrrrr}\hline
\cell{Problem} & \multicolumn{10}{c}{\cell{The parameter values}} \\ \hline
$\lambda^{-1}$  & \cell{$10^{-3}$} & \cell{$111.1$} & \cell{$222.2$} & \cell{$333.3$} & \cell{$444.4$} & \cell{$555.6$} & \cell{$666.7$} & \cell{$777.8$} & \cell{$888.9$} & \cell{$10^3$} \\ \hline
\multicolumn{11}{c}{\cell{The accuracy of problem \texttt{a1a}}} \\ \hline
\cell{$(\mathrm{1P2D})$}      & \cellbf{0.7539}  &  \cellbf{0.8717} & \cellbf{0.8717} & \cellbf{0.8710} & \cellbf{0.8710} & \cellbf{0.8710} &  \cellbf{0.8710}  & \cellbf{0.8710} & \cellbf{0.8710}  & \cellbf{0.8710} \\
\cell{\texttt{LibSVM}}  & \cellbf{0.7539}  & \cell{0.8692}  & \cell{0.8698} & \cell{0.8698} & \cell{0.8698} & \cell{0.8698} & \cell{0.8698}   & \cell{0.8698} &  \cell{0.8679} & \cell{0.8698} \\ 
\hline
\multicolumn{11}{c}{\cell{The CPU time [in second] of problem \texttt{a1a}}} \\ \hline
\cell{$(\mathrm{1P2D})$} & \cell{4.4045} & \cell{4.3769} & \cell{4.4246} & \cell{4.4941} & \cell{4.6238} & \cell{4.5175} & \cell{4.4836} & \cell{4.4719} & \cell{4.7179} & \cell{4.8097} \\ 
\cell{\texttt{LibSVM}} & \cell{0.2549} & \cell{2.1909} & \cell{4.3884} & \cell{5.8583} & \cell{8.3662} & \cell{11.2350} & \cell{11.7036} & \cell{12.9832} & \cell{17.1424} & \cell{17.4362} \\ \hline
\multicolumn{11}{c}{\cell{The accuracy of problem \texttt{news20}}} \\ \hline
\cell{$(\mathrm{1P2D})$} & \cell{0.5001} & \cell{0.9987} & \cell{0.9987} & \cell{0.9987} & \cell{0.9987} & \cell{0.9987} & \cell{0.9987} & \cell{0.9987} & \cell{0.9987} & \cell{0.9987}  \\
\cell{\texttt{LibSVM}} & \cell{0.5001} & \cell{0.9987} & \cell{0.9987} & \cell{0.9987} & \cell{0.9987} & \cell{0.9988} & \cell{0.9988} & \cell{0.9988} & \cell{0.9988} & \cell{0.9988} \\ \hline
\multicolumn{11}{c}{\cell{The CPU time [in second] of problem \texttt{news20}}} \\ \hline
\cell{$(\mathrm{1P2D})$} & \cell{762.31} & \cell{1023.22} & \cell{994.64} & \cell{1043.06} & \cell{984.24} & \cell{989.70} & \cell{1064.33} & \cell{1073.94} & \cell{984.47} & \cell{1018.35} \\
\cell{\texttt{LibSVM}} & \cell{890.26} & \cell{1440.28} & \cell{1449.23} & \cell{1439.77} & \cell{1434.27} & \cell{1518.56} & \cell{1560.38} & \cell{1557.48} & \cell{1535.19} & \cell{1530.71} \\ \hline
\end{tabular}
%\end{footnotesize}
\end{center}
%\vskip-0.4cm
\end{table*}

As can be seen from these results that both solvers give relatively the same objective values, the accuracy for these two problems, while the computational of $(\mathrm{1P2D})$ is much lower than \texttt{LibSVM}. We note that \texttt{LibSVM} was implemented in C++ while $(\mathrm{1P2D})$ is simply a Matlab code. \texttt{LibSVM} becomes slower when the parameter $\lambda$ getting smaller due to the active-set strategy.
The  $(\mathrm{1P2D})$ algorithm is almost independent of the regularization parameter $\lambda$, which is different from active-set methods. In addition, the performance of $(\mathrm{1P2D})$ can be improved by taking account its parallelization ability, which has not been exploited yet in our Matlab implementation.

\begin{figure*}[ht!]
%\vskip-0.3cm
\begin{center}
%\centerline{\includegraphics[scale=0.42]{figs/a1a_plot}}
\centerline{\includegraphics[width=15.5cm, height=4.5cm]{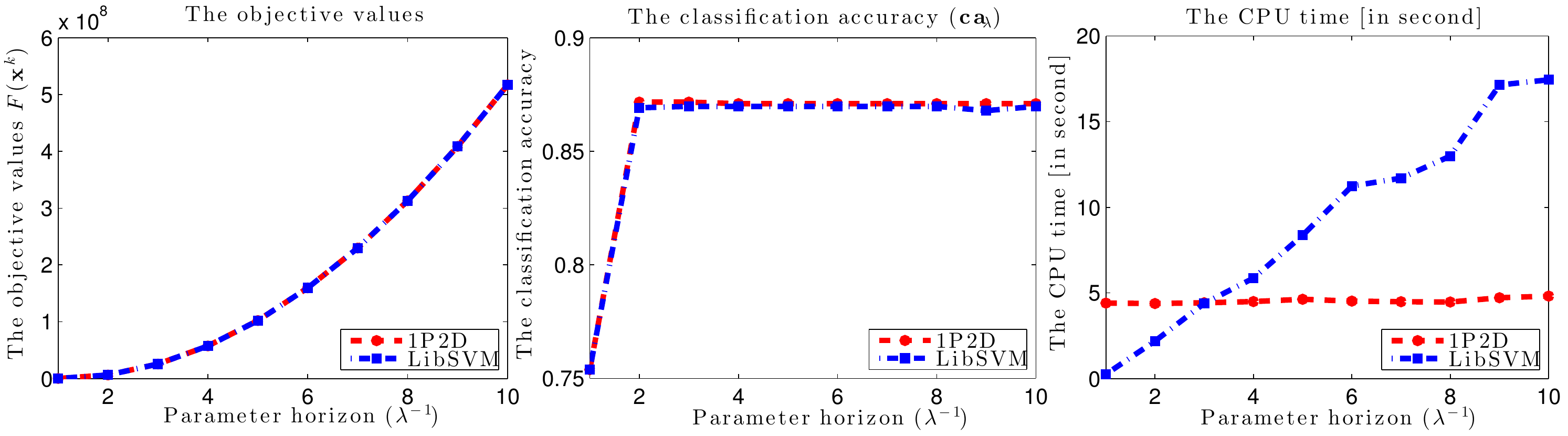}}
\caption{The results of two algorithms on the small $\texttt{a1a}$ problem}\label{fig:a1a_plot}
\end{center}
%\vskip-0.4cm
\end{figure*}

\begin{figure*}[ht!]
%\vskip-1cm
\begin{center}
%\centerline{\includegraphics[scale=0.42]{figs/news20_plot}}
\centerline{\includegraphics[width=15.5cm, height=4.5cm]{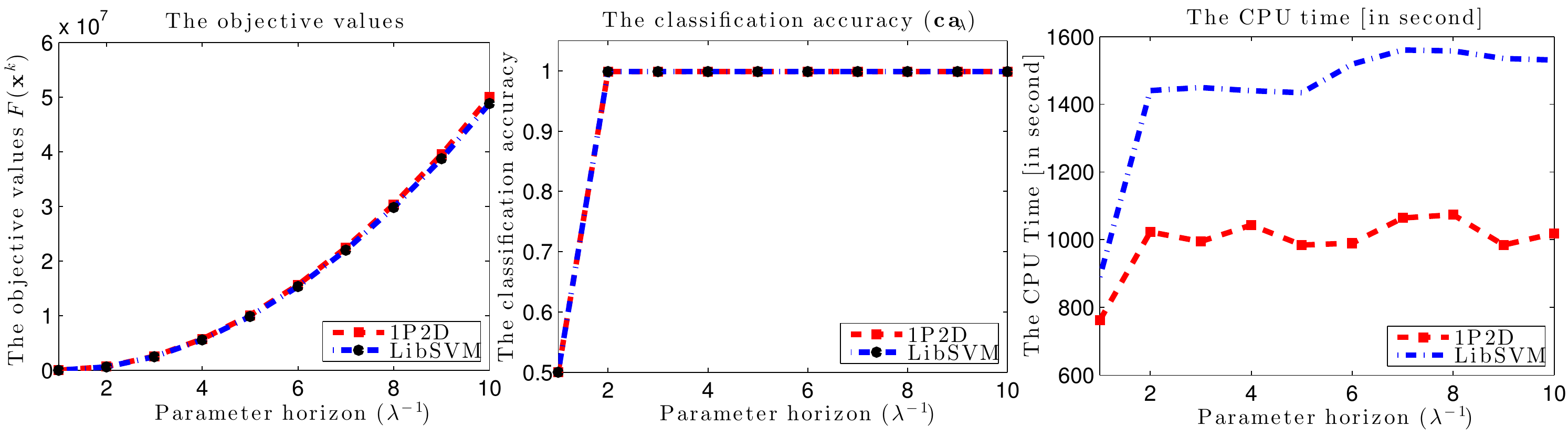}}
\caption{The results of two algorithms on the large-scale $\texttt{news20}$ problem}\label{fig:news20_plot}
\end{center}
%\vskip-1cm
\end{figure*}

To immediately see the performance without looking at the numbers in Table \ref{tbl:svm_results}, we plot the results in Figures \ref{fig:a1a_plot} and \ref{fig:news20_plot} for two separate problems, respectively.

\section{Conclusions}\label{sec: conclusions}
%We have presented a unified primal-dual algorithmic framework for solving a broad class of constrained convex optimization. 
We introduce a model-based excessive gap (MEG) technique for constructing and analyzing first-order methods that numerically approximate an optimal solution of \eqref{eq:separable_convex}. Thanks to a combination of smoothing strategies and MEG, we introduce, to the best of our knowledge, the first algorithmic schemes for \eqref{eq:separable_convex} that theoretically obtain optimal convergence rates directly without averaging the iterates and that seamlessly handle the $p$-decomposability structure. Surprisingly, our analysis techniques enable inexact characterizations, which is important for the augmented Lagrangian versions with lower-iteration counts. We expect a deeper understanding of MEG and different smoothing strategies to help us in tailoring adaptive update strategies for our schemes (as well as several other connected and well-known schemes) 
in order to further improve the empirical performance.  
%
%
%Scalability, rigorous convergence rates, minimal assumptions, and good practical performance are the hallmarks of the primal-dual decomposition algorithms we describe in this paper. 
%The basic algorithm only applies the proximity operators of the functions and the linear operators or their adjoints with care---only once per iteration provided that these functions possess a tractable proximity operator. Moreover, by a suitable choice of the prox-functions, the algorithms can be implemented in a parallel and distributed manner. We also demonstrate that our framework covers important existing algorithms as its variants, and provides new convergence guarantees. % inspirited from the main framework.
%
%
%Interestingly, the additional operator computations for the adaptive kicking can significantly enhance convergence as it can decrease the numerator of the convergence guarantee drastically. Finally, the class of optimization problems \eqref{eq:separable_convex} captures with two terms is already large. Notably, our algorithmic framework retains its rigorous convergence rates when $f$ is separable to any number of terms, as long as the individual proximal maps are tractable. 

\section*{Acknowledgments}
%\vskip-0.2cm
This work is supported in part by the European Commission under the grants MIRG-268398 and ERC Future Proof, and by the Swiss Science Foundation under the grants SNF 200021-132548, SNF 200021-146750 and SNF CRSII2-147633.
\appendix
\section{The proofs of technical statements}\label{sec:appendix}
This appendix provides the technical proofs of Lemmas and Theorems introduced in the main text.

% A.1. The proof of Lemma 3.5.
\subsection{The proof of Lemma \ref{le:excessive_gap_aug_Lag_func}: Bounds on the objective residual and feasibility gap.}
By induction, it follows from Definition \ref{de:decreasing_sequence} that $G_k(\bar{\wb}^k) \leq \omega_kG_0(\bar{\wb}^0) - \Psi_k$, where $\omega_k := \prod_{j=0}^{k-1}(1-\tau_j)$ and $\Psi_k :=  \Psi_0 + \sum_{j=1}^{k-1}\prod_{l=0}^{j-1}(1-\tau_l)\Psi_j$. Using the definition \eqref{eq:smoothed_gap_func} of $G_k$ and the definition \eqref{eq:g_gamma} of $g_{\gamma}$, we can reexpress $G_k$ as $G_k(\bar{\wb}^k) = f(\bar{\xb}^k) - g_{\gamma_k}(\bar{\yb}^k) +  (1/(2\beta_k))\Vert\Ab\bar{\xb}^k - \bb\Vert_2^2$.
This expression leads to 
\begin{equation}\label{eq:main_bound_gap}
f(\bar{\xb}^k) - g_{\gamma_k}(\bar{\yb}^k) \leq \omega_kG_{0}(\bar{\wb}^0) - \Psi_k -  (1/(2\beta_k))\norm{\Ab\bar{\xb}^k - \bb}_2^2,
\end{equation}
which is indeed \eqref{eq:main_bound_gap}.

Now, we notice that under Assumption A.\ref{as:A1}, the solution set $\Yc^{\star}$ of the dual problem \eqref{eq:dual_prob} is also nonempty and bounded. Moreover, the strong duality holds, i.e., $f^{\star} = g^{\star}$. Any point $(\xb^{\star}, \yb^{\star})\in\Xc^{\star}\times\Yc^{\star}$ is a primal-dual solution to \eqref{eq:separable_convex}-\eqref{eq:dual_prob}, and is also a saddle point of $\Lc$, i.e., $\Lc(\xb^{\star},\yb) \leq \Lc(\xb^{\star}, \yb^{\star}) \leq \Lc(\xb, \yb^{\star})$ for all $\xb\in \Xc$ and $\yb\in\R^m$. 
These inequalities lead to the following estimate
\begin{equation}\label{eq:pd_lower_bound}
f(\xb) - g(\yb) \geq f(\xb) -  f^{\star}  \geq -\norm{\yb^{\star}}_2\norm{\Ab\xb - \bb}_2, ~~\forall \xb\in\Xc, ~\yb\in\R^m,
\end{equation}
which is exactly \eqref{eq:pd_lower_bound}.
Now, we combine \eqref{eq:smoothed_gap_of_g}, \eqref{eq:main_bound_gap} and \eqref{eq:pd_lower_bound} to get the following:
\begin{align}
- \norm{\yb^{\star}}_2\Vert \Ab\bar{\xb}^k - \bb\Vert_2 \leq f(\bar{\xb}^k) - f^{\star} \leq f(\bar{\xb}^k) - g(\bar{\yb}^k)  \leq S_k -  (1/(2\beta_k))\norm{\Ab\bar{\xb}^k - \bb}_2^2 \leq S_k,
\label{eq:main_bound_gap_2a1}
\end{align}
where $S_k := \omega_kG_0(\bar{\wb}^0) + \gamma_kD_{\Xc}^{\Sb} - \Psi_k$. This bound is exactly \eqref{eq:main_bound_gap_2a}.

Finally, we prove \eqref{eq:main_bound_gap_2b}.
Let $t := \Vert\Ab\bar{\xb}^k - \bb\Vert_2$. It follows from \eqref{eq:main_bound_gap_2a1} that $-\Vert\yb^{\star}\Vert_2t \leq S_k -  (1/(2\beta_k))t^2$. This inequation of $t$ leads to $t \leq \beta_k\Big[\norm{\yb^{\star}}_2 + \sqrt{\norm{\yb^{\star}}_2^2 + 2\beta_k^{-1}S_k}\Big]$ provided $\beta_k\Vert\yb^{\star}\Vert_2^2 + 2S_k \geq 0$. This estimate is indeed \eqref{eq:main_bound_gap_2b}.
\Eproof
% End of the proof.

% A.2. The proof of Lemma 4.1.
\subsection{Convergence analysis: The proof of  Theorem \ref{th:convergence_A1}.}
Our proof of  Theorem \ref{th:convergence_A1} takes the following outline:
\begin{enumerate}
\item We prove two key lemmas: Lemma \ref{le:maintain_excessive_gap1} and Lemma \ref{le:maintain_excessive_gap2}. These lemmas provide conditions to update the step-size $\tau_k$.
\item We show how to find starting points for Algorithm \ref{alg:pd_alg} using Lemma \ref{le:starting_point}.
\item We provide an update rule for the step-size parameter $\tau_k$ in Lemma \ref{le:update} based on the conditions of Lemmas \ref{le:maintain_excessive_gap1} and  \ref{le:maintain_excessive_gap2}.
\item We combine the above results to finalize the proof of Theorem \ref{th:convergence_A1}.
\end{enumerate}

\subsubsection{The proof of Lemma \ref{le:maintain_excessive_gap1}: The condition for selecting step-size $\tau_k$ in \eqref{eq:pd_scheme_2p}.} 
Let us denote by $d_k(\wb) := \gamma_kd_b(\Sb\xb,\Sb\xb_c) + (\beta_k/2)\Vert\yb\Vert_2^2$ and $\bar{\xb}^{\star}_k := \xb^{\star}_{\gamma_k}(\bar{\yb}^k)$.
If we define 
\begin{equation}\label{eq:H_k0}
H_k(\wb) := f(\bar{\xb}^k) - f(\xb) + F(\bar{\wb}^k)^T(\bar{\wb}^k - \wb) - d_{k}(\wb),
\end{equation}
the objective function in \eqref{eq:smoothed_gap_func}, then by the definition of $G_{k+1}$ and $\Wc := \Xc\times \R^m$, we have
\begin{align}\label{eq:lm41_proof1}
G_{k+1}(\bar{\wb}^{k+1}) := \max_{\wb\in\Wc} H_{k+1}(\wb).
\end{align}
The proof is divided in the following steps:

\noindent\textit{\underline{Step 1:} Splitting $H_k$ and $H_{k+1}$.}
Using the definition of $F$ in \eqref{eq:VIP}, we can write $F(\bar{\wb}^{k})^T(\bar{\wb}^{k} - \wb) =  (\Ab\bar{\xb}^{k} - \bb)^T\yb - (\Ab\xb - \bb)^T\bar{\yb}^{k}$.
Plugging this expression into \eqref{eq:H_k0} we obtain
\begin{equation}\label{eq:H_k}
H_k(\wb) := f(\bar{\xb}^k) - f(\xb) + (\Ab\bar{\xb}^{k} - \bb)^T\yb - (\Ab\xb - \bb)^T\bar{\yb}^{k} -  d_{k}(\wb).
\end{equation}
Similarly to \eqref{eq:H_k}, we also have $H_{k+1}(\wb) = f(\bar{\xb}^{k+1}) - f(\xb) + (\Ab\bar{\xb}^{k+1} - \bb)^T\yb - (\Ab\xb - \bb)^T\bar{\yb}^{k+1}$. Using this expression and $\bar{\yb}^{k+1} = (1-\tau_k)\bar{\yb}^k + \tau_k\hat{\yb}^k$ in \eqref{eq:pd_scheme_2p} we get
\begin{align*}
H_{k+1}(\wb) = f(\bar{\xb}^{k+1}) - f(\xb) + (\Ab\bar{\xb}^{k+1} - \bb)^T\yb -  (1 \!-\! \tau_k)(\Ab\xb \!-\! \bb)^T\bar{\yb}^{k} \!-\! \tau_k(\Ab\xb \!-\! \bb)^T\hat{\yb}^k - d_{k+1}(\wb).
\end{align*}
By adding and then subtracting $(1-\tau_k)[f(\bar{\xb}^k) - f(\xb)]$ into this inequality, we obtain
\begin{align}\label{eq:lm41_proof2}
H_{k+1}(\wb) &= (1-\tau_k)\big[ f(\bar{\xb}^k - f(\xb) - (\Ab\xb - \bb)^T\bar{\yb}^k \big] + (\Ab\bar{\xb}^{k+1} - \bb)^T\yb  - \tau_k(\Ab\xb  - \bb)^T\hat{\yb}^k \nonumber\\
&+ f(\bar{\xb}^{k+1}) - (1-\tau_k)f(\bar{\xb}^k) - \tau_kf(\xb)  - d_{k+1}(\wb).
\end{align}

\noindent\textit{\underline{Step 2:} Estimating a lower bound for $G_{k}$.}
By using the definition \eqref{eq:g_gamma} of $g_{\gamma}$, we have
\begin{align*}
f(\xb) + (\Ab\xb - \bb)^T\bar{\yb}^k + \gamma_kd_b(\Sb\xb,\Sb\xb_c) \geq g_{\gamma_k}(\bar{\yb}^k) + \gamma_kd_b(\Sb\xb,\Sb\xb^{\star}_{\gamma_k}(\bar{\yb}^k)).
\end{align*}
Using this inequality and  $\displaystyle\max_{\yb\in\R^m}\Big\{(\Ab\bar{\xb}^k - \bb)^T\yb - (\beta/2)\norm{\yb}_2^2\Big\} = (1/(2\beta))\norm{\Ab\bar{\xb}^k - \bb}_2^2$, we can show that
\begin{align}\label{eq:lm41_proof5}
G_k(\bar{\wb}^k) &:={\!\!} \max_{\wb\in\Wc}\!\!\Big\{\! f(\bar{\xb}^{k}) \!-\! f(\xb) \!-\! (\Ab\xb \!\!-\! \bb)^T\!\bar{\yb}^k \!-\! \gamma_kd_b(\Sb\xb,\Sb\xb_c) \!+\! (\Ab\bar{\xb}^k \!\!-\! \bb)^T\!\yb \!-\! (\beta_k/2)\!\norm{\yb}_2^2 \!\Big\} \nonumber\\
& \geq f(\bar{\xb}^k) - f(\xb) - (\Ab\xb - \bb)^T\bar{\yb}^k - \gamma_kd_b(\Sb\xb,\Sb\xb_c) + (1/(2\beta_k))\norm{\Ab\bar{\xb}^k - \bb}_2^2 + \gamma_kd_b(\Sb\xb,\Sb\bar{\xb}^{\star}_k). 
\end{align}
From the second line $\hat{\yb}^k := \beta_{k+1}^{-1}(\Ab\hat{\xb}^k - \bb)$ of \eqref{eq:pd_scheme_2p}, we also have the following equality
\begin{align}\label{eq:lm41_proof6}
\norm{\Ab\bar{\xb}^k - \bb}_2^2 &= \norm{\Ab\hat{\xb}^k-\bb}_2^2 + 2(\Ab\hat{\xb}^k - \bb)^T\Ab(\bar{\xb}^k - \hat{\xb}^k) + \norm{\Ab(\hat{\xb}^k - \bar{\xb}^k)}^2_2\nonumber\\
&=\norm{\Ab\hat{\xb}^k-\bb}_2^2 + 2\beta_{k+1}(\hat{\yb}^k)^T\Ab(\bar{\xb}^k - \hat{\xb}^k) + \norm{\Ab(\hat{\xb}^k - \bar{\xb}^k)}^2_2.
\end{align}
Since $\beta_{k+1} = (1-\tau_k)\beta_k$ due to \eqref{eq:update_param1}, substituting \eqref{eq:lm41_proof6} into \eqref{eq:lm41_proof5} we obtain
\begin{align}\label{eq:lm41_proof7}
f(\bar{\xb}^k) \!-\! f(\xb)  \!-\!  (\Ab\xb \!-\! \bb)^T\!\bar{\yb}^k \!-\! \gamma_kd_b(\Sb\xb,\Sb\xb_c)  &\!\leq\!  G_k(\bar{\wb}^k) \!-\! (1/(2\beta_{k\!+\!1}))\norm{\Ab\hat{\xb}^k \!-\! \bb}_2^2 \!-\! (\hat{\yb}^k)^T\!\Ab(\bar{\xb}^k \!-\! \hat{\xb}^k) \\
& \!-\! \gamma_kd_b(\Sb\xb,\Sb\bar{\xb}^{\star}_k)  \!-\! (1/(2\beta_{k\!+\!1}))\big[\norm{\Ab(\bar{\xb}^k \!-\! \hat{\xb}^k)}_2^2 \!-\! \tau_k\norm{\Ab\bar{\xb}^k \!-\! \bb}^2_2\big]. \nonumber
\end{align}

\noindent\textit{\underline{Step 3:} Estimating an upper bound for $H_{k+1}$.}
First, from the update rule \eqref{eq:update_param1}, we have $\gamma_{k+1} = (1-c_k\tau_k)\gamma_k \geq (1-\tau_k)\gamma_k$ for any $c_k \leq 1$ and $\beta_{k+1} = (1-\tau_k)\beta_k$. Hence, we can show that
\begin{align}\label{eq:lm41_proof3}
d_{k\!+\!1}(\wb) = \gamma_{k\!+\!1}d_b(\Sb\xb,\Sb\xb_c) \!+\! (\beta_{k\!+\!1}/2)\norm{\yb}_2^2 \geq (1 \!-\! \tau_k)\gamma_kd_b(\Sb\xb,\Sb\xb_c)  +  (\beta_{k+1}/2)\norm{\yb}_2^2.
\end{align}
Second, by using $\hat{\yb}^k := \beta_{k+1}^{-1}(\Ab\hat{\xb}^k - \bb)$ in \eqref{eq:pd_scheme_2p}, we also have the following equality
\begin{align}\label{eq:lm41_proof8}
(\Ab\xb - \bb)^T\hat{\yb}^k &= (\hat{\yb}^k)^T\Ab(\hat{\xb}^k - \bb) + (\hat{\yb}^k)^T\Ab(\xb - \hat{\xb}^k) = (1/\beta_{k+1})\norm{\Ab\hat{\xb}^k - \bb}_2^2 + (\hat{\yb}^k)^T\Ab(\xb - \hat{\xb}^k).
\end{align}
Third, substituting \eqref{eq:lm41_proof3}, \eqref{eq:lm41_proof7} and \eqref{eq:lm41_proof8} into \eqref{eq:lm41_proof2}, we can upperbound the estimate $H_{k+1}$ as
\begin{align}\label{eq:lm41_proof9}
H_{k+1}(\wb) &\overset{\tiny\eqref{eq:lm41_proof3}}{\leq} (1-\tau_k)\big[f(\bar{\xb}^k) - f(\xb) - (\Ab\xb \!-\! \bb)^T\bar{\yb}^{k} - \gamma_kd_b(\Sb\xb,\Sb\xb_c)\big] + (\Ab\bar{\xb}^{k+1} - \bb)^T\yb - (\beta_{k+1}/2)\Vert\yb\Vert_2^2 \nonumber\\
& - \tau_k(\Ab\xb  - \bb)^T\hat{\yb}^k  + f(\bar{\xb}^{k+1}) - (1-\tau_k)f(\bar{\xb}^k) - \tau_kf(\xb)\nonumber\\
&\overset{\tiny\eqref{eq:lm41_proof7}+\eqref{eq:lm41_proof8}}{\leq}  (1-\tau_k)G_k(\bar{\wb}^k) + (\Ab\bar{\xb}^{k + 1}  - \bb)^T\yb - (\beta_{k+1}/2)\norm{\yb}_2^2 + \big[f(\bar{\xb}^{k+1}) - (1-\tau_k)f(\bar{\xb}^k) - \tau_kf(\xb) \big] \nonumber\\
&- \frac{(1-\tau_k)}{2\beta_{k+1}}\big[\norm{\Ab\hat{\xb}^k-\bb}_2^2 + \norm{\Ab(\bar{\xb}^k - \hat{\xb}^k)}_2^2 - \tau_k\norm{\Ab\bar{\xb}^k - \bb}^2_2\big]  - (1-\tau_k) (\hat{\yb}^k)^T\Ab(\bar{\xb}^k-\hat{\xb}^k) \nonumber\\
& - \tau_k\big[(1/\beta_{k+1})\norm{\Ab\hat{\xb}^k - \bb}_2^2 - (\hat{\yb}^k)^T\Ab(\xb -\hat{\xb}^k)\big] - (1-\tau_k)\gamma_kd_b(\Sb\xb,\Sb\bar{\xb}^{\star}_k)  \nonumber\\
&= (1-\tau_k)G_k(\bar{\wb}^k)  + \big[f(\bar{\xb}^{k+1}) - (1-\tau_k)f(\bar{\xb}^k) - \tau_kf(\xb)\big]  + (\Ab\bar{\xb}^{k + 1}  - \bb)^T\yb - (\beta_{k+1}/2)\norm{\yb}_2^2\nonumber\\
&+ (\hat{\yb}^k)^T\Ab\big[(1-\tau_k)\bar{\xb}^k + \tau_k\xb - \hat{\xb}^k\big] -  (1/(2\beta_{k+1}))\norm{\Ab\hat{\xb}^k - \bb}_2^2 - (1-\tau_k)\gamma_kd_b(\Sb\xb,\Sb\bar{\xb}^{\star}_k)\nonumber\\
& -  (1/(2\beta_{k\!+\!1}))\left[(1\!-\!\tau_k)\norm{\Ab(\bar{\xb}^k \!-\! \hat{\xb}^k)}_2^2 \!-\! (1 \!-\!\tau_k)\tau_k\norm{\Ab\bar{\xb}^k \!-\! \bb}^2_2 \!+\! \tau_k\norm{\Ab\hat{\xb}^k \!-\! \bb}_2^2\right].%\nonumber
\end{align}

\noindent\textit{\underline{Step 4:} Refining the upper bound of $H_{k+1}$.}
Let $\ub := (1-\tau_k)\bar{\xb}^k + \tau_k\xb \in \Xc$ and
\begin{equation}\label{eq:T3_term}
\mathcal{T}_{[3]} := (1/(2\beta_{k + 1}))\big[(1 - \tau_k)\norm{\Ab(\bar{\xb}^k - \hat{\xb}^k)}_2^2  - (1 - \tau_k)\tau_k\norm{\Ab\bar{\xb}^k - \bb}^2_2 +\tau_k\norm{\Ab\hat{\xb}^k - \bb}_2^2\big].
\end{equation}
First, by the convexity of $f$ we have $f(\ub) \leq (1-\tau_k)f(\bar{\xb}) + \tau_kf(\xb)$.
Second, from the first line of $(\mathrm{2P1D})$ we have $\ub - \hat{\xb}^k = \tau_k(\xb -\bar{\xb}^{\star}_k)$.
Third, by the strong convexity of $d_b$ and the condition \eqref{eq:pd_condition1}, 
we can estimate
\begin{align*}
(1-\tau_k)\gamma_kd_b(\Sb\xb, \Sb\bar{\xb}^{\star}_k) &\geq \bar{L}^g\beta_{k+1}^{-1}\tau_k^2d_b(\Sb\ub, \Sb\bar{\xb}^{\star}_k) \geq   (\bar{L}^g/2)\beta_{k+1}^{-1}\tau_k^2\Vert \Sb(\xb- \bar{\xb}^{\star}_k\Vert_2^2 \geq (\bar{L}^g/2)\beta_{k+1}^{-1}\norm{\Sb(\ub - \hat{\xb}^k)}_2^2.
 \end{align*}
Finally, substituting these expressions into \eqref{eq:lm41_proof9} we obtain
\begin{align}\label{eq:lm41_proof11}
H_{k+1}(\wb) & \leq (1-\tau_k)G_k(\bar{\wb}^k) + f(\bar{\xb}^{k+1}) - f(\ub) +  (\Ab\bar{\xb}^{k + 1}  - \bb)^T\yb - (\beta_{k+1}/2)\norm{\yb}_2^2\nonumber\\
& - (\hat{\yb}^k)^T\Ab(\ub - \hat{\xb}^k) - (\bar{L}^g/2)\beta_{k+1}^{-1}\norm{\Sb(\ub - \hat{\xb}^k)}_2^2 - (1/(2\beta_{k+1}))\norm{\Ab\hat{\xb}^k - \bb}_2^2  - \mathcal{T}_{[3]}.
\end{align}

\noindent\textit{\underline{Step 5:} Final touches on the upper bound of $H_{k+1}$.}
By the third line of \eqref{eq:pd_scheme_2p}, we have $\bar{\xb}^{k+1} := \mathrm{prox}_{\Sb f}(\hat{\xb}^k, \hat{\yb}^k;\beta_{k+1})$.
If we define $\mathcal{H}_{\beta_{k+1}}(\ub) := f(\ub) + (\hat{\yb}^k)^T\Ab(\ub - \hat{\xb}^k) + (\bar{L}^g/(2\beta_{k+1}))\norm{\Sb(\ub - \hat{\xb}^k)}_2^2$, then, by \eqref{eq:prox_oper2}, we have
\begin{align}\label{eq:lm41_proof12}
\mathcal{H}_{\beta_{k+1}}(\ub) \geq \mathcal{H}_{\beta_{k+1}}(\bar{\xb}^{k+1}), ~~\forall \ub\in\Xc.
\end{align}
On the other hand, since $\displaystyle\max_{\yb\in\R^m}\big\{(\Ab\bar{\xb}^{k+1} - \bb)^T\yb - (\beta_{k+1}/2)\norm{\yb}_2^2\big\} = (1/(2\beta_{k+1}))\norm{\Ab\bar{\xb}^{k+1} - \bb}_2^2$, one has
\begin{align}\label{eq:lm41_proof13}
(\Ab\bar{\xb}^{k+1} - \bb)^T\yb - (\beta_{k+1}/2)\norm{\yb}_2^2 \leq (1/(2\beta_{k+1}))\norm{\Ab\bar{\xb}^{k+1} - \bb}_2^2, ~~\forall \yb\in\R^m.
\end{align}
Substituting \eqref{eq:lm41_proof13} and \eqref{eq:lm41_proof12} into \eqref{eq:lm41_proof11} we get
\begin{align}\label{eq:lm41_proof14}
H_{k+1}(\wb) & \leq (1-\tau_k)G_k(\bar{\wb}^k) + f(\bar{\xb}^{k+1}) - f(\bar{\xb}^{k+1}) - (1/(2\beta_{k+1}))\norm{\Ab\hat{\xb}^k - \bb}_2^2\nonumber\\
& - (\hat{\yb}^k)^T\Ab(\bar{\xb}^{k+1} - \hat{\xb}^k) - (\bar{L}^g/(2\beta_{k+1}))\Vert\Sb(\bar{\xb}^{k+1} - \hat{\xb}^k)\Vert_2^2  +  (1/(2\beta_{k+1}))\norm{\Ab\bar{\xb}^{k+1} - \bb}_2^2 - \mathcal{T}_{[3]}.
\end{align}
By the condition \eqref{eq:condition_dp} with $\hat{\xb} = \hat{\xb}^k$, $\xb = \bar{\xb}^{k+1}$ and $\hat{\yb}^k := \beta_{k+1}^{-1}(\Ab\hat{\xb}^k - \bb)$, we have
\begin{align}\label{eq:lm41_proof15}
(2\beta_{k\!+\!1})^{-1}\norm{\Ab\hat{\xb}^k \!-\! \bb}_2^2 + (\hat{\yb}^k)^T\Ab(\bar{\xb}^{k\!+\!1} \!-\! \hat{\xb}^k)  \!+\! (2\beta_{k\!+\!1})^{-1}\bar{L}^g\norm{\Sb(\bar{\xb}^{k\!+\!1} \!-\! \hat{\xb}^k)}_2^2 \geq  (2\beta_{k\!+\!1})^{-1}\norm{\Ab\bar{\xb}^{k\!+\!1} \!-\! \bb}_2^2.
\end{align}
Substituting this inequality into \eqref{eq:lm41_proof14} we finally get
\begin{align}\label{eq:lm41_proof16}
H_{k+1}(\wb) & \leq (1-\tau_k)G_k(\bar{\wb}^k) - \mathcal{T}_{[3]}.
\end{align}
\noindent\textit{\underline{Step 6:} We simplify $\mathcal{T}_{[3]}$ and prove \eqref{eq:descent_gap}.}
From the definition \eqref{eq:T3_term} of $\mathcal{T}_{[3]}$, we can estimate
\begin{align}\label{eq:lm41_proof17}
\mathcal{T}_{[3]} &:=  (2\beta_{k+1})^{-1}\big[ (1-\tau_k)\norm{\Ab(\bar{\xb}^k - \hat{\xb}^k)}_2^2 - \tau_k(1-\tau_k)\norm{\Ab\bar{\xb}^k - \bb}^2_2  + \tau_k\norm{\Ab\hat{\xb}^k - \bb}^2_2\big]\nonumber\\
& =(2\beta_{k+1})^{-1}\norm{(\Ab\hat{\xb}^k - \bb) - (1-\tau_k)(\Ab\bar{\xb}^k - \bb)}^2_2\nonumber\\
&\overset{\tiny{\eqref{eq:pd_scheme_2p}\textrm{(line 1)}}}{=} (2\beta_{k+1})^{-1}\tau_k^2\norm{\Ab\bar{\xb}^{\star}_k - \bb}_2^2.
\end{align}
Substituting \eqref{eq:lm41_proof17} into \eqref{eq:lm41_proof16} and taking the maximization over $\Wc$ we obtain 
\begin{align*}
G_{k+1}(\bar{\wb}^{k+1}) = \max_{\wb\in\mathcal{W}}H_{k+1}(\wb)  \leq (1-\tau_k)G_k(\bar{\wb}^k) - (\tau_k^2/(2\beta_{k+1}))\norm{\Ab\bar{\xb}^{\star}_k - \bb}_2^2,
\end{align*}
which is indeed \eqref{eq:descent_gap}.
\Eproof
% End of the proof.

% Lemma 4.2.
\subsubsection{The proof of Lemma \ref{le:maintain_excessive_gap2}: The condition for selecting step-size $\tau_k$ in \eqref{eq:pd_scheme_2d}.} 
Let us denote by $\bar{\yb}_k^{\star} := \yb^{\star}_{\beta_k}(\bar{\xb}^k)$, $\bar{\xb}^{\star}_k := \xb^{\star}_{\gamma_{k+1}}(\bar{\yb}^k)$ and $\hat{\xb}_k^{\star} := \xb^{\star}_{\gamma_{k+1}}(\hat{\yb}^k)$. 

\noindent\textit{\underline{Step 1:} Estimate $G_k$.}
By using $H_k$ as in the proof of Lemma \ref{le:maintain_excessive_gap1}[\eqref{eq:lm41_proof2}], we have
\begin{align}\label{eq:lm42_proof42_1a}
G_k(\bar{\wb}^k) &= \max_{\wb\in \Wc} \set{ f(\bar{\xb}^k) - f(\xb) + (\Ab\bar{\xb}^k - \bb)^T\yb - (\Ab\xb - \bb)^T\bar{\yb}^k - d_{k}(\wb) } \nonumber\\
&\geq f(\bar{\xb}^k) + \max_{\xb\in\Xc}\big\{ - f(\xb)  - (\Ab\xb - \bb)^T\bar{\yb}^k - \gamma_kd_b(\Sb\xb,\xb_c)  \big\}  + \max_{\yb\in\R^m}\big\{(\Ab\bar{\xb}^k - \bb)^T\yb - \frac{\beta_k}{2}\norm{\yb}^2_2\big\}.
\end{align}
Now, since $\sb^T\yb - (\beta/2)\norm{\yb}_2^2 = (1/(2\beta))\norm{\sb}_2^2 - (\beta/2)\Vert\yb - (1/\beta)\sb\Vert_2^2$ for all $\yb,\sb\in\R^m$, we have
\begin{align*}
(\Ab\bar{\xb}^k \!-\! \bb)^T\yb - (\beta_k/2)\norm{\yb}_2^2 + (\beta/2)\Vert \yb \!-\! \bar{\yb}^{\star}_k\Vert_2^2 \leq \max_{\yb\in\R^m}\set{(\Ab\bar{\xb}^k - \bb)^T\!\yb - (\beta_k/2)\norm{\yb}_2^2}.
\end{align*}
Substituting this estimate into \eqref{eq:lm42_proof42_1a} we get
\begin{align}\label{eq:lm42_proof42_1}
G_k(\bar{\wb}^k) &\geq f(\bar{\xb}^k) + \max_{\xb\in\Xc}\big\{ - f(\xb)  - (\Ab\xb - \bb)^T\bar{\yb}^k - \gamma_kd_b(\Sb\xb,\xb_c)  \big\} \nonumber\\
&  +(\Ab\bar{\xb}^k \!-\! \bb)^T\yb - (\beta_k/2)\norm{\yb}_2^2 + (\beta/2)\Vert \yb \!-\! \bar{\yb}^{\star}_k\Vert_2^2.
\end{align}

\noindent\textit{\underline{Step 2:} Properties of $g_{\gamma}$.}
Let $\varphi_{\gamma}(\yb) := \displaystyle\max_{\xb\in\Xc} \big\{ - f(\xb)  - (\Ab\xb - \bb)^T\bar{\yb}^k - \gamma d_b(\Sb\xb,\Sb\xb_c) \big\}$. It is clear that $\varphi_{\gamma_k}(\yb) \equiv -g_{\gamma_k}(\yb)$, which is convex and smooth. Hence, by Definition \ref{def:smooth_func}, we have
\begin{equation}\label{eq:lm42_proof42_2}
\left\{\begin{array}{ll}
\varphi_{\gamma}(\yb) &{\!\!\!\!\!\!}\geq \varphi_{\gamma}(\hat{\yb}) + \nabla\varphi_{\gamma}(\hat{\yb})^T(\yb - \hat{\yb}),\\
\varphi_{\gamma}(\yb) &{\!\!\!\!\!\!}\leq \varphi_{\gamma}(\hat{\yb}) + \nabla\varphi_{\gamma}(\hat{\yb})^T(\yb - \hat{\yb}) + (L^g_{\gamma}/2)\norm{\yb-\hat{\yb}}_2^2, ~~\forall \yb,\hat{\yb}\in\R^m,\\
\varphi_{\bar{\gamma}}(\yb) &{\!\!\!\!\!\!}\geq \varphi_{\gamma}(\yb) + (\gamma - \bar{\gamma})d_b(\Sb\xb^{\star}_{\gamma}(\yb),\Sb\xb_c), ~~\forall~\bar{\gamma},\gamma > 0.
\end{array}\right.
\end{equation}
Here the first inequality follows from the convexity of $\varphi_{\gamma}$, while the second follows from the Lipschitz continuity of $\nabla{\varphi}_{\gamma}$. We prove the third inequality.
The function $s(\xb, \gamma) := -f(\xb) - \yb^T(\Ab\xb - \bb) - \gamma d_b(\Sb\xb, \Sb\xb_c)$ is concave with respect to $\xb$ and linear with respect to $\gamma$. It is clear that $\varphi_{\gamma}(\yb) = \max_{\xb\in\Xc}\{s(\xb,\gamma)\}$, which is convex with respect to $\gamma$ \cite{Boyd2004}. Moreover, its derivative with respect to $\gamma$ is given by $-d_b(\Sb\xb^{\star}_{\gamma}(\yb), \Sb\xb_c) \geq 0$. This function is nonincreasing, which leads to the first inequality of \eqref{eq:lm42_proof42_2}.

\noindent\textit{\underline{Step 3:} A refinement of $G_k$.}
By the definition of $\hat{\xb}^{\star}_k$ and $\nabla{\varphi}_{\gamma_{k+1}}(\hat{\yb}^k) = \bb - \Ab\hat{\xb}^{\star}_k$, we can express
\begin{align}\label{eq:lm42_proof42_3}
-f(\hat{\xb}^{\star}_k) - (\Ab\hat{\xb}^{\star}_k - \bb)^T\yb &= -f(\hat{\xb}^{\star}_k) - (\Ab\hat{\xb}^{\star}_k - \bb)^T\hat{\yb}^k - (\Ab\hat{\xb}^{\star}_k - \bb)^T(\yb - \hat{\yb}^k) \nonumber\\
& = \varphi_{\gamma_{k\!+\!1}}(\hat{\yb}^k) + \nabla{\varphi_{\gamma_{k+1}}}(\hat{\yb}^k)^T(\yb - \hat{\yb}^k) + \gamma_{k+1}d_b(\Sb\hat{\xb}^{\star}_k,\Sb\xb_c).
\end{align}
Multiplying \eqref{eq:lm42_proof42_1} by $1-\tau_k$ and then using the definition of $\varphi_{\gamma}$, we have
\begin{align}\label{eq:lm42_proof42_4a}
(1-\tau_k)G_k(\bar{\wb}^k) &= (1-\tau_k)\Big[ \varphi_{\gamma_k}(\bar{\yb}^k) + f(\bar{\xb}^k) + (\Ab\bar{\xb}^k - \bb)^T\yb   + \frac{\beta_k}{2}\norm{\yb - \bar{\yb}^{\star}_k}^2_2 - \frac{\beta_k}{2}\norm{\yb}_2^2\Big].
\end{align}
Using the third inequality of \eqref{eq:lm42_proof42_3} with $\bar{\gamma} = \gamma_k$ and $\gamma = \gamma_{k+1} = (1-c_k\tau_k)\gamma_k$ and \eqref{eq:lm42_proof42_3} into \eqref{eq:lm42_proof42_4a} we obtain
\begin{align}\label{eq:lm42_proof42_4}
(1\!-\!\tau_k)G_k(\bar{\wb}^k) &\!\geq\! (1 \!-\!\tau_k)\big[\varphi_{\gamma_{k\!+\!1}}(\bar{\yb}^k) \!+\! f(\bar{\xb}^k)  \!+\! (\Ab\bar{\xb}^k \!-\! \bb)^T\yb  \!+\! \frac{\beta_k}{2}\norm{\yb \!-\! \bar{\yb}^{\star}_k}^2_2 \!-\! \frac{\beta_k}{2}\norm{\yb}_2^2 - \tau_kc_k\gamma_{k}d_b(\Sb\bar{\xb}_k^{\star},\Sb\xb_c)\big] \nonumber\\
& +\! \tau_k\Big[\! f(\hat{\xb}^{\star}_k) \!+\! (\Ab\hat{\xb}^{\star}_k \!-\! \bb)^T\!\yb +\! \varphi_{\gamma_{k\!+\!1}}(\hat{\yb}^k) \!+\!\nabla\varphi_{\gamma_{k\!+\!1}}(\hat{\yb}^k)^T\!(\yb \!-\! \hat{\yb}^k) \!+\! \gamma_{k\!+\!1}d_b(\Sb\hat{\xb}^{\star}_k,\!\Sb\xb_c)\!\Big].
\end{align}
Now, using the first line of \eqref{eq:lm42_proof42_3}, we have $\varphi_{\gamma_{k+1}}(\bar{\yb}^k) \geq \varphi_{\gamma_{k+1}}(\hat{\yb}^k) + \nabla{\varphi_{\gamma_{k+1}}}(\hat{\yb}^k)^T(\bar{\yb}^k - \hat{\yb}^k)$.
On the other hand, by the convexity of $f$ and the second line of \eqref{eq:pd_scheme_2d}, we easily get $f(\bar{\xb}^{k+1}) = f((1-\tau_k)\bar{\xb}^k + \tau_k\hat{\xb}^{\star}_k) \leq (1-\tau_k)f(\bar{\xb}^k) + \tau_kf(\hat{\xb}^{\star}_k)$. Using these inequalities and $\bar{\xb}^{k+1} = (1-\tau_k)\bar{\xb}^k + \tau_k\hat{\xb}_k^{\star}$ into \eqref{eq:lm42_proof42_4} we can further estimate
\begin{align}\label{eq:lm42_proof42_4c}
(1-\tau_k)G_k(\bar{\wb}^k)  &\geq \varphi_{\gamma_{k+1}}(\hat{\yb}^k) + \nabla\varphi_{\gamma_{k\!+\!1}}(\hat{\yb}^k)^T\big[ (1-\tau_k)(\bar{\yb}^k + \tau_k\yb - \hat{\yb}^k\big] + f(\bar{\xb}^{k+1})  + (\Ab\bar{\xb}^{k+1} - \bb)^T\yb \nonumber\\
& + \frac{(1-\tau_k)\beta_k}{2}\big[\norm{\yb - \bar{\yb}^{\star}_k}^2_2 -  \norm{\yb}_2^2\big] + \tau_k\gamma_{k+1}d_b(\Sb\hat{\xb}^k_{\star}, \Sb\xb_c) - (1\!-\!\tau_k)\tau_kc_k\gamma_{k}d_b(\Sb\bar{\xb}^{\star}_k,\Sb\xb_c). 
\end{align}
\noindent\textit{\underline{Step 4:} We prove \eqref{eq:descent_gap}.}
Let $\vb := (1-\tau_k)\bar{\yb}^k + \tau_k\yb\in\R^m$. Using the first line of \eqref{eq:pd_scheme_2d}, we can show that $\vb - \hat{\yb}^k := \tau_k(\yb - \bar{\yb}^{\star}_k)$. 
Substituting $\vb$ into \eqref{eq:lm42_proof42_4c} and taking the maximization over $\R^m$, we get
\begin{align}\label{eq:lm42_proof42_4d}
(1-\tau_k)G_k(\bar{\wb}^k) & \geq \max_{\vb\in\R^m}\Big\{\varphi_{\gamma_{k+1}}(\hat{\yb}^k) + \nabla\varphi_{\gamma_{k+1}}(\hat{\yb}^k)^T(\vb - \hat{\yb}^k) + (\beta_{k+1}/\tau_k^2)\norm{\vb - \hat{\yb}^k}^2_2\Big\} \nonumber \\
&+ f(\bar{\xb}^{k+1}) + (\Ab\bar{\xb}^{k+1} - \bb)^T\yb - (\beta_{k+1}/2)\norm{\yb}^2_2 + \mathcal{T}_{[4]},
\end{align}
where $\mathcal{T}_{[4]} := \tau_k(1-\tau_k)\gamma_{k}\big[d_b(\Sb\hat{\xb}^k_{\star}, \Sb\xb_c) - c_kd_b(\Sb\bar{\xb}^{\star}_k, \Sb\xb_c)\big]$.

From the condition $\beta_{k+1}\gamma_{k+1} \geq \bar{L}^g\tau_k^2$ in \eqref{eq:pd_condition2}, we have $\beta_{k+1}\tau_k^{-2} \geq \bar{L}^g\gamma_{k+1}^{-1} = L^g_{\gamma_{k+1}}$. 
Using this inequality, the second line of \eqref{eq:pd_scheme_2d} and the second inequality of \eqref{eq:lm42_proof42_2} with $\gamma = \gamma_{k+1}$, $\hat{\yb} = \hat{\yb}^k$ and $\yb = \bar{\yb}^{k+1}$, we can further refine \eqref{eq:lm42_proof42_4d} as
\begin{align}\label{eq:lm42_proof42_5a}
(1\!-\!\tau_k)G_k(\bar{\wb}^k) &\!\geq\! \varphi_{\gamma_{k+1}}(\bar{\yb}^{k+1}) + f(\bar{\xb}^{k+1}) +  (\Ab\bar{\xb}^{k+1} - \bb)^T\yb   - (\beta_{k+1}/2)\norm{\yb}^2_2  + \mathcal{T}_{[4]}\nonumber\\
&\geq f(\bar{\xb}^{k + 1})  - f(\xb) + (\Ab\bar{\xb}^{k + 1} - \bb)^T\yb  - (\Ab\xb - \bb)^T\bar{\yb}^{k+1}  - d_{k+1}(\wb) +  \mathcal{T}_{[4]}\\
& = H_{k+1}(\wb) + \mathcal{T}_{[4]}.\nonumber
\end{align}
Since the left-hand side of \eqref{eq:lm42_proof42_5a} is constant, by maximizing over $\wb\in\Wc$ the right-hand side of this inequality, we finally get $G_{k+1}(\bar{\wb}^{k+1}) \leq (1-\tau_k)G_k(\bar{\wb}^k) - \mathcal{T}_{[4]}$, which is indeed \eqref{eq:descent_gap}.
\Eproof
% End of the proof.

%% 4.3. Find a starting point.
\subsubsection{The proof of Lemma \ref{le:starting_point}: Finding starting points for Algorithm \ref{alg:pd_alg}.}
From the definition of $H_k$ in the proof of Lemma \ref{le:maintain_excessive_gap1}[\eqref{eq:lm41_proof2}] and the definition of $g_{\gamma}$, we can show that
\begin{align}\label{eq:lm43_proof1}
G_0(\bar{\wb}^0) &= \max_{\wb\in\Wc}\set{f(\bar{\xb}^0) - f(\xb) - (\Ab\xb - \bb)^T\bar{\yb}^0 + (\Ab\bar{\xb}^0 - \bb)^T\yb - d_0(\wb)}\nonumber\\
&=\max_{\yb\in\R^m}\Big\{ f(\bar{\xb}^0) + (\Ab\bar{\xb}^0 - \bb)^T\yb - \frac{\beta_0}{2}\norm{\yb}_2^2 \Big\} - \min_{\xb\in\Xc}\Big\{f(\xb) + (\Ab\xb - \bb)^T\bar{\yb}^0 + \gamma_0d_b(\Sb\xb,\Sb\xb_c)\Big\} \nonumber\\
&=\max_{\yb\in\R^m}\Big\{ f(\bar{\xb}^0) + (\Ab\bar{\xb}^0 - \bb)^T\yb - \frac{\beta_0}{2}\norm{\yb}_2^2  - g_{\gamma_0}(\bar{\yb}^0)\Big\}.
\end{align}
By  the definition of $\bar{\xb}^0$ and  $\yb^c := 0^m$, we have $g_{\gamma_0}(\yb^c) = f(\bar{\xb}^0) + \gamma_0d_b(\Sb\bar{\xb}^0,\Sb\xb_c)$. 
By Definition \ref{def:smooth_func}, $\nabla{g}_{\gamma}(\cdot)$ is $\bar{L}^g/\gamma_0$-Lipschitz continuous, by \cite[Theorem 2.1.5]{Nesterov2004},  we have $g_{\gamma_0}(\bar{\yb}^0) \geq g_{\gamma_0}(\yb^c) + \nabla{g_{\gamma_0}}(\yb^c)^T(\bar{\yb}^0 - \yb_c) - \frac{\bar{L}^g}{2\gamma_0}\norm{\bar{\yb}^0 - \yb_c}^2_2$. Moreover, $\nabla{g_{\gamma_0}}(\yb^c) = \Ab\bar{\xb}^0 - \bb$ for $\bar{\yb}^0 = (1/\beta_0)\left(\Ab\bar{\xb}^0 - \bb\right)$. 
Hence, $g_{\gamma_0}(\bar{\yb}^0) \geq f(\bar{\xb}^0) + (\Ab\bar{\xb}^0 - \bb)^T\bar{\yb}^0 - \frac{\bar{L}^g}{2\gamma_0\beta_0^2}\Vert\Ab\bar{\xb}^0 - \bb\Vert_2^2$. Using this inequality into \eqref{eq:lm43_proof1}, we can further estimate
\begin{align*}
G_0(\bar{\wb}^0) &\leq \max_{\yb\in\R^m}\Big\{ (\Ab\bar{\xb}^0 - \bb)^T\yb - \frac{\beta_0}{2}\norm{\yb}_2^2 + \frac{\bar{L}^g}{2\gamma_0\beta_0^2}\norm{\Ab\bar{\xb}^0 - \bb}_2^2 - \gamma_0d_b(\Sb\bar{\xb}^0, \Sb\xb_c) \Big\} \nonumber\\
&\leq -\frac{1}{2\beta_0}\Big(2 - \frac{\bar{L}^g}{\beta_0\gamma_0}\Big)\norm{\Ab\bar{\xb}^0 - \bb}_2^2 - \gamma_0d_b(\Sb\bar{\xb}^0, \Sb\xb_c),
\end{align*}
which leads to $G_0(\bar{\wb}^0) \leq - \gamma_0d_b(\Sb\bar{\xb}^0, \Sb\xb_c)$ provided that  $\beta_0\gamma_0 \geq \bar{L}^g$.
The statement \eqref{eq:starting_point2} of Lemma \ref{le:starting_point} can be proved similarly.
\Eproof
% End of the proof.

%% 4.4. Update the parameters.
\subsubsection{The proof of Lemma \ref{le:update}: Update rule for  step-size parameter $\tau_k$.}
For any $c_{k+1} \leq 1$ and $a_k\geq 0$, we have $1 + c_{k+1} + 2a_k \leq 1 + c_{k+1} + \sqrt{4a_k^2 + (1-c_{k+1})^2} \leq 2a_k + 1-c_{k+1} + 1 + c_{k+1} = 2a_k + 2$.
From \eqref{eq:tau_update_rule}, we can easily show that $a_k + (c_{k+1} + 1)/2 \leq a_{k+1} \leq a_k + 1$.
By induction, we can derive from this estimate that
\begin{equation*}
a_0 + k/2 + (1/2)\sum_{i=1}^{k}c_i \leq a_{k} \leq a_0 + k.
\end{equation*} 
On the other hand, from \eqref{eq:tau_update_rule} we have $a_0 := \big(\sqrt{(1+c_0)^2 + 4(1-c_0)} + 1 + c_0\big)/2$.
Combining two last expressions and $s_k := \sum_{i=1}^{k}c_i$, we obtain \eqref{eq:tau_rate}.
The estimate \eqref{eq:beta_gamma_rate} follows from the relation $\beta_{k+1}\gamma_{k+1} = \bar{L}^g\tau_k^2 = \bar{L}^ga_k^{-2}$ and  \eqref{eq:tau_rate}. 

Now, let us consider the case $c_k = 0$ for all $k\geq 0$. Then, the update rule for $\gamma_k$ becomes $\gamma_{k+1} := \gamma_k = \gamma_0 = \bar{L}^g/\beta_0$ for all $k\geq 0$. 
Moreover, we have $a_0 = (1 + \sqrt{5})/2$. Then the first line of \eqref{eq:beta_rate1} follows directly from \eqref{eq:beta_gamma_rate} and $1 < a_0 < 2$.

If $c_k = 1$ for all $k\geq 0$ then $a_0 = 2$ and $\tau_0 = 0.5$. Moreover, we have $(1-\tau_{k+1})^2\tau_k^2 = \tau_{k+1}^2$, which leads to $(1-\tau_{k+1}) = \tau_{k+1}/\tau_k$. Therefore, $\beta_{k+1} = \beta_0\prod_{i=0}^k(1-\tau_i) = \beta_0(1-\tau_0)\prod_{i=1}^k\frac{\tau_i}{\tau_{i-1}} = \beta_0(1-\tau_0)\frac{\tau_k}{\tau_0} = \beta_0a_k^{-1}$. 
Moreover, from \eqref{eq:tau_rate} we have $k+ 2 = k + a_0 \leq a_k \leq k + a_0 = k + 2$.
Combining the last inequality and this equality we obtain the second line of \eqref{eq:beta_rate1}.
\Eproof
%% End of the proof.

%% A.6. The proof of Theorem 5.1.
\subsubsection{The full-proof of Theorem \ref{th:convergence_A1}.}
Under Assumption A.\ref{as:A1}, by the well-known properties of augmented Lagrangian function $\mathcal{L}_{\gamma}$, see, e.g. \cite{Bertsekas1996d}, we have 
\begin{equation*}
\mathcal{L}_{\gamma}(\xb^{\star}, \yb) \leq \mathcal{L}_{\gamma}(\xb^{\star},\yb^{\star}) \equiv \mathcal{L}(\xb^{\star},\yb^{\star}) = f^{\star} = g^{\star} \leq \mathcal{L}_{\gamma}(\xb,\yb^{\star})
\end{equation*}
 for all $\xb\in\Xc, \yb\in\mathbf{R}^m$, $(\xb^{\star}, \yb^{\star})\in\Wc^{\star}$ and $\gamma > 0$.
This expression leads to
\begin{align*}
\tilde{g}_{\gamma}(\yb) \leq   f(\xb) + (\Ab\xb - \bb)^T\yb^{\star} + (\gamma/2)\norm{\Ab\xb - \bb}_2^2 \leq f(\xb) + \norm{\yb^{\star}}_2\norm{\Ab\xb - \bb}_2 + (\gamma/2)\norm{\Ab\xb - \bb}_2^2.
\end{align*}
Hence, for any $\yb^{\star}\in\Yc^{\star}$, we obtain
\begin{align}\label{eq:th51_proof1}
f(\xb) - \tilde{g}_{\gamma}(\yb) \geq f(\xb) - f^{\star} \geq - \norm{\yb^{\star}}_2\norm{\Ab\xb - \bb}_2 - (\gamma/2)\norm{\Ab\xb - \bb}_2^2, ~~\forall  \xb\in\Xc, \yb\in\dom{g_{\gamma}}.
\end{align}
Let $t := \norm{\Ab\bar{\xb}^k - \bb}_2$.  
By combining \eqref{eq:th51_proof1} and \eqref{eq:main_bound_gap} we obtain $\frac{(1-\gamma_k\beta_k)}{\beta_k}t^2 - 2\norm{\yb^{\star}}_2t - 2(\omega_kG_0(\bar{\wb}^0) - \Psi_k) \leq 0$. 
Since $\gamma_k\beta_k \leq \bar{L}^g\tau_{k-1}^2 < \bar{L}^g \equiv 1$, we can show that 
\begin{equation}\label{eq:th51_proof2}
\Vert\Ab\bar{\xb}^k - \bb\Vert_2 \leq \Big(\frac{\beta_k}{1-\beta_k\gamma_k}\Big)\Big[\norm{\yb^{\star}}_2 + \Big(\norm{\yb^{\star}}_2^2 + \frac{2(\omega_kG_0(\bar{\wb}^0) - \Psi_k)(1-\beta_k\gamma_k)}{\beta_k}\Big)^{1/2}\Big].
\end{equation}
To prove \eqref{eq:convergence_A1}, we note that by setting $c_k := 0$ for all $k\geq 0$ in Lemma \ref{le:update}, we can derive $\frac{\beta_k}{1 - \gamma_k\beta_k} \leq \frac{4\sqrt{\bar{L}^g}}{k^2 - 4} \leq \frac{4}{(k+1)^2}$ for $k \geq 0$. In addition, $\omega_kG_0(\bar{\wb}^0) - \Psi_k \leq 0$ due to Lemma \ref{le:starting_point}. 
Using these estimates into \eqref{eq:th51_proof2}, we obtain $\Vert\Ab\bar{\xb}^k - \bb\Vert_2 \leq \frac{8D^{\star}_{\mathcal{Y}}}{(k+1)^2}$, which is the first inequality of \eqref{eq:convergence_A1}.

From \eqref{eq:main_bound_gap} and \eqref{eq:th51_proof1}  we have $f(\bar{\xb}^k) - f^{\star} \leq f(\xb) - \tilde{g}_{\gamma_k}(\bar{\yb}^k) \leq 0$. This inequality and \eqref{eq:th51_proof1} implies the second inequality of \eqref{eq:convergence_A1}.

Next, we prove \eqref{eq:convergence_A2a}. By Lemma \ref{le:excessive_gap_aug_Lag_func} we have $\Vert\Ab\bar{\xb}^k - \bb\Vert_2 \leq \beta_k\norm{\yb^{\star}} + \sqrt{\beta_k^2\norm{\yb^{\star}}^2 + 2\beta_k\gamma_kD_k} \leq 2\beta_kD^{\star}_{\mathcal{Y}} + \sqrt{2\gamma_k\beta_kD_{\Xc}^{\mathbb{I}}}$, where $\yb^{\star}$ is one minimum norm element of $\mathcal{Y}^{\star}$.  
By Lemma \ref{le:update}, we have $\beta_k\gamma_k  = \frac{\bar{L}^g}{(k+1)^2}$ and $\beta_k = \frac{\sqrt{\bar{L}^g}}{k+1}$. Combines these equalities, we obtain the first inequality in \eqref{eq:convergence_A2a}. The second inequality of \eqref{eq:convergence_A2a} follows from Lemma \ref{le:excessive_gap_aug_Lag_func}  and $\beta_k = \frac{\sqrt{\bar{L}^g}}{k+1}$.

To prove \eqref{eq:convergence_A2b}, we first see from Lemma \ref{le:starting_point} and Theorem \ref{le:maintain_excessive_gap2} that the sequence $\set{(\bar{\xb}^k, \bar{\yb}^k)}$ generated by Algorithm \ref{alg:pd_alg} maintains the condition \eqref{eq:descent_gap}. By Lemma \ref{le:excessive_gap_aug_Lag_func} and Lemma \ref{le:update} we have
\begin{equation*}
\Vert\Ab\bar{\xb}^k - \bb\Vert_2 \leq \frac{8\bar{L}_g}{\gamma_0(k+1)^2}\norm{\yb^{\star}} + \frac{\sqrt{8\bar{L}_gD_k}}{(k+1)}.
\end{equation*}
By the definition of $D_{\mathcal{Y}^{\star}}$ and the choice of $\gamma_0$, we obtain from this inequality the first estimate of \eqref{eq:convergence_A2b}.
The second estimate of \eqref{eq:convergence_A2b} immediately follows from \eqref{eq:main_bound_gap_3} and the choice of $\gamma_0$.
\Eproof
%% End of the proof.

% Proof of Corollary 2.1.
\subsection{The proof of Corollary \ref{co:strong_convex_convergence}: Strong convexity case.}
For simplicity of presentation, we divide this proof into few steps.

\noindent\textit{\underline{Step 1:} The proof of Corollary \ref{co:strong_convex_convergence} for the $(\mathrm{1P2D}_s)$ scheme. }
The proof of the two first estimates in Corollary \ref{co:strong_convex_convergence} for $(\mathrm{1P2D}_s)$ can be done similarly to \cite[Theorem 4]{TranDinh2012a}, where we can  show that $-\frac{4L^g_f}{(k+2)^2}(D^{\star}_{\mathcal{Y}})^2 \leq f(\bar{\xb}^k) - g(\bar{\yb}^k) \leq 0$ and $\norm{\Ab\bar{\xb}^k - \bb}_2 \leq \frac{4L^g_f}{(k+2)^2}D^{\star}_{\mathcal{Y}}$. 
However, we have $-\norm{\yb^{\star}}\norm{\Ab\xb - \bb}_2 \leq f(\xb) - f^{\star} \leq f(\xb) - g(\yb)$ for $\xb\in\Xc$, $\yb\in\R^m$ and $\yb^{\star}\in\Yc^{\star}$ due to \eqref{eq:pd_lower_bound}. 
The first inequality implies the second inequality of Corollary \ref{co:strong_convex_convergence}.

\noindent\textit{\underline{Step 2:} The proof of Corollary \ref{co:strong_convex_convergence} for the $(\mathrm{2P1D}_s)$ scheme.}
Next, we prove the first two estimates in Corollary \ref{co:strong_convex_convergence} for the scheme $(\mathrm{2P1D}_s)$. 
Let $\hat{\yb}^k := \beta_k^{-1}(\Ab\hat{\xb}^k - \bb)$.
By applying the same argument as the proof of \eqref{eq:lm41_proof11} in Lemma \ref{le:maintain_excessive_gap1} to the scheme $(\mathrm{2P1D}_s)$, we obtain
\begin{align}\label{eq:scvx_proof1}
H_{k+1}(\wb) & \leq (1-\tau_k)G_k(\bar{\wb}^k) + f(\bar{\xb}^{k+1}) - f(\ub) +  (\Ab\bar{\xb}^{k + 1}  - \bb)^T\yb - (\beta_{k}/2)\norm{\yb}_2^2\nonumber\\
& - (\hat{\yb}^k)^T\Ab(\ub - \hat{\xb}^k) - (1/(2\beta_{k}))\norm{\Ab\hat{\xb}^k - \bb}_2^2  - ((1-\tau_k)\sigma_f/2)\norm{\xb - \xb^{*}(\bar{\yb}^k)}_2^2 - \mathcal{T}_{[3]},
\end{align}
where $\mathcal{T}_{[3]}$ is defined by \eqref{eq:T3_term}.

Let us assume that  $\beta_k(1-\tau_k)\sigma_f \geq \norm{\Ab}_2^2\tau_k^2$. By using this relation, $\xb - \xb^{\star}(\bar{\yb}^k) = \tau^{-1}_k(\ub - \hat{\xb}^k)$ and \eqref{eq:lm41_proof13}, we can further modify \eqref{eq:scvx_proof1} as
\begin{align}\label{eq:scvx_proof1b}
H_{k + 1}(\wb)  & \leq (1 - \tau_k)G_k(\bar{\wb}^k) + f(\bar{\xb}^{k+1}) - \big[ f(\ub)  +  (\hat{\yb}^k)^T\Ab(\ub -  \hat{\xb}^k) + (\norm{\Ab}_2^2/(2\beta_k))\Vert\ub - \hat{\xb}^k\Vert_2^2\big]\nonumber\\
& + (1/(2\beta_k))\Vert\Ab\bar{\xb}^{k+1} - \bb\Vert_2^2  -  (1/(2\beta_{k}))\norm{\Ab\hat{\xb}^k - \bb}_2^2 - \mathcal{T}_{[3]}.
\end{align}
Using the second line $\bar{\xb}^{k+1} = \mathrm{prox}_{\mathbb{I}f}(\hat{\xb}^k,\hat{\yb}^k; \beta_k)$ of $(\mathrm{2P1D}_s)$, we have 
\begin{equation*}
f(\ub)  +  (\hat{\yb}^k)^T\Ab(\ub -  \hat{\xb}^k) + \frac{\norm{\Ab}_2^2}{2\beta_k}\Vert\ub - \hat{\xb}^k\Vert_2^2 \geq f(\bar{\xb}^{k+1}) + (\hat{\yb}^k)^T\Ab(\bar{\xb}^{k+1} -  \hat{\xb}^k) + \frac{\norm{\Ab}_2^2}{2\beta_k}\Vert\bar{\xb}^{k+1} - \hat{\xb}^k\Vert_2^2.
\end{equation*}
Substituting this inequality into \eqref{eq:scvx_proof1b} we reach
\begin{align}\label{eq:scvx_proof1c}
H_{k + 1}(\wb)  & \leq (1 \!-\! \tau_k)G_k(\bar{\wb}^k) \!-\!  (1/(2\beta_{k}))\norm{\Ab\hat{\xb}^k \!-\! \bb}_2^2  \!-\! (\hat{\yb}^k)^T\Ab(\bar{\xb}^{k\!+\!1} \!\!-\!  \hat{\xb}^k) \!-\! (\norm{\Ab}_2^2/(2\beta_k))\Vert\bar{\xb}^{k\!+\!1} \!-\! \hat{\xb}^k\Vert_2^2\nonumber\\
& + (1/(2\beta_k))\Vert\Ab\bar{\xb}^{k+1} - \bb\Vert_2^2  - \mathcal{T}_{[3]}.
\end{align}
Now, we use the expression \eqref{eq:lm41_proof6} for $\bar{\xb}^k := \bar{\xb}^{k+1}$, we can estimate
\begin{align*}
(1/(2\beta_k)\norm{\Ab\bar{\xb}^{k+1} - \bb}_2^2 &\leq  (1/(2\beta_k))\norm{\Ab\hat{\xb}^k-\bb}_2^2 + (\hat{\yb}^k)^T\Ab(\bar{\xb}^{k+1} - \hat{\xb}^k) + (\norm{\Ab}^2_2/(2\beta_k))\norm{\bar{\xb}^{k+1} - \hat{\xb}^k}^2_2.
\end{align*}
Substituting this inequality into \eqref{eq:scvx_proof1c}, we finally obtain 
\begin{align*}
H_{k + 1}(\wb) &\leq (1 - \tau_k)G_k(\bar{\wb}^k) - \mathcal{T}_{[3]}.
\end{align*}
Since $\mathcal{T}_{[3]} = (2\beta_{k+1})^{-1}\tau_k^2\Vert\Ab\xb^{\star}(\bar{\yb}) - \bb\Vert_2^2$ due to \eqref{eq:lm41_proof17}, by maximizing the last inequality over $\wb\in\Wc$, we obtain  $G_{k+1}(\bar{\wb}^{k+1}) \leq (1 -\tau_k)G_k(\bar{\wb}^k) - (2\beta_{k+1})^{-1}\tau_k^2\Vert\Ab\xb^{\star}(\bar{\yb}) - \bb\Vert_2^2$. 
This inequality shows that the condition \eqref{eq:descent_gap} satisfies with $\psi_k := (2\beta_{k+1})^{-1}\tau_k^2\Vert\Ab\xb^{\star}(\bar{\yb}) - \bb\Vert_2^2 \geq 0$.

To complete the proof, we derive the condition on updating $\tau_k$ from $\frac{(1-\tau_k)\sigma_f}{\tau_k^2} \geq \frac{\norm{\Ab}_2^2}{\beta_k}$. 
Indeed, since $\beta_{k+1} = (1-\tau_k)\beta_k$, we have $\frac{(1-\tau_{k+1})\sigma_f}{\tau_{k+1}^2} \geq \frac{\norm{\Ab}_2^2}{\beta_{k+1}}$ by induction. 
Combining the two last conditions with equality, we obtain  $(1-\tau_{k+1})\tau_k^2 = \tau_{k+1}^2$. This relation leads to $\tau_{k+1} = \tau_k\big(\sqrt{\tau_k^2 + 4} - \tau_k\big)/2$ as given in Corollary \ref{co:strong_convex_convergence}. 
Now, we use the same argument as the proof of $(\mathrm{1P2D}_s)$  to obtain the worst-case bounds in Corollary \ref{co:strong_convex_convergence}.

\vskip0.05cm
\noindent\textit{\underline{Step 3:} The proof for the bound on $\{\bar{\xb}^k\}$ in Corollary \ref{co:strong_convex_convergence}.}
Finally, we prove the last estimate of \eqref{eq:strong_cvx_main_estimate}.
Indeed, by the strong convexity of $f$, we have $f(\bar{\xb}^k) - f^{\star} \geq \xi_f(\xb^{\star})^T(\bar{\xb}^k - \xb^{\star}) + \frac{\sigma_f}{2}\norm{\bar{\xb}^k - \xb^{\star}}_2^2$, where $\xi_f(\xb^{\star}) \in \partial{f}(\xb^{\star})$ is one subgradient of $f$ at $\xb^{\star}$.
On the other hand, since $\xb^{\star}$ is the optimal solution of \eqref{eq:separable_convex}, using the optimality condition of this problem, we have $(\xi_{f}(\xb^{\star}) + \Ab^T\yb^{\star})^T(\xb - \xb^{\star}) \geq 0$ for any $\xb\in\Xc$ and $\yb^{\star}\in\mathcal{Y}^{\star}$ and $\Ab\xb^{\star} = \bb$. Using these expressions, we can show that 
\begin{equation*}
f(\bar{\xb}^k) - f^{\star} \geq \frac{\sigma_f}{2}\Vert\bar{\xb}^k - \xb^{\star}\Vert_2^2 - (\Ab\bar{\xb}^k - \bb)^T\yb^{\star} \geq \frac{\sigma_f}{2}\Vert\bar{\xb}^k - \xb^{\star}\Vert_2^2 - \Vert\yb^{\star}\Vert_2\Vert\Ab\bar{\xb}^k - \bb\Vert_2.
\end{equation*}
This estimate leads to $\Vert\bar{\xb}^k - \xb^{\star}\Vert^2_2 \leq \frac{2}{\sigma_f}[f(\bar{\xb}^k) - f^{\star}] + \frac{2\Vert\yb^{\star}\Vert_2}{\sigma_f}\Vert\Ab\bar{\xb}^k - \bb\Vert_2 \leq \frac{16L^g_f}{\sigma_f(k+2)^2}(D^{\star}_{\mathcal{Y}})^2$,  which is indeed the third estimate in Corollary \ref{co:strong_convex_convergence}.
\Eproof
%% End of the proof.

%%%%%%%%%%%%%%%%%%%%%%%%%%%%%%%%%%%%%%%%%%%%%%%%%%%%%%
% The proof of Lemma 5.1.
\subsection{The proof of Lemma \ref{le:strongly_concave}: Dual function $g_{\gamma}$ is strongly convex}
Let $\xb^{\star}_{p,\gamma}(\yb)$ be the solution of the minimization problem in \eqref{eq:g_gamma^p}. Since this problem is an unconstrained convex minimization, we can write its optimality condition as
\begin{equation}\label{eq:opt_cond_prob_p}
\Ab_p^T\yb + \nabla{f_p}(\xb^{\star}_{p,\gamma}(\yb)) + \gamma\nabla d_p(\xb^{\star}_{p,\gamma}(\yb), \xb_p^c) = 0.
\end{equation} 
Moreover, we have $\nabla{g}_{\gamma}^p(\yb) := \Ab_p\xb^{\star}_{p,\gamma}(\yb)$. 
Since $\nabla{f_p}$ is $L_{f_p}$-Lipschitz gradient and $\nabla d_p(\cdot, \xb_p^c)$ is $1$-Lipschitz continuous,  the function $\psi_p(\cdot) := \nabla{f_p}(\cdot) + \gamma \nabla d_p(\cdot, \xb_p^c)$ is $(L_{f_p} + \gamma)$-Lipschitz continuous.
Using  Baillon-Haddad's theorem \cite[Corollary 18.16]{Bauschke2011}, we obtain that $\psi_p(\cdot)$ is $(L_{f_p} + \gamma)^{-1}$-co-coercive, i.e.,:
\begin{equation}\label{eq:cocoercive}
(\psi_p(\xb_p) - \psi_p(\hat{\xb}_p))^T(\xb_p - \hat{\xb}_p) \geq (L_{f_p} + \gamma)^{-1}\Vert \psi_p(\xb_p) - \psi_p(\hat{\xb}_p)\Vert_2^2, ~~\forall \xb_p, \hat{\xb}_p\in\R^{n_p}.
\end{equation}
Now, let $g_{\gamma}^p$ be defined by \eqref{eq:g_gamma^p}, we estimate the term $\mathcal{A} := (\nabla{g}_{\gamma}^p(\yb) - \nabla{g}_{\gamma}^p(\hat{\yb}))^T(\yb - \hat{\yb})$ as follows:
\begin{align*}
\begin{array}{ll}
(\nabla{g}_{\gamma}^p(\yb) - \nabla{g}_{\gamma}^p(\hat{\yb}))^T(\yb - \hat{\yb}) &=  (\Ab_p\xb^{\star}_{p,\gamma}(\yb) - \Ab_p\xb^{\star}_{p,\gamma}(\hat{\yb}))^T(\yb - \hat{\yb}) \nonumber\\
& =  (\yb - \hat{\yb})^T\Ab_p(\xb^{\star}_{p,\gamma}(\yb) - \xb^{\star}_{p,\gamma}(\hat{\yb})) \nonumber\\
&\overset{\tiny\eqref{eq:opt_cond_prob_p}}{=} - \left(\psi_p(\xb^{\star}_{p,\gamma}(\yb)) - \psi_p(\xb^{\star}_{p,\gamma}(\hat{\yb}))\right)^T(\xb^{\star}_{p,\gamma}(\yb) - \xb^{\star}_{p,\gamma}(\hat{\yb})) \nonumber\\
&\overset{\tiny\eqref{eq:cocoercive}}{\leq} - (L_{f_p} + \gamma)^{-1}\Vert \psi_p(\xb^{\star}_{p,\gamma}(\yb)) - \psi_p(\xb^{\star}_{p,\gamma}(\hat{\yb}) \Vert_2^2 \nonumber\\
&\overset{\tiny\eqref{eq:opt_cond_prob_p}}{=} -(L_{f_p}+\gamma)^{-1}\Vert\Ab_p^T(\yb - \hat{\yb}\Vert_2^2 \nonumber\\
&\leq - (L_{f_p}+\gamma)^{-1}\lambda_{\min}(\Ab_p^T\Ab_p)\Vert\yb - \hat{\yb}\Vert_2^2.
\end{array}
\end{align*}
This inequality shows that $g_{\gamma}^p$ is strongly concave with the parameter $\sigma_{g_{\gamma}^p} := (L_{f_p}+\gamma)^{-1}\lambda_{\min}(\Ab_p^T\Ab_p) > 0$.
Since $g_{\gamma}(\cdot) = \sum_{i=1}^{p-1}g_{\gamma}^i(\cdot) + g_{\gamma}^p(\cdot)$, it is also strongly convex with the same parameter $\sigma_{g_{\gamma}^p} > 0$.
\Eproof
%% End of the proof.

\subsection{The proof of Corollary \ref{co:Lipschitz_convergence}: The Lipschitz gradient case.}
From Lemma \ref{le:strongly_concave}, we note that $\varphi_{\gamma} = -g_{\gamma}$ satisfies $\varphi_{\gamma}(\yb) \geq \varphi(\hat{\yb}) + \nabla{\varphi_{\gamma}}(\hat{\yb})^T(\yb - \hat{\yb}) + (\sigma_g/2)\norm{\yb - \hat{\yb}}_2^2$, where $\sigma_g := (L_{f_p}+\gamma_0)^{-1}\lambda_{\min}(\Ab_p^T\Ab_p) \leq (L_{f_p}+\gamma_k)^{-1}\lambda_{\min}(\Ab_p^T\Ab_p)$ for all $k\geq 0$ due to $\gamma_k \leq \gamma_0$. 
Using this inequality instead of the second inequality of \eqref{eq:lm42_proof42_2} and $\gamma_{k+1} - \gamma_k = -\tau_k\gamma_{k+1}$, we obtain from \eqref{eq:lm42_proof42_5a} that
\begin{align}\label{eq:cor2_proof1}
(1-\tau_k)G_k(\bar{\wb}^k) & \geq H_{k+1}(\wb) + \bar{\mathcal{T}}_{[4]},
\end{align}
where 
\begin{align*}
\bar{\mathcal{T}}_{[4]} &:= (\gamma_{k+1}/2)\big[\tau_k\norm{\Sb(\hat{\xb}_k^{\star} - \xb_c)}_2^2 - (1-\tau_k)\tau_k\norm{\Sb(\bar{\xb}^{\star}_k - \xb_c)}_2^2 + \sigma_g(1-\tau_k)\norm{\Sb(\bar{\xb}^{\star}_k - \hat{\xb}_k^{\star})}_2^2\big] \\
&\geq (\underline{\sigma}_g\gamma_{k+1}/2)\norm{\Sb(\hat{\xb}^{\star}_k - \xb_c) - (1-\tau_k)\Sb(\bar{\xb}_k^{\star} - \xb_c)}_2^2.
\end{align*}
Here $\bar{\xb}^{\star}_k := \xb_{\gamma_{k+1}}^{\star}(\bar{\yb}^k)$ and $\underline{\sigma}_g := \min\set{\sigma_g, 1} > 0$.
We note that $\bar{\mathcal{T}}_{[4]} \geq 0$, taking the maximization both sides in \eqref{eq:cor2_proof1} w.r.t. $\wb\in\Wc$, we obtain $(1-\tau_k)G_k(\bar{\wb}^k) \geq G_{k+1}(\bar{\wb}^{k+1}) + \psi_k$, where $\psi_k := (\underline{\sigma}_g\gamma_{k+1}/2)\norm{\Sb(\hat{\xb}^{\star}_k - \xb_c) - (1-\tau_k)\Sb(\bar{\xb}_k^{\star} - \xb_c)}_2^2 \geq 0$.
Finally, the proof of the estimates \eqref{eq:Lipschitz_main_estimate} in Corollary \ref{co:Lipschitz_convergence} can be done similarly as the proof of Theorem \ref{th:convergence_A1}(c).
\Eproof
%% End of the proof.

%%%%%%%%%%%%%%%%%%%%%%%%%%%%%%%%%%%%%%%%%%%%%%%%
%% A.6. The proof of Theorem 5.1.
\subsection{The proof of Theorem \ref{th:inexact_convergence}: Inexact augmented Lagrangian method}
We divide the prove into few steps as follows.

\noindent\textit{\underline{Step 1:} Approximate smoothed gap function.}
Let us define an approximate gap function $G_{\gamma\beta}^{\delta}$ of the exact smoothed gap function $G_{\gamma\beta}$ in \eqref{eq:smoothed_gap_func} as follows:
\begin{equation}\label{eq:inexact_smoothed_gap_func}
G^{\delta}_{\gamma\beta}(\bar{\wb}) := \delta\textrm{-}\max_{\wb\in\Wc}\set{f(\bar{\xb}) - f(\xb) + F(\wb)^T(\bar{\wb} - \wb) - d_{\gamma\beta}(\wb)},
\end{equation}
where the approximation only involves in $\xb$ in the sense of \eqref{eq:prox_oper2_inexact}, i.e.:
\begin{equation}\label{eq:gap_distance}
G_{\gamma\beta}(\bar{\wb}) \leq G_{\gamma\beta}^{\delta}(\bar{\wb}) + (\gamma/2)\delta^2.
\end{equation}

\noindent\textit{\underline{Step 2:} The first estimate of $G_k$.}
Let $\varphi_{\gamma}$ be defined by \eqref{eq:lm42_proof42_2}, $\hat{\xb}^{\delta}_k := \xb^{\delta}_{\gamma}(\hat{\yb}^k)$, $\varphi_{\gamma}^{\delta}(\yb) := - f(\xb^{\delta}_{\gamma}(\yb)) - (\Ab\xb^{\delta}_{\gamma}(\yb) - \bb)^T\yb - \gamma d_b(\Sb\xb^{\delta}_{\gamma}(\yb),\Sb\xb_c)$ and $\nabla{\varphi_{\gamma}^{\delta}}(\yb) := \bb - \Ab\xb^{\delta}_{\gamma}(\yb)$. 
Then, by \eqref{eq:prox_oper2_inexact} we have
\begin{equation}\label{eq:th63_proof1}
\varphi_{\gamma}(\yb) - \varphi_{\gamma}^{\delta}(\yb) \leq \gamma\delta^2/2~~~\textrm{and}~~~\norm{\nabla{\varphi_{\gamma}^{\delta}}(\yb) - \nabla{\varphi_{\gamma}}(\yb)}_2\leq \delta.
 \end{equation}
 Since $\varphi_{\gamma_k}(\bar{\yb}^k) \geq \varphi_{\gamma_k}(\hat{\yb}^k) + \nabla{\varphi}_{\gamma_k}(\hat{\yb}^k)^T(\bar{\yb}^k - \hat{\yb}^k)$ and $f(\hat{\xb}_k^{\delta}) + (\Ab\hat{\xb}_k^{\delta} - \bb)^T\yb + (\gamma_k/2)\norm{\Ab\hat{\xb}_k^{\delta} - \bb}_2^2 + \varphi_{\gamma_k}^{\delta}(\hat{\yb}^k) + \nabla{\varphi}_{\gamma_k}^{\delta}(
\hat{\yb}^k)^T(\yb - \hat{\yb}^k) = 0$, it follows from \eqref{eq:lm42_proof42_4a} and $\beta_{k+1} = (1-\tau_k)\beta_k$ that
\begin{align*}
(1-\tau_k)G_k(\bar{\wb}^k) &\geq (1-\tau_k)\big[\varphi_{\gamma_k}(\hat{\yb}^k) + \nabla\varphi_{\gamma_k}(\hat{\yb}^k)^T(\bar{\yb}^k - \hat{\yb}^k) + f(\bar{\xb}^k) + (\Ab\bar{\xb}^k - \bb)^T\yb\big]\nonumber\\
&+\tau_k\big[f(\hat{\xb}_k^{\delta}) + (\Ab\hat{\xb}_k^{\delta} - \bb)^T\yb + \varphi_{\gamma_k}^{\delta}(\hat{\yb}^k) + \nabla{\varphi}_{\gamma_k}^{\delta}( \hat{\yb}^k)^T(\yb - \hat{\yb}^k) \big] \\
&+ (\beta_{k+1}/2)\norm{\yb - \bar{\yb}^{\star}_k}_2^2 - (\beta_{k+1}/2)\norm{\yb}_2^2 + (\tau_k\gamma_k/2)\norm{\Ab\hat{\xb}_k^{\delta} - \bb}_2^2.\nonumber
 \end{align*}
 Now, using \eqref{eq:th63_proof1}, the third line $\bar{\xb}^{k+1} = (1-\tau_k)\bar{\xb}^k + \tau_k\tilde{\xb}^{\delta_k}_{\gamma_k}(\hat{\yb}^k)$ of $(\mathrm{i1P2D})$, $\ub := (1-\tau_k)\bar{\yb}^k + \tau_k\yb$ and $\ub - \hat{\yb}^k = \tau_k(\yb - \bar{\yb}^{\star}_k)$ , we can further estimate
\begin{align}\label{eq:th63_proof2b}
(1-\tau_k)G_k(\bar{\wb}^k) & \geq \varphi_{\gamma_k}^{\delta}(\hat{\yb}^k) + \nabla\varphi_{\gamma_k}^{\delta}(\hat{\yb}^k)^T(\ub - \hat{\yb}^k)  +\frac{\beta_{k+1}}{2\tau_k^2}\norm{\ub - \hat{\yb}^k}_2^2 + \frac{\tau_k\gamma_k}{2}\norm{\Ab\hat{\xb}_k^{\delta} - \bb}_2^2 + f(\bar{\xb}^{k+1}) \nonumber\\
&+(1-\tau_k)\big[\nabla{\varphi}_{\gamma_k}( \hat{\yb}^k) - \nabla{\varphi}_{\gamma_k}^{\delta}( \hat{\yb}^k)\big]^T(\bar{\yb}^k - \hat{\yb}^k) + (\Ab\bar{\xb}^{k+1} - \bb)^T\yb - \frac{\beta_{k+1}}{2}\norm{\yb}_2^2. 
\end{align}

\noindent\textit{\underline{Step 3:} The second estimate of $G_k$.}
From the fourth line of $(\mathrm{i1P2D})$ we have $\varphi_{\gamma_k}^{\delta}(\hat{\yb}^k) + \nabla\varphi_{\gamma_k}^{\delta}(\hat{\yb}^k)^T(\ub - \hat{\yb}^k) + \frac{\bar{L}^g}{2\gamma_k}\Vert\ub-\hat{\yb}^k\Vert_2^2 \geq \varphi_{\gamma_k}^{\delta}(\hat{\yb}^k) + \nabla\varphi_{\gamma_k}^{\delta}(\hat{\yb}^k)^T(\bar{\yb}^{k+1} - \hat{\yb}^k) + \frac{\bar{L}^g}{2\gamma_k}\Vert\bar{\yb}^{k+1} -\hat{\yb}^k\Vert_2^2$.
Using this inequality, $\bar{L}^g = 1$ and the condition $\beta_{k+1}\gamma_k   \geq \bar{L}^g\tau_k^2 = \tau_k^2$ we can show that
\begin{align}\label{eq:th63_proof2c}
(1-\tau_k)G_k(\bar{\wb}^k) & \geq \varphi_{\gamma_k}^{\delta}(\hat{\yb}^k) + \nabla\varphi_{\gamma_k}^{\delta}(\hat{\yb}^k)^T(\bar{\yb}^{k+1} - \hat{\yb}^k) +\frac{\bar{L}^g}{2\gamma_k}\norm{\bar{\yb}^{k+1} - \hat{\yb}^k}_2^2 + f(\bar{\xb}^{k+1}) + (\Ab\bar{\xb}^{k+1} - \bb)^T\yb \nonumber\\
& - \frac{\beta_{k+1}}{2}\norm{\yb}_2^2 + \frac{\tau_k\gamma_k}{2}\norm{\Ab\hat{\xb}_k^{\delta} - \bb}_2^2 +(1-\tau_k)\big[\nabla{\varphi}_{\gamma_k}( \hat{\yb}^k) - \nabla{\varphi}_{\gamma_k}^{\delta}( \hat{\yb}^k)\big]^T(\bar{\yb}^k - \hat{\yb}^k). 
\end{align}
By using \eqref{eq:th63_proof1} and the first inequality of \eqref{eq:lm42_proof42_2} we can write
\begin{align}\label{eq:th63_proof2d}
\mathcal{T}_{[4]} &:= \varphi_{\gamma_k}^{\delta}(\hat{\yb}^k) +  \nabla{\varphi}_{\gamma_k}^{\delta}(\hat{\yb}^k)^T(\bar{\yb}^{k + 1} - \hat{\yb}^k) + \frac{\bar{L}^g}{2\gamma_k}\norm{\bar{\yb}^{k + 1} - \hat{\yb}^k}_2^2 \nonumber\\
&\geq \varphi_{\gamma_k}(\hat{\yb}^k) + \nabla{\varphi}_{\gamma_k}(\hat{\yb}^k)^T(\bar{\yb}^{k+1} - \hat{\yb}^k) + \frac{\bar{L}^g}{2\gamma_k}\norm{\bar{\yb}^{k+1} - \hat{\yb}^k}_2^2 + \big[\nabla{\varphi}_{\gamma_k}^{\delta}(\hat{\yb}^k) - \nabla{\varphi}_{\gamma_k}(\hat{\yb}^k)\big]^T(\bar{\yb}^{k+1} - \hat{\yb}^k) - (\gamma_k\delta^2/2) \nonumber\\
&\geq \varphi_{\gamma_k}(\bar{\yb}^{k+1}) + \big[\nabla{\varphi}_{\gamma_k}^{\delta}(\hat{\yb}^k) - \nabla{\varphi}_{\gamma_k}(\hat{\yb}^k)\big]^T(\bar{\yb}^{k+1} - \hat{\yb}^k) - (\gamma_k\delta^2/2). 
\end{align}
Substituting \eqref{eq:th63_proof2d} into \eqref{eq:th63_proof2c} and  then using \eqref{eq:th63_proof1} and the definition of $\varphi_{\gamma}(\cdot)$ we get
\begin{align}\label{eq:th63_proof2e}
(1- \tau_k)G_k(\bar{\wb}^k) & \geq \varphi_{\gamma_k}(\bar{\yb}^{k+1}) + f(\bar{\xb}^{k+1}) + (\Ab\bar{\xb}^{k+1} - \bb)^T\yb - (\beta_{k+1}/2)\norm{\yb}_2^2\nonumber\\
&\!\! - \delta_k\norm{(1\!-\!\tau_k)\bar{\yb}^k \!+\! \tau_k\hat{\yb}^k - \bar{\yb}^{k+1}}_2 + \frac{\tau_k\gamma_k}{2}\Vert\Ab\hat{\xb}_k^{\delta_k} \!-\! \bb\Vert_2^2 \!-\! (\gamma_k\delta^2_{k}/2)\nonumber\\
&\geq f(\bar{\xb}^{k\!+\!1}) \!-\! f(\xb) \!+\! (\Ab\bar{\xb}^{k\!+\!1} \!-\! \bb)^T\yb - (\Ab\xb \!-\! \bb)^T\bar{\yb}^{k\!+\!1} - (\beta_{k\!+\!1}/2)\norm{\yb}_2^2 - (\gamma_k/2)\norm{\Ab\xb - \bb}_2^2 - \mathcal{T}_{[5]}\nonumber\\
&\geq H_{k+1}(\bar{\wb}^{k+1}) - \mathcal{T}_{[5]}, 
\end{align}
provided that $\gamma_{k+1}\geq \gamma_k$, where $\mathcal{T}_{[5]} := \delta_k\norm{(1\!-\!\tau_k)\bar{\yb}^k \!+\! \tau_k\hat{\yb}^k - \bar{\yb}^{k+1}}_2 + (\gamma_k\delta^2_{k})/2 - (\tau_k\gamma_k/2)\Vert\Ab\hat{\xb}_k^{\delta_k} \!-\! \bb\Vert_2^2$.

\noindent\textit{\underline{Step 4:} Simplify $\mathcal{T}_{[5]}$ to obtain \eqref{eq:descent_gap}.}
Using the definition of $\hat{\yb}^k$ and $\bar{\yb}^{k+1}$ we can further estimate $\mathcal{T}_{[5]}$ as 
\begin{align}\label{eq:th63_proof4}
\mathcal{T}_{[5]} := (1\!-\!\tau_k)\tau_k\delta_k\norm{\bar{\yb}^k \!-\! \hat{\yb}^k}_2 \!+\!  (\gamma_k\delta^2_{k})/2 \!+\! \gamma_k\delta_{k}\Vert\Ab\hat{\xb}_k^{\delta_k} \!-\! \bb\Vert_2 \!-\! (\tau_k\gamma_k/2)\Vert\Ab\hat{\xb}_k^{\delta_k} \!-\! \bb\Vert_2^2.
\end{align}
Taking the maximization of \eqref{eq:th63_proof2e} over $\wb\in\Wc$ we finally obtain
\begin{align}\label{eq:th63_proof5}
G_{k+1}(\bar{\wb}^{k+1}) \leq (1 - \tau_k)G_k(\bar{\wb}^k) - \psi_k,
\end{align}
where $\psi_k := (\tau_k\gamma_0/2)\Vert\Ab\hat{\xb}_k^{\delta_k} \!-\! \bb\Vert_2^2 -(1\!-\!\tau_k)\tau_k\delta_k\norm{\bar{\yb}^k \!-\! \bar{\yb}^{\star}_k}_2 - \gamma_0\delta_{k}\Vert\Ab\hat{\xb}_k^{\delta_k} \!-\! \bb\Vert_2 - (\gamma_0\delta^2_{k})/2$.

\noindent\textit{\underline{Step 5:} Prove \eqref{eq:convergence_inexact}.}
We note that $\Vert\Ab\hat{\xb}_k^{\delta_k} \!-\! \bb\Vert_2\leq D_{\Xc}^{\Ab}$ and $\gamma_0 = \bar{L}^g = 1$, which lead to $\psi_k \geq -q_k\delta_k$, where $q_k := (1-\tau_k)\tau_k\norm{\bar{\yb}^k \!-\! \bar{\yb}^{\star}_k}_2 + (D_{\Xc}^{\Ab} + 1)/2$.
In this case \eqref{eq:th63_proof5} leads to $G_{k+1}(\bar{\wb}^{k+1}) \leq (1 - \tau_k)G_k(\bar{\wb}^k) + q_k\delta_k$. Therefore, if we choose $\delta_k$ so that $q_k\delta_k  \leq q_{k-1}\delta_{k-1}$. Then, by induction and $\prod_{i=0}^k(1-\tau_k) \leq \frac{4}{(k+2)^2}$ due to Lemma \ref{le:update}, the last estimate leads to 
\begin{equation}\label{eq:th63_proof6}
G_k(\bar{\wb}^k) \leq \omega_kG_0(\bar{\wb}^0) + q_0\delta_0 + 4\sum_{j=1}^{k-1}\frac{q_j\delta_j}{(j+1)^2} \leq \omega_kG_0(\bar{\wb}^0) + 4q_0\delta_0\zeta(2).
\end{equation}
Here $\zeta(s) := \sum_{j=1}^{\infty}j^{-s}$ is the zeta-function. 
We note that the starting point $\bar{\wb}^0$ is also computed up to the accuracy $\delta_0$, i.e. $G_0(\bar{\wb}^0) \leq (\gamma_0\delta_0^2/2)$ and $\zeta(2) < 1.64494$. 
Plugging these into \eqref{eq:th63_proof6} we have $G_k(\bar{\wb}^k) \leq 7q_0\delta_0$.
Combining this inequality and Lemma \ref{le:excessive_gap_aug_Lag_func} we obtain the second estimate of \eqref{eq:convergence_inexact}.
Finally, by using the bound $G_k(\bar{\wb}^k) \leq 7q_0\delta_0$, it follows from \eqref{eq:th51_proof2} that $\norm{\Ab\bar{\xb}^k - \bb}_2 \leq \frac{4}{(k+1)^2}\left[2D^{\star}_{\mathcal{Y}} + \sqrt{\frac{14q_0\delta_0}{(k+1)^2}}\right]$, which is the first estimate in \eqref{eq:convergence_inexact}.
%\end{proof}
\Eproof

\bibliographystyle{plain}

\end{document}